\documentclass[sn-mathphys]{sn-jnl}

\usepackage{anyfontsize}
\usepackage{amssymb}
\usepackage{cmap}
\usepackage{hyperxmp}
\usepackage{mathtools}
\usepackage{tikz-cd}
\usepackage{program}
\mathtoolsset{showonlyrefs}
\usepackage{enumerate}
\usepackage{todonotes}

\usepackage{lineno}

\jyear{2024}%

\theoremstyle{thmstyleone}%

\theoremstyle{thmstyletwo}%

\theoremstyle{thmstylethree}%

\newcommand{\RR}{\mathbb{R}} 

\DeclareMathOperator{\img}{im} %

\DeclareMathOperator{\rank}{rank} %

\def\cech{\v{C}ech }

\begin{document}

\title[ ]{Topological Data Analysis and Topological Deep Learning Beyond Persistent Homology--A Review}

\author[1]{\fnm{Zhe} \sur{Su}}
\author[2]{\fnm{Xiang} \sur{Liu}}
\author[3]{\fnm{Layal} \sur{Bou Hamdan}}
\author[3]{\fnm{Vasileios} \sur{Maroulas}}
\author[4]{\fnm{Jie} \sur{Wu}}
\author[5]{\fnm{Gunnar} \sur{Carlsson}}
\author[2,6,7]{\fnm{Guo-Wei} \sur{Wei}\textsuperscript{*}}

\affil[1]{\orgdiv{Department of Mathematics and Statistics}, \orgname{Auburn University}, \orgaddress{\state{AL 36849}, \country{USA}}}
\affil[2]{\orgdiv{Department of Mathematics}, \orgname{Michigan State University}, \orgaddress{\state{MI 48824}, \country{USA}}}
\affil[3]{\orgdiv{Department of Mathematics}, \orgname{University of Tennessee, Knoxville}, \orgaddress{\state{TN 37996-1320}, \country{USA}}}
\affil[4]{\orgname{Beijing Institute of Mathematical Sciences and Applications}, 
\orgaddress{\state{Beijing 101408}},
\orgaddress{\country{China}}}
\affil[5]{\orgdiv{Department of Mathematics}, \orgname{Stanford University}, \orgaddress{\state{CA  94305}, \country{USA}}}
\affil[6]{\orgdiv{Department of Biochemistry and Molecular Biology}, \orgname{Michigan State University}, \orgaddress{\state{MI 48824}, \country{USA}}}
\affil[7]{\orgdiv{Department of Electrical and Computer Engineering}, \orgname{Michigan State University}, \orgaddress{\state{MI 48824}, \country{USA}}}

\footnotetext[1]{*Corresponding author. Email: weig@msu.edu}


\abstract{
Topological data analysis (TDA) is a rapidly evolving field in applied mathematics and data science that leverages tools from topology to uncover robust, shape-driven insights in complex datasets. The main workhorse is persistent homology, a technique rooted in algebraic topology. Paired with topological deep learning (TDL) or topological machine learning, persistent homology has achieved tremendous success in a wide variety of applications in science, engineering, medicine, and industry. However, persistent homology has many limitations due to its high-level abstraction, insensitivity to non-topological changes, and reliance on point cloud data. This paper presents a comprehensive review of TDA and TDL beyond persistent homology. It analyzes how persistent topological Laplacians and Dirac operators provide spectral representations to capture both topological invariants and homotopic evolution. Other formulations are presented in terms of sheaf theory, Mayer topology, and interaction topology. For data on differentiable manifolds, techniques rooted in differential topology, such as persistent de Rham cohomology, persistent Hodge Laplacian, and Hodge decomposition, are reviewed. For one-dimensional (1D) curves embedded in 3-space, approaches from geometric topology are discussed, including multiscale Gauss-link integrals, persistent Jones polynomials, and persistent Khovanov homology. This paper further discusses the appropriate selection of topological tools for different input data, such as point clouds, sequential data, data on manifolds, curves embedded in 3-space, and data with additional non-geometric information. A review is also given of various topological representations, software packages, and machine learning vectorizations. Finally, this review ends with concluding remarks. 
}



\maketitle

\newpage

{\setcounter{tocdepth}{4} \tableofcontents}

\newpage

\section{Introduction}

Topological data analysis (TDA), rooted in topological techniques, is an emerging approach in applied mathematics and data science. TDA extracts topological information from data that cannot be obtained using statistical, physical, or other mathematical methods. Since topology is one of the most abstract fields in mathematics, TDA inherently reduces the dimensionality of high-dimensional data. Additionally, it simplifies geometric complexity and provides high-order data representations, extending conventional zeroth-order algorithms. Unlike classical topology or homology, TDA employs multiscale analysis, bridging topology and geometry through filtration.

A cornerstone of TDA is persistent homology \cite{edelsbrunner2002topological, zomorodian2005computing}, an algebraic topology tool for analyzing point clouds or discrete data. Persistent homology constructs a family of topological spaces across spatial scales, capturing topological invariants of varying dimensions. For further details, see \cite{carlsson2009topology, edelsbrunner2008persistent}. Its standard representations, i.e., persistence barcodes \cite{ghrist2008barcodes} and persistence diagrams \cite{edelsbrunner2002topological}, encode the birth, death, and persistence of features (e.g., connected components, holes, cavities) during filtration. Initially used for qualitative analysis, persistent homology has evolved into a robust tool for quantitative analysis and machine learning \cite{xia2014persistent,xia2015persistent,cang2015topological}. Applications span image processing \cite{bae2017beyond,clough2020topological}, neuroscience \cite{dabaghian2012topological}, computational chemistry \cite{townsend2020representation}, biology \cite{cang2017topologynet,cang2018integration,gameiro2015topological,kovacev2016using}, nanomaterials \cite{lee2017quantifying,xia2015persistent}, crystalline materials \cite{jiang2021topological}, and complex networks \cite{horak2009persistent}.
 
Despite its utility, persistent homology has notable constraints \cite{wei2025persistent}. It primarily captures topological invariants, often overlooking geometric shape evolution during filtration when no topological changes occur. 
Moreover, it struggles with localized topological information in data and cannot directly analyze data on differential manifolds or 1D curves embedded in 3-space (e.g., knots, links).

Efforts have been made to develop persistent topological Laplacians for various data structures, providing representations that can be integrated into machine learning models. For point cloud data and graphs, persistent combinatorial Laplacians were introduced by Wang et al. \cite{wang2020persistent} and have been extensively studied in recent years \cite{wang2021hermes, memoli2022persistent, gulen2023generalization, liu2024algebraic}. The harmonic spectra of persistent combinatorial Laplacians capture changes in topological invariants, which could recover the topological information captured by persistent homology, while the non-harmonic spectra encode geometric shape changes in the filtration process \cite{memoli2022persistent}. In addition to persistent combinatorial Laplacians, a variety of persistent topological Laplacians have been proposed for other topological domains. These include persistent sheaf Laplacians \cite{wei2025persistent} for labeled point cloud data, persistent path Laplacians \cite{wang2023persistent} for directed graphs and networks, persistent hypergraphy Laplacians \cite{liu2021persistent} and hyperdigraphy Laplacians \cite{chen2023persistent} for graphs and networks with complex relationships, persistent directed flag Laplacian for flag complexes \cite{jones2025persistent, zia2025persistent}. Persistent topological Laplacians have been successfully applied to protein-ligand binding prediction \cite{meng2021persistent, su2024persistent, chen2024multiscale}, interactomic network modeling \cite{du2024multiscale}, gene expression analysis \cite{cottrell2023plpca},  deep mutational scanning \cite{chen2023topological}, and protein engineering \cite{qiu2023persistent}.
In particular, one of the most impressive achievements of persistent topological Laplacians is the prediction of emerging dominant SARS-COV-2 variants BA4 and BA5 \cite{chen2022persistent}.

Another important spectral approach is quantum persistent homology or persistent Dirac, first introduced by Ameneyro et al. \cite{ameneyro2024quantum}. In quantum theory, the Dirac operator is for momentum, whereas the Laplacian is for kinetic energy, which reveals their relationship. Persistent Dirac operators have been studied in many topological domains with various applications \cite{wee2023persistent, suwayyid2024persistent}.  

For data on (differential) manifolds, a natural approach is persistent de Rham cohomology, rooted in differential topology or differential geometry. However, persistent de Rham cohomology has rarely been implemented independently. In contrast, evolutionary de Rham-Hodge theory, which includes persistent de Rham cohomology, was first proposed by Chen et al. for data on differential manifolds in Lagrangian representation, i.e., tetrahedral meshes \cite{chen2021evolutionary}. This approach draws on differential geometry or differential topology. Recently, persistent Hodge Laplacians have been developed in the Eulerian representation, i.e., Cartesian grid defined by level-set functions \cite{su2024persistent}. The similarity and difference between combinatorial Laplacians and Hodge Laplacians have been thoroughly discussed in \cite{ribando2024graph}. However, the term ``Hodge Laplacians" has also been used in the combinatorial setting, which can be somewhat misleading. Essentially, although the orthogonal decomposition can be carried out for both combinatorial Laplacians and Hodge Laplacians, only the Hodge decomposition of vector fields on differential manifolds can give rise to curl-free, divergence-free, and harmonic components, i.e., the so-called Helmholtz-Hodge decomposition \cite{ladyzhenskaya1969mathematical,zhao20193d}.  
Persistent Sheaf theory \cite{yegnesh2016persistence,wei2025persistentsheaf} and the persistent interaction topology \cite{liu2023interaction}  can effectively describe the local topology in the data.  
 For manifolds with boundary, a more refined topology-preserving 5-component Hodge decomposition was developed to further split the harmonic component into three subcomponents under certain boundary conditions \cite{poelke2016boundary, poelke2017hodge, zhao20193d,su2024hodge, su2024topology}. Hodge decomposition is an important tool in analyzing vector fields by extracting different dynamic features, and has many applications in computational fluid dynamics~\cite{yang2021clebsch, Yin2023FluidCohomology}, geometric modeling~\cite{wang2021computing}, and spectral data analysis~\cite{keros2023spectral}. Recently, it has been applied to the study of single-cell RNA velocities~\cite{su2024hodge} and has been integrated with deep learning models for medical image analysis \cite{liu2025manifold}.

For data on 1D curves embedded in 3-space, natural approaches stem from low-dimensional topology. Multiscale Gauss linking integral method was proposed to deal with knot-like data \cite{shen2024knot}. Additionally, multiscale Jones polynomial and persistent Jones polynomial have been introduced \cite{song2025multiscale}. Recently, persistent Khovanov homology has been developed \cite{shen2024evolutionary, liu2024persistent}.  
While Khovanov homology bridges geometric topology and algebraic topology, persistent Khovanov homology extends TDA to 1D data embedded in 3-space. 

Cang and Wei introduced the first integration of persistent homology and deep neural networks, including topological deep learning (TDL), in 2017 \cite{cang2017topologynet}. 
By leveraging the high-level abstraction of topology, TDL dramatically reduces dimensionality, simplifies geometric complexity, captures high-order interactions, and offers an interpretable learning framework. Persistent homology-based machine learning models have been successfully applied in various fields in science, engineering, medicine, industry, defense, etc.    \cite{papamarkou2024position,hensel2021survey}.  
However, topology inherently simplifies data, which implies the irreversible loss of certain information. Therefore, competitive performance from TDL mostly involves intrinsically complex data, such as those found in biological sciences \cite{cang2018representability}.  Perhaps some of the most compelling examples of applications in which TDL consistently demonstrates its advantages over other competing methods are the victories of TDL in the D3R Grand Challenges  \cite{nguyen2019mathematical,nguyen2020mathdl}, the discovery of SARS-CoV-2 evolution mechanisms \cite{chen2020mutations,wang2021mechanisms},
and the successful forecasting of SARS-CoV-2 variants BA.2  \cite{chen2022omicronBA2}
and BA.4/BA.5  \cite{chen2022persistent} about two months in advance. 

Machine learning typically requires all samples to have a uniform input vector. Therefore, featurization or vectorization of TDA information is often needed to convert topological representations, such as persistence barcodes, diagrams, images, landscapes, Betti curves, etc., into fixed-length feature vectors.
Many attempts have been made to define kernels \cite{reininghaus2015stable,kusano2016persistence,carriere2017sliced,le2018persistence}. Others focus on developing featurization techniques, such as statistical summaries \cite{cang2015topological,asaad2022persistent,cang2017topologynet,cang2018representability,atienza2019persistent,adcock2013ring,kalivsnik2019tropical,zielinski2019persistence,zielinski2021persistence}, persistence bins \cite{cang2017topologynet}, topological images \cite{cang2018representability}, and Wasserstein metrics \cite{cang2018representability}.

Numerous reviews have been conducted on persistent homology and its featurization \cite{edelsbrunner2008persistent,pun2022persistent,ali2023survey,chung2022persistence,barnes2021comparative,conti2022topological}, as well as a dedicated survey on persistent topological Laplacians \cite{wei2025persistent}. However, there is no review of the generalizations of persistent homology and their vectorization. The aim of this paper is to present a high-level overview of various TDA approaches beyond persistent homology and their machine learning featurization methods, without too many technical details or formulas. This work is intended to be accessible to new researchers in the field and to provide a broad perspective on various techniques and their applications.

The rest of the review is organized as follows. Section~\ref{sec:algebraicTDA} is devoted to algebraic topology approaches and techniques defined on point-cloud data. These approaches include persistent homology and its generalizations to various topological domains, persistent topological Laplacians on various topological domains, persistent Dirac operators, persistent sheaf theory, persistent Mayer topology, and persistent interactions. Section~\ref{secDifferentialT} reviews methods that originate in differential topology or differential geometry. We discuss the de Rham-Hodge theory for data on manifolds with boundary, a commonly used setting in various applications. Persistent de Rham-Hodge theory, including persistent de Rham cohomology and persistent Hodge Laplacian, is discussed. Section~\ref{secGeometricT} is devoted to various geometric topology-rooted TDA methods, such as multiscale Gauss linking integrals, persistent Jones polynomial, persistent Khovanov homology, and Khovanov Laplacian.
Vectorization techniques are reviewed in Section~\ref{sec:MachineLearningFeaturization}. We discuss commonly used methods, such as statistical summaries, Betti curves, persistence landscapes, and persistence images. We also provide relevant featurization techniques applicable to spectral methods, including eigenvalues, eigenvectors, and Hodge decompositions. In addition, we discuss the choice of TDA methods for effective analysis of different types of datasets. Finally, in Section~\ref{sec:conclusion}, we provide concluding remarks.

\section{Algebraic topology approaches }
\label{sec:algebraicTDA}

Topological data analysis for point cloud data has grown rapidly in recent years, providing a range of methods built upon tools in algebraic topology.
These methods offer multiscale characterizations of point cloud data by extracting topological and geometrical features through a filtration process. The resulting features provide rich information on the underlying data and have been successfully integrated into many downstream analyses, including classification, clustering, and other machine learning tasks. 

The most fundamental method, persistent homology \cite{edelsbrunner2002topological, zomorodian2005computing}, captures the persistence of topological features, such as connected components, holes, and cavities, on different scales in a filtration process. However, it has several limitations as it focuses mostly on the topology without considering the geometric evolution in the filtration or additional structures associated with data points, and it can be computationally expensive for large datasets, see Section~\ref{sec:perspectives}.

Recently developed methods, such as persistent combinatorial Laplacians \cite{wang2020persistent} with their variations \cite{wang2023persistent,liu2021persistent,chen2023persistent,jones2025persistent, zia2025persistent}, and the persistent Dirac operator \cite{ameneyro2024quantum, ameneyro2023quantum}, along with their variations \cite{suwayyid2024persistentB, suwayyid2024persistent}, generalize persistent homology and address some of its limitations. Persistent combinatorial Laplacians \cite{wang2020persistent} integrate the multiscale analysis with combinatorial Laplacians, allowing for the analysis of both topological and geometric or combinatorial structures across different scales. Their harmonic spectra fully recover the topological information from persistent homology, while the non-harmonic spectra offer additional geometric or combinatorial features of the data. The persistent Dirac \cite{ameneyro2024quantum, ameneyro2023quantum}, whose square yields the persistent combinatorial Laplacians, provides efficient algorithms based on quantum computing with exponential speed-up compared to classical methods in topological data analysis.
There are also other formulations based on alternative homology theories, such as the persistent sheaf Laplacians \cite{wei2025persistentsheaf}, which further extend the persistent combinatorial Laplacians, allowing the analysis of labeled or weighted datasets; the Mayer topology \cite{shen2024persistent, suwayyid2024persistent}; and the persistent interaction topology \cite{liu2023interaction, liu2024persistent}. These tools offer new perspectives for capturing topological features that are not accessible through standard homology.

In the following, we first review persistent homology with its variations and limitations, then turn to the persistent combinatorial Laplacians and their various variations for analyzing point cloud data with additional structures. We then describe the persistent Dirac operator and its variations. Finally, we discuss other formulations, including the persistent sheaf topology, the persistent Mayer topology, and the persistent interaction topology.

\subsection{Persistent homology}
\label{sec:persistentHomology}

Persistent homology is a multiscale extension of simplicial homology, which is a tool in algebraic topology that classifies simplicial complexes based on their homology groups. The homology groups measure the topological features in the underlying simplicial complex, such as connected components, holes, and cavities. Their ranks, given by Betti numbers, count the number of these topological features. For example, the $0$th Betti number provides the number of connected components, the $1$st Betti number gives the number of circles or holes, and the $2$nd Betti number represents the number of cavities. 

However, simplicial homology does not directly apply to point-cloud data, as such data do not have a well-defined topological structure unless explicitly constructed. In addition, it is sensitive to scales, as different choices of scales can lead to different homology groups. Persistent homology addresses these limitations by introducing filtration and persistence, providing a more robust characterization of the underlying data.

\subsubsection{Simplicial homology}

The simplicial complexes are made of simplices, where a $k$-simplex is defined to be the convex hull of $k\!+\!1$ affinely independent points $v_0, v_1, \cdots, v_k$, and has dimension $k$. Examples include a point ($0$-simplex), a line segment ($1$-simplex), a triangle ($2$-simplex), and a tetrahedron ($3$-simplex). These points $v_0, v_1, \cdots, v_k$ are called the vertices of the $k$-simplex, and the simplices formed by a subset of the set of all these vertices are called faces.

A simplicial complex $K$ is then defined to be a finite set of simplices satisfying the gluing conditions, where all faces of a simplex are also in $K$, and the intersection of any two simplices in $K$ is either empty or a common face of both. The dimension of $K$ is given by the largest dimension of its simplices. A linear combination of all $k$-simplices of $K$ with coefficients in a field (commonly $\mathbb{Z}_2$) defines a $k$-chain. The set of all $k$-chains, together with addition, forms an abelian group $\mathcal{C}_k(K)$, called the $k$-th chain group. The set of all $k$-simplices of $K$ serves as a basis of $\mathcal{C}_k(K)$.

To define the boundary operator $\partial$, each simplex of $K$ needs to be oriented (not necessary when using $\mathbb{Z}_2$ coefficients as $-1 = +1$). An oriented $k$-simplex is represented as an ordered list of its vertices, denoted by $[v_0, v_1, \cdots, v_k]$. The orientation is crucial for constructing the chain complex and defining homology, as it ensures that $\partial\partial = 0$. The $k$-th boundary operator $\partial_k: \mathcal{C}_k(K)\to \mathcal{C}_{k-1}(K)$ is a linear map that acts on each $k$-simplex by mapping it to an alternating sum of its $(k\!-\!1)$-faces, where each is obtained by removing one vertex from the $k$-simplex, i.e.,
\begin{align}
    \partial_k[v_0, v_1, \cdots, v_k] = \sum_{i} (-1)^i[v_0, \cdots, \hat{v_i}, \cdots, v_k]
\end{align}
The boundary operator satisfies $\partial_k\partial_{k+1} = 0$. A chain complex is given by a sequence of chain groups connected by the boundary operator:
\begin{align}
    \cdots\xrightarrow{\partial_{k+2}}\mathcal{C}_{k+1}
    \xrightarrow{\partial_{k+1}}\mathcal{C}_{k} \xrightarrow{\: \partial_{k} \:}\mathcal{C}_{k-1} 
    \xrightarrow{\partial_{k-1}}\cdots
\end{align}
The $k$-th simplicial homology group of $K$ is defined as the quotient of the kernel of the boundary operator $\partial_k$ modulo the image of $\partial_{k+1}$:
\begin{align}
    H_k(K) = \ker{\partial_k }/\img{\partial_{k+1}},
\end{align}
where $\ker{\partial_k}$ and $\img{\partial_{k+1}}$ are called the group of $k$-cycles and the group  of $k$-boundaries, respectively.
The homology group $H_k(K)$ captures the $k$-cycles that are not $k$-boundaries, i.e., the $k$-dimensional holes in the simplicial complex. Its rank, known as the $k$-th Betti number $\beta_k$, counts the number of these holes. 

\subsubsection{Filtration and persistence}

A filtration is given by a nested sequence of simplicial complexes. For point cloud data, one needs to construct such a filtration by building simplicial complexes to represent data points and their connections in various dimensions. The three most commonly used complexes are the Vietoris-Rips complex \cite{vietoris1927hoheren}, the \v{C}ech complex \cite{edelsbrunner2010computational}, and the Alpha complex \cite{edelsbrunner2011alpha}. There are also other constructions such as witness complex \cite{de2004topological}, Cover complex \cite{gudhi:CoverComplex}, and Tangential complex \cite{gudhi:TangentialComplex}.

The Vietoris-Rips complex \cite{vietoris1927hoheren} is constructed by considering pairwise intersections of balls centered at each data point within a fixed distance, leading to a simplex whenever its vertices have pairwise overlapping regions. The \cech complex \cite{edelsbrunner2010computational} requires a nonempty intersection of all balls in a set, meaning that a simplex is included only if all corresponding balls interSection The Alpha complex can be constructed in a similar manner to the \cech complex. However, each ball is restricted to the Voronoi cell associated with each point. The Alpha complex \cite{edelsbrunner2011alpha} is a subset of the \cech complex but has the same homology. By varying the radius of the balls, one obtains a filtration of simplicial complexes. The \cech complex involves checking for a large number of intersections and thus comes with significantly higher computational costs for large datasets, making the Vietories-Rips complex the most common choice in real applications. An alternative to the \cech complex is the Alpha complex, which is particularly computationally more efficient when the ambient space has dimension $2$ or $3$ \cite{edelsbrunner2010computational}.

Let $K$ be the largest simplicial complex constructed from a point cloud. Then a filtration of this point cloud is a nested sequence of the subcomplexes of $K$:
\begin{align}
    \emptyset\subset K^0\subset K^1 \cdots\subset K^m = K.
\end{align}
Denote by $i^{i, j}$ the inclusion maps from $K_i$ to $K_j$ with $i<j$. Its induced map $f_k^{i, j}: H_k(K^i)\to H_k(K^j)$ on the homology groups tracks how the $k$-dimensional holes evolve across the filtration, more precisely, if the $k$-dimensional homology classes in $ H_k(K^i)$ still exist in $ H_k(K^j)$. 
Denote by $\partial_{k}^i: \mathcal{C}_{k}(K^i)\to \mathcal{C}_{k-1}(K^i)$ the $k$-th boundary operator on simplicial complex $K^i$. The $k$-th persistent homology group is defined to be the image of $f_k^{i, j}$:
\begin{align}
    H_k^{i, j} = \img f_k^{i,j} = \ker\partial_k^i/\left(\img\partial_{k+1}^j\cap\ker\partial_k^i\right),
\end{align}
and its rank leads to the corresponding $k$-th persistent Betti number $\beta_k^{i, j} = \rank H_k^{i, j}$. The birth and death of a topological feature can then be defined by examining the elements in the persistent homology groups. A feature is said to be born at $K^i$ if it first appears at $K^i$, and it is said to die at $K^j$ if it merges with a feature at $K^j$ that already exists before $K^j$. The persitence of $\gamma$, denoting the lifetime of a topological feature, is then given by $j\!-\!i$, and considered infinite if it never dies in the filtration.

These birth and death information of topological features gives rise to various representations of persistent homology, such as persistence barcodes \cite{ghrist2008barcodes}, persistence diagrams \cite{edelsbrunner2002topological}, persistent Betti numbers, Betti curves, persistence landscapes \cite{bubenik2015statistical}, persistence images \cite{adams2017persistence}. 
More details about these featurization techniques will be discussed in Section~\ref{sec:MachineLearningFeaturization}.

Persistent homology has been successfully applied in various fields, ranging from image processing to complex networks \cite{bae2017beyond,clough2020topological, dabaghian2012topological, townsend2020representation, cang2017topologynet, cang2018integration,gameiro2015topological,kovacev2016using, lee2017quantifying, jiang2021topological,horak2009persistent}. The most significant achievements would be its integration with deep learning models, which has led to exceptional performance in the D3R grand challenges \cite{nguyen2019mathematical,nguyen2020mathdl} and in predicting the evolutionary mechanism of SARS-CoV-2 \cite{chen2020mutations,wang2021mechanisms}.

\subsubsection{Other variations}

In addition to the standard formulation of persistent homology for simplicial complexes, several variations have been proposed for analyzing data with different structures or for addressing specific application needs. We briefly review these methods below.

\paragraph{Multiparameter persistent homology} The multiparameter persistent homology was first explored in \cite{carlsson2007theory}, while a variation in size theory was introduced in \cite{frosini1999size}. The framework extends classical homology by considering two or more parameters instead of one, and it is suitable for analyzing data with outliers or variations in density, data equipped with real-valued functions and functional data with large local noise. For more introductory resources, see \cite{karaguler2021survey, dey2022computational, schenck2022algebraic, botnan2023introduction}. Discussions on the computation and software can be found in \cite{botnan2023introduction}.

\paragraph{Zig-zag persistent homology} The Zig-zag persistent homology was introduced in \cite{carlsson2010zigzag} to handle a more general setting where the sequence of spaces is connected with maps in arbitrary directions, in contrast to the classical persistent homology where the inclusion maps follow in the same direction. The framework is powerful since it allows simplicial complexes in the filtration to evolve with the insertion and deletion of simplices, which is suitable for analyzing dynamic point clouds or dynamic graphs with varying vertex connections \cite{holme2012temporal, dey2021computing, dey2021updating, chen2021z}. Various software for computing the zigzag persistence have been proposed \cite{carlsson2009zigzag, dey2022fast}. The framework has also been extended to the multiparameter setting \cite{dey2021updating, dey2024computing}. 

\paragraph{Parametrized homology}
Parametrized homology \cite{carlsson2019parametrized} extends levelset zigzag persistent homology to a continuous-parameter setting. This framework leverages the theory of rectangle measures \cite{chazal2016structure} to define four distinct continuous-parameter persistence diagrams, each corresponding to one of the four bar types in the levelset zigzag barcode. These diagrams encode homological features and capture detailed information about how each feature disappears at both endpoints of its defining interval.

\paragraph{Localized homology}
Localized homology \cite{zomorodian2008localized} is proposed for identifying the locations of topological attributes such as tunnels and voids. Given a topological space with a cover, the localized homology tool can be used to determine the location (the specific cover element) of topological features within the space. Specifically, it constructs a larger blowup complex \cite{segal1968classifying} that is homotopy equivalent to the original space and has clear boundaries between cover elements. The homology bases are then computed on the blowup complex to localize the features to specific cover elements. An efficient algorithm is also provided for practical applications. 

\paragraph{Evolutionary homology}
Cang et al. \cite{cang2020evolutionary} introduced Evolutionary Homology (EH), a method generating time-dependent topological invariants ("evolutionary barcodes") for individual components within physical systems. EH works by coupling differential equations or chaotic oscillators to model system interactions, creating coupled dynamics. On the resulting trajectories of these oscillators, they define simplices, simplicial complexes, algebraic groups, and topological persistence using a time evolution-based filtration. This process reveals the topology-function relationship of individual components and is particularly effective for analyzing local topology within point cloud data. EH was successfully applied to realistic problems in biophysics.  

\paragraph{Multiresolution persistent homology} Multiresolution persistent homology was proposed for analyzing excessively large datasets \cite{xia2015multiresolution,xia2015multiresolution2}. The idea is to match the scale of interest in the data by employing the rigidity-density functions modulated through resolution parameters. By tuning these parameters, the original large-scale datasets can be represented at an appropriate lower resolution, enabling efficient topological analysis.

\paragraph{Multidimensional persistence}
The multidimensional persistence was first introduced in 2007 for analyzing the topology of multifiltrations \cite{carlsson2007theory}. Later, another version of multidimensional persistence \cite{xia2015multidimensional} was proposed for transforming both dynamic and single molecular data into image-based representations. For dynamical data, the image representation is constructed by repeatedly applying persistent homology to each data frame across time steps. For single molecular data, persistent homology is applied at multiple spatial scales to the volumetric density representation of the molecule, and the resulting topological features are aggregated to form the final image.

\paragraph{Weighted persistent homology} The weighted persistent homology is used to study the weighted data, such as weighted simplicial complexes, where each simplex is assigned a weight. There are primarily two types of methods for capturing the weight information: one approach uses  weights to refine the filtration process \cite{petri2013topological,edelsbrunner2013persistent}
, while the other revises the boundary operators of the chain complexes by incorporating the weights \cite{dawson1990homology,ren2018weighted}. The weighted persistent homology has been applied to various biomolecular data analysis applications \cite{meng2020weighted,pun2020weighted}.

\paragraph{Persistent path homology} Persistent path homology has been introduced for analyzing directed networks that incorporate asymmetric directional information \cite{chowdhury2018persistent}, which integrates the filtration method in TDA with the path homology theory. Path homology \cite{grigor2012homologies}, recently renamed GLMY homology \cite{ivanov2024simplicial}, has been extensively developed in both theory \cite{grigor2023homotopy,grigor2022advances,li2024singular,di2024path,li2024primitive} and applications \cite{wu2023metabolomic,chen2023path,gong2024topological,feng2025network}.

\paragraph{Persistent hypergraph homology} Persistent hypergraph homology \cite{bressan2019embedded} is a generalization of the standard persistent simplicial homology to hypergraphs. Using embedded homology theory, hypergraph homology is suitable for addressing systems with incomplete information \cite{liu2021hypergraph}, such as the coauthor network, where a three-person collaboration exists but not all pairwise collaborations exist. A stability result \cite{ren2020stability} and a discrete Morse theory \cite{ren2021discrete} for hypergraphs have been established, and an algorithmic implementation for persistent hypergraph homology \cite{liu2024computing} has also been developed.

\paragraph{Persistent super-hypergraph homology} Persistent super-hypergraph homology \cite{grbic2022aspects} provides a unified and generalized framework that encompasses persistent simplicial homology, persistent hypergraph homology, and persistent $\Delta$-set homology as special cases. Compared with a simplicial complex and a hypergraph, a distinct characteristic of super-hypergraph is that each super-edge is not uniquely determined by its vertices, implying the presence of additional internal structure. This makes it a flexible framework to capture the complexity of various systems \cite{liu2024intcomplex,feng2024hypernetwork}.

\paragraph{Graph-complex-based persistent homology}
In TDA, modeling data with an appropriate simplicial complex is as critical as data representation in machine learning. The most commonly used constructions are the Vietoris-Rips complex \cite{zomorodian2010fast} and the Alpha complex \cite{edelsbrunner2011alpha}, maainly due to the availability of well-established and user-friendly software packages. Several studies have explored the use of other graph-based complexes, such as the flag complex \cite{horak2009persistent,lutgehetmann2020computing}, Dowker complex \cite{liu2022dowker}, neighborhood complex \cite{liu2021neighborhood}, Hom-complex \cite{liu2022hom}, and others \cite{jonsson2008simplicial,grbic2022aspects}. Selecting the most suitable complex for a given application, rather than just using the popular choices, is worth further investigation.

\paragraph{Extended persistent homology}
Distances between persistence diagrams are often used to compare topological features of different topological spaces. However, persistent diagrams of topological spaces with different Betti numbers can lead to infinite distances under any filtration. Extended persistent homology \cite{cohen2009extending,turner2024extended} addresses this limitation by extending the filtration process so that all homology classes are assigned a birth and a death.

\paragraph{Join persistent homology}
Join persistent homology \cite{wang2025join} is based on the join operation in topology. The idea is that different topological features can be captured from multiple perspectives by employing various types of unit in the join process. This allows the method to reveal diverse structural information within the data. Join persistent homology also provides a solution to the issue of infinite distance between persistence diagrams of spaces with different Betti numbers.

\paragraph{Persistent cohomology}
Persistent cohomology is a sister counterpart to persistent homology, offering computational advantages in time and space. In practice, people often compute the persistent cohomology and subsequently derive persistent homology via the duality \cite{de2011dualities}, as the former can be more efficient to compute. An example is the GUDHI package \cite{maria2014gudhi}. Many applications rely on the original persistent cohomology algorithm \cite{de2009persistent}. Several studies have explored incorporating domain-specific information into cohomology classes to enhance the characterization of background-related information in applications \cite{cang2020persistent,liu2021hypergraph}.

\paragraph{Cayley-persistence}
Cayley-persistence is formulated based on the group grading \cite{bi2022cayley}. Similar to the standard $\mathbb{Z}$-graded persistence, it has the decomposition results. The key difference is its ability to capture more flexible persistence patterns. It allows the emergence of persistence not only on the intervals and grids, but also on more complex structures introduced by the underlying group grading.

\paragraph{Morse theory}
Morse theory \cite{morse1925relations} provides a framework for analyzing and characterizing the topology of shapes. The central idea is to represent a shape as a topological space equipped with a Morse function that captures critical topological features through its critical points.  Discrete Morse theory is integrated with persistent homology for unweighted and undirected networks
\cite{kannan2019persistent}.  Computational algorithm for Morse theory has been developed \cite{gunther2012efficient}. 
 The persistent Homology of Morse decompositions has been applied to combinatorial dynamics \cite{dey2019persistent}. 
Among the various descriptors derived from Morse theory, the Reeb graph is probably the most widely used \cite{reeb1946points,biasotti2008reeb,ge2011data,zhao2018protein}. A key property of the Reeb graph is its ability to encode the topological structure of a shape into a one-dimensional graph, independent of the intrinsic dimension of the underlying manifold.

\paragraph{Conley index}
Much of dynamical system analysis involves the existence and structure of invariant sets. Conley index theory, developed by Conley and his followers \cite{conley1978isolated,mischaikow1999conley,chen2008efficient}, is an important topological tool for studying such invariant sets. Rather than analyzing arbitrary invariant sets, Conley index theory focuses on isolated invariant sets, which are more robust to continuous perturbations than general invariant sets. For an isolated invariant set $S$, a quotient space can be constructed from its index pair, the homotopy Conley index of $S$ is the homotopy type of this quotient space, and the cohomological Conley index of $S$ is the Alexander-Spanier cohomology of the quotient space with integer coefficients. The Conley index is invariant under certain deformations of the dynamical system and can be used to establish the existence of fixed points and periodic orbits of the dynamical system. Recently, the Conley index has been extended to the persistence setting for data analysis \cite{dey2020persistence,dey2022persistence}.

\paragraph{Persistent equivariant cohomology}
For a dataset with a group action on it, persistent equivariant cohomology captures not only the underlying shape of the dataset but also properties of the group action, namely, the symmetries inherent in the data. In \cite{adams2024persistent}, the authors provide an explicit description of the persistent equivariant cohomology of the circle action on the Vietoris-Rips metric thickening of the circle. In particular, they derive explicit formulas for the persistent equivariant cohomology of these Vietoris-Rips thickening. When equivariant cohomology is based on a free group action, it can be viewed as a special case of twisted cohomology, with the twist arising from the group action. Various studies have explored the theory of twisted (co)homology in different contexts \cite{li2017twisted,grigor2016cohomology,zhang2022twisted}.

\paragraph{Persistence over posets}
Persistence over posets \cite{kim2023persistence} is a general framework where persistence modules are defined on posets, with multiparameter persistence as a special case where the indexing poset is $\mathbf{R}^n$.
As the most known invariant for multiparameter persistence, the rank invariant \cite{carlsson2007theory} has evolved into the generalized rank invariant \cite{patel2018generalized} for persistence over posets by naturally extending its domain. In particular, for any interval decomposable persistence module of finite-dimensional vector spaces, its barcode can be extracted from the generalized rank invariant by the principle of inclusion-exclusion \cite{kim2021generalized}. This implies that the generalized rank invariant is a complete invariant for interval decomposable persistence modules, while the rank invariant is not a complete invariant for interval decomposable multiparameter persistence modules \cite{botnan2020rectangle}.

\subsubsection{Additional perspectives}
\label{sec:perspectives}

Despite the great utility of persistent homology, there are a number of directions in where it is desirable to extend it.  These extensions parallel much of the work that has been done in mainstream algebraic topology, i.e., the topology of actual spaces and manifolds rather than point clouds.  Here are a few of them.  
\begin{itemize}

\item{{\bf Spectral analysis:} Note that persistent homology fails to capture the homotopic geometric shape changes within a filtration, meaning that it cannot distinguish between different shapes with the same topology \cite{hernandez2025persistence}. Several frameworks have been developed to address this limitation based on the spectral analysis of Laplacians \cite{wei2025persistent}. These methods not only recover topological information from persistent homology but also capture additional information that describes the homotopic geometric shape changes in the filtration.}

\item{{\bf Analytic methods:} From the early days of the subject, it was understood that in the case of manifolds, there are approaches to computing a version of homology (homology with $\mathbb{R}$ coefficients) using differential forms. These methods rely on {\em de Rham cohomology} \cite{warner1983foundations}. The advantage of these methods is that they tie the homological qualitative information together with the analysis on the manifold. Another advantage is that the {\em Hodge theorem} \cite{warner1983foundations} provides a choice, within a cohomology class, of a harmonic representative, which yields clarity about the significance of the class.  A key ingredient within this analysis is the Hodge Laplacian, a partial differential operator on differential forms.    }

\item{{\bf Sheaves:} After the initial development of homological methods in the early part of the 20th century, by 1950 it became understood that in order to gain more power in the use of the homological methods for the understanding of spaces, one needed the more complex notion of {\em cohomology with coefficients in a sheaf} pioneered by J. Leray \cite{Leray1946sheaf}.  This notion became critical to the extension of topological techniques to algebraic geometry and number theory \cite{weil1949numbers, deligne1974conjecture, deligne1980conjecture}.}

\item{{\bf Local to global computation: } The direct computation of topological invariants using the linear algebraic methods, is unmanageable except in the smallest toy examples. For this reason, algebraic topologists were forced to develop methods for parallelizing the computations, and for using structures inherent in the problem, such as a fibration structure.  Examples of these ideas are the use of Mayer-Vietoris long exact sequences and spectral sequences, the Serre spectral sequence for a fibration, Eilenberg-Moore spectral sequences \cite{mccleary2001user}.}

\item{{\bf A priori knowledge:} Often a topological problem is equipped with a priori information, such as a map to a reference base space or some choices of subspaces within the space.  Homological methods have been refined to capture this additional information in more subtle invariants.  }

\end{itemize}

\subsection{Persistent combinatorial Laplacians}

To address the limitation that persistent homology fails to capture the homotopic geometric changes during a filtration, recent advances have introduced persistent combinatorial Laplacians \cite{wang2020persistent}, and some other variations of persistent Laplacians \cite{wang2023persistent,liu2021persistent,chen2023persistent,jones2025persistent, zia2025persistent,wei2025persistentsheaf} that encode in spectra both the topological and the geometric or combinatorial structures of data. Compared to persistent homology, these persistent Laplacians enable deeper analysis of datasets. Their harmonic spectra (i.e., zero eigenvalues) fully recover the topological outputs from the persistent homology, while the non-harmonic spectra (i.e., nonzero eigenvalues) capture additional geometric or combinatorial information about the shape of data. This additional information is particularly useful in describing the homotopic geometric shape changes in a filtration, which the persistent homology fails to detect. A comprehensive survey of persistent Laplacians is given by Wei and Wei \cite{wei2025persistent}. 

These persistent topological Laplacians can be viewed as multiscale formulations of combinatorial Laplacians \cite{eckmann1944harmonische}, which originate from spectral graph theory and generalize the notion of graph Laplacians \cite{kirchhoff1847ueber} to higher-dimensional simplicial complexes. The study of the spectra of combinatorial Laplacians is crucial in understanding the shape of data and has gained much attention in recent years \cite{goldberg2002combinatorial, gundert2014higher, horak2013spectra}.

Below, we begin with an overview of combinatorial Laplacians, then discuss the persistent combinatorial Laplacians, and finally cover several other variations of persistent topological Laplacians.

\subsubsection{Combinatorial Laplacians}

The combinatorial Laplacians extend the classical graph Laplacians to simplicial complexes \cite{eckmann1944harmonische}. Compared to graph Laplacians, which only encode connections between nodes, the combinatorial Laplacians are constructed using the boundary operators that capture the relationships between simplices even in higher dimensions, including vertices, edges, triangles, and higher-dimensional structures.

The $k$-th combinatorial Laplacian $L_k$ on the $k$-th chain group $\mathcal{C}_k(K)$ is given as follows:
\begin{align}\label{eq.Lk.combinatorial}
    L_k = \partial_{k+1}\partial^*_{k+1} + \partial^*_k\partial_k,
\end{align}
where $\partial^*_k$ is the adjoint operator of the boundary operator $\partial_k$ with respect to an inner product defined on the chain group $\mathcal{C}_k(K)$. Note that here the combinatorial Laplacian $L_k$ is defined on the chain group, a corresponding construction can also be made on the cochain group, see~\cite{horak2013spectra}.
The combinatorial Laplacian $L_k$ is positive semi-definite symmetric, and thus its spectrum consists of only non-negative real eigenvalues. In the case that $K$ is a graph, $L_0 = \partial_1\partial^*_1$ reduces to the graph Laplacian. In addition, the combinatorial Laplacians lead to a combinatorial version of Hodge decomposition:
\begin{align}
    \mathcal{C}_k = \img\partial_{k+1}\oplus\ker L_k\oplus\img\partial^*_k,
\end{align}
see \cite{lim2020hodge} for further details. 

In the case where the inner product on the chain group $\mathcal{C}_k(K)$ is defined by the Kronecker delta $\delta_{i,j}$, i.e., all basis elements in $\mathcal{C}_k(K)$ are considered mutually orthogonal, the boundary operator $\partial_k$ can be represented as a matrix of size $N_{k-1}\times N_{k}$ and its adjoint $\partial^*_k$ given as the transpose of the matrix. Here $N_k$ denotes the number of $k$-simplices.

It has been shown in \cite{eckmann1944harmonische} that the kernel of the $k$-th Laplacian on a simplicial complex $K$ corresponds to its $k$-th simplicial homology group:
\begin{align}
    \ker L_k = H_k(K),
\end{align}
thereby capturing the topological information of the simplicial complex. Its dimension is given by the $k$-th Betti number. In addition, the non-harmonic spectra capture the geometrical or combinatorial features of the underlying structure. In particular, as noted in graph theory, the first nonzero eigenvalue describes the connectedness of a graph.

\subsubsection{Persistent combinatorial Laplacians}

To incorporate the ideas of persistence into the combinatorial Laplacians, Wang et al. \cite{wang2020persistent} introduced the persistent combinatorial Laplacians, which is the first work, both theoretically and practically, to establish a connection between the persistent homology and the combinatorial Laplacians. A notion of persistent Laplacian also appeared in a seminar note \cite{Lieutier2014HarmonicForms}. A comprehensive theoretical study of persistent combinatorial Laplacians can be found in \cite{memoli2022persistent} and a stability result is found in  \cite{liu2024algebraic}.

Given a pair of simplicial complexes $K^i\subset K^j$ in a filtration, the $k$-th persistent combinatorial Laplacian on the $k$-th chain group $\mathcal{C}_k(K^i)$ is defined as follows:
\begin{align}
    L_k^{i, j} = \partial_{k+1}^{i, j}\left(\partial_{k+1}^{i, j}\right)^* + \left(\partial_k^i\right)^*\partial_k^i,
\end{align}
where $\partial_k^i$ is the $k$-th boundary operator on $K^i$, $\partial_{k+1}^{i, j}$ is called the persistent boundary operator, given by the restriction of $\partial_{k+1}^j$ to a subgroup $\mathcal{C}_{k+1}^{i, j}$ of the $(k\!+\!1)$-th chain group $\mathcal{C}_{k+1}(K^j)$ so that the images of its elements are in $\mathcal{C}_k(K^i)$, and $\left(\partial_k^i\right)^*$ and $\left(\partial_{k+1}^{i, j}\right)^*$ are their adjoint operators, respectively. In the case that $i = j$, the $k$-th persistent boundary operator equals the standard boundary operator, i.e., $\partial_{k}^{i, i} = \partial_{k}^{i}$, and thus the $k$-th persistent combinatorial Laplacian reduces to the standard $k$-th combinatorial Laplacian on $K^i$. 
The persistent combinatorial Laplacian $L_k^{i, j}$, analogous to the combinatorial Laplacian, is positive semi-definite symmetric, with its eigenvalues real and non-negative. This operator also results in a persistent Hodge decomposition \cite{Lieutier2014HarmonicForms} as follows:
\begin{align}
    \mathcal{C}_k(K^i) = \img\partial_{k+1}^{i,j}\oplus\ker L_k^{i,j}\oplus\img\left(\partial^i_k\right)^*.
\end{align}

The implementation of persistent combinatorial Laplacians has been made available as an open-source software package, called HERMES \cite{wang2021hermes}. The algorithm is based on a practical construction of filtration of alpha complexes \cite{edelsbrunner2011alpha}. By constructing a single $k$-th boundary operator for the final complete alpha complex (i.e., the Delaunay tessellation \cite{delaunay1934sphere}) and certain projection matrices, the matrix representations of $\left(\partial_{k+1}^{i, j}\right)^*$ and $(\partial_{k}^{i})^*$ are both given simply by the transpose of the matrices representing $\partial_{k+1}^{i, j}$ and $\partial_{k}^{i}$. An alternative algorithm for computing persistent combinatorial Laplacians can be found in \cite{memoli2022persistent}. 

An immediate result is that the kernel of the $k$-th persistent combinatorial Laplacians is isomorphic to the $k$-th persistent homology group \cite{Lieutier2014HarmonicForms, memoli2022persistent}:
\begin{align}
    \ker L_k^{i,j} = H_k^{i,j},
\end{align}
and thus its dimension is the $k$-th persistent Betti number. The harmonic spectra of persistent combinatorial Laplacians, therefore, recover the topological information from the persistent homology, while the non-harmonic spectra provide additional geometric or combinatorial evolutionary features of the shape of the data.

Persistent combinatorial Laplacians have been applied to various biological problems, such as protein thermal stability \cite{wang2020persistent}, protein-ligand binding \cite{meng2021persistent}, and protein-protein binding problems \cite{wee2022persistent}, etc. In particular, they have also been successfully used to forecast the SARS-CoV-2 variants BA.4 and BA.5 \cite{chen2022persistent}. The advantage of persistent combinatorial Laplacians over persistent homology has been demonstrated in protein engineering \cite{qiu2023persistent}.

\subsubsection{Other variations}
\label{sec:persistentCombinatorialLaplacians.othervariations}

We now outline various variations of persistent combinatorial Laplacians proposed for analyzing point-cloud data with other additional structures, such as directional information or more complex relationships.

\paragraph{Persistent path Laplacians}
Persistent path Laplacians \cite{wang2023persistent} was introduced by Wang and Wei for analyzing point data with directional information, such as directed graphs and networks. The framework builds upon the path homology \cite{grigor2012homologies}, which was later extended to a persistent framework on directed networks \cite{chowdhury2018persistent}. A weighted path homology for weighted digraphs was introduced in \cite{lin2019weighted}. More recently, an efficient algorithm for computing the $1$-dimensional persistent path homology was proposed in \cite{dey2022efficient}. 

\paragraph{Persistent hyper(di)graphy Laplacians}
Persistent hypergraphy Laplacians \cite{liu2021persistent} and hyperdigraphy Laplacians \cite{chen2023persistent} are proposed for complex systems with and without directional information, respectively. These constructions are based on the infimum and supremum chain complexes in embedded homology \cite{bressan2019embedded}, enabling the definition of the corresponding infimum and supremum Laplacians. Particularly, when the hypergraph reduces to a simplicial complex, the infimum and supremum chain complexes coincide, recovering the simplicial Laplacian as a special case.

\paragraph{Persistent hopping path-Laplacian}
Persistent hopping path Laplacian \cite{liu2023persistent} generalizes the $k$-hopping path Laplacian of graphs \cite{estrada2012path} in two aspects: first, by extending from graphs to simplicial complexes and hypergraphs, and second, by incorporating persistence into the framework. By selecting different hopping parameters, this method provides a comprehensive description of pairwise, many-body, direct, and indirect interactions within the data structures.

\paragraph{Persistent directed flag Laplacian}
Persistent directed flag Laplacian \cite{jones2025persistent} is constructed based on the flag complex of directed graphs, which only includes those complete subgraphs that satisfy specific order constraints. By integrating the directed flag complex with the standard simplicial Laplacian, it provides an efficient tool for analyzing complex networks that exhibits asymmetric directional relationships. Persistent directed flag Laplacian has been applied to complex biomolecular data \cite{zia2025persistent}.  

\paragraph{Persistent Laplacian for simplicial maps}
Persistent Laplacian for simplicial maps generalizes the standard persistent Laplacian defined on a filtration with inclusion maps \cite{gulen2023generalization}. It maintains a connection to persistent homology, where the nullity of the persistent Laplacian equals the persistent Betti numbers associated with the simplicial map-induced space sequence. A key difference is that its persistence depends not only on the sequence of simplicial complexes, as in the standard setting, but also on the structures of the simplicial maps themselves.  

\paragraph{Multiscale Hochschild Laplacian}
Multiscale Hochschild Laplacian \cite{he2025multi} offers an efficient tool for analyzing digraph data by integrating Hochschild Laplacians with multiscale analysis. Particularly, the Hochschild Laplacian is based on the truncated path algebras of the digraphs, as the full path algebra of a digraph with directed cycles is infinite-dimensional, making it unsuitable for computation and practical applications. Similar to the combinatorial Laplacians, the kernel of the Hochschild Laplacian is isomorphic to the Hochschild cohomology in the corresponding dimension.  

\subsection{Persistent Dirac operator}

Recently, quantum algorithms have been introduced for the efficient computation of topological and geometrical features to address the computational complexity of persistent homology, especially for large point-cloud data \cite{lloyd2016quantum,ameneyro2024quantum, bianconi2021topological,baccini2022weighted, krishnagopal2023topology}. These methods allow the data to be accessed in quantum parallel on a quantum computer and the combinatorial simplicial complexes to be encoded as quantum states, enabling an exponential speed-up over the classical topological data analysis techniques. Central to these methods is the Dirac operator \cite{dirac1928quantum} on simplicial complexes, whose square yields exactly the combinatorial Laplacians of all degrees. It captures the same information as combinatorial Laplacians but provides additional information as it is more sensitive to local topological and combinatorial or geometrical features as a linear operator. The Dirac operator has been found particularly useful in analyzing topological signals \cite{calmon2023dirac}, and has also been applied to molecular representations \cite{wee2023persistent}.

The persistent Dirac operator   \cite{ameneyro2024quantum, ameneyro2023quantum} extends the Dirac operator to the persistent setting for analyzing point-cloud data. Its square yields persistent combinatorial Laplacians and thus can be seen as a generalization of the persistent combinatorial Laplacian technique, capturing the same topological information while encoding more local topological and geometric structures of the underlying data. The Dirac operator has also been extended for tracking the persistence of paths and hypergraphs \cite{suwayyid2024persistentB}. The Persistent Mayer Dirac \cite{suwayyid2024persistent} has been developed based on $N$-chain complexes, expanding the scope and improving the applicability of persistent Dirac-based techniques.

In the following, we briefly review the Dirac operator and the persistent Dirac operator. We also discuss the other variations for data with additional structures. 

\subsubsection{Dirac operator}

In the combinatorial case, the Dirac operator also connects simplices of consecutive dimensions in a simplicial complex $K$ as the combinatorial Laplacians. It is defined as the sum of the boundary operator and its adjoint, given by $D = \partial + \partial^*$. It can be viewed as the square root of the combinatorial Laplacian defined in the full chain complex, i.e., the direct sum of all chain groups $\mathcal{C}(K) = \bigoplus_k\mathcal{C}_k(K)$. When expressed in the Hermitian representation, its square leads to a block matrix with diagonals given by the combinatorial Laplacians in different degrees: 
\begin{align}\label{eq.Lk.dirac}
    D^2 = \begin{pmatrix}
        0 & \partial_{1} & 0 & \cdots\\
        \partial^*_{1} & 0 & \partial_{2} & \cdots\\
        0 & \partial^*_{2} & 0 & \cdots\\
        \vdots & \vdots & \vdots & \ddots
    \end{pmatrix}^2 = 
    \begin{pmatrix}
        L_1 & 0 & 0 & \cdots\\
        0 & L_2 & 0 & \cdots\\
        0 & 0 & L_3 & \cdots\\
        \vdots & \vdots & \vdots & \ddots
    \end{pmatrix}
\end{align}
Compared to the combinatorial Laplacian, which preserves the degree and maps a $k$-chain to a $k$-chain, the Dirac operator maps a $k$-chain to both $(k\!-\!1)$- and $(k\!+\!1)$-chains. In addition, the Dirac operator is a first-order differential operator, which means that it is more sensitive to local geometric and topological structures compared to the Laplacians, which are second-orders. The kernel of the Dirac operator encodes the Betti numbers of a simplicial complex across all degrees, and the nonzero eigenvalues are the square roots of those of the combinatorial Laplacians.

The key advantage of using the Dirac operator in topological data analysis lies in its suitability for efficient quantum computations. Lloyd et al. \cite{lloyd2016quantum} introduced the first quantum algorithm for estimating Betti numbers by using Dirac operators, with an exponential speed-up over the classical algorithm in topological data analysis. The method constructs a quantum state that represents the simplicial complex at each scale, and then use a quantum phase algorithm \cite{nielsen2010quantum}
to find the dimension of the kernel of the Dirac operator $D$ in Hermitian representation.
The method provides a multiscale version of Betti numbers without explicitly tracking the persistence in filtration.

\subsubsection{Persistent Dirac operator}

The Dirac operator was extended by Ameneyro et al. \cite{ameneyro2024quantum} to the persistent setting for the first time in a quantum implementation for computing the persistence of point clouds using the Dirac operator across different scales. 
The framework utilizes a single $k$-th boundary operator $\partial_k$ encoded using 1-qubit Pauli X gates defined on the set of all possible $k$-simplices in the filtration, and also a projection matrix $P_k^i$ from this set to the set of $k$-simplices in the simplicial complex $K^i$. The $k$-th persistent boundary operator is then given as $\partial_k^{i,j} = P_{k\!-\!1}^i\partial_k P_k^j$ for a pair of simplicial complexes $K^i\subset K^j$ stored as quantum states in the filtration. The $k$-th persistent Dirac operator as Hermitian is defined in \cite{ameneyro2024quantum} as
\begin{align}
    D_k^{i,j} = \begin{pmatrix}
        0 & \partial_{k}^{i,i} & 0\\[3pt]
        \left(\partial_{k}^{i,i} \right)^* & 0 & \partial_{k\!+\!1}^{i,j}\\[3pt]
        0 & \left(\partial_{k\!+\!1}^{i,j}\right)^* & 0
    \end{pmatrix} -  
    \xi\begin{pmatrix}
        P^i_{k\!-\!1} & 0 & 0\\[3pt]
        0 & -P^i_{k} & 0\\[3pt]
        0 & 0 & P^j_{k\!+\!1}
    \end{pmatrix}
\end{align}
with $\xi\in\RR$, and its square is given by 
\begin{align}
    \left(D_k^{i,j}\right)^2 = \begin{pmatrix}
        \partial_{k}^{i,i}\left(\partial_{k}^{i,i}\right)^* & 0 & \partial_{k}^{i,i}\partial_{k\!+\!1}^{i,j}\\[3pt]
        0 & L_k^{i,j} & 0\\[3pt]
        \left(\partial_{k\!+\!1}^{i,j}\right)^*\left(\partial_{k}^{i,i}\right)^* & 0 & \left(\partial_{k\!+\!1}^{i,j}\right)^*\partial_{k\!+\!1}^{i,j}
    \end{pmatrix} +  
    \xi^2\begin{pmatrix}
        P^i_{k\!-\!1} & 0 & 0\\[3pt]
        0 & P^i_{k} & 0\\[3pt]
        0 & 0 & P^j_{k\!+\!1}
    \end{pmatrix}
\end{align}
Based on its definition, the square of the persistent Dirac does not lead to a diagonal block matrix as in the standard case, but just one of the blocks provides the $k$-th persistent combinatorial Laplacian shifted by $\xi^2$. It follows that the eigenvalues of the $k$-th persistent combinatorial Laplacian correspond to a $\xi^2$ shift of the squared eigenvalues of the $k$-th Dirac operator. In the computation of persistent Betti numbers, the shift $\xi$ is chosen to be the positive eigenvalue of the Dirac operator. As noted in \cite{ameneyro2024quantum}, setting $\xi=0$ results in overcounting of the states in the kernel of the Laplacian and thus an incorrect estimate of Betti numbers. 

The persistent Dirac operator has been applied for analyzing time series data \cite{ameneyro2023quantum}. Another quantum implementation of computing the persistent Betti numbers was proposed by Hayakawa \cite{hayakawa2022quantum}, which is based on the block-encoding of persistent combinatorial Laplacians. While these methods provide efficient computations of the persistent topological features, their realization requires fault tolerant quantum computers, which are unlikely to be available in the near future for realistic big data applications.

\subsubsection{Other variations}

Here we provide a few variations of the persistent Dirac operator recently developed for data with different structures. 

\paragraph{Persistent path Dirac and hypergraph Dirac operators}
The path Dirac and hypergraph Dirac operators, along with their persistence extensions, were introduced in \cite{suwayyid2024persistentB} for extracting valuable information from path complexes and hypergraphs. Their application to molecular structure analysis has also been explored. 

\paragraph{Persistent Mayer Dirac operators}
The Mayer Dirac and persistent Mayer Dirac operators on $N$-chain complexes were proposed in \cite{suwayyid2024persistent} and applied to molecular structure analysis. Note that Mayer homology \cite{mayer1942new} differs from the standard homology. In this context, the boundary operator satisfies $\partial^N = 0$ with $N\ge 2$ instead of the usual condition $\partial^2 = 0$. For a description of tools based on Mayer topology, see Section~\ref{sec:algebraicTopo.Mayer}.

\paragraph{Khovanov Dirac}
The Khovanov Dirac \cite{jones2025khovanov} was proposed for analyzing the global topological features of knots and links. It is based on Khovanov homology, a powerful topological invariant in geometric topology for studying Knots and links. This topic will be discussed further later in Section~\ref{sec:geometricTopo.Khovanov}.

\subsection{Other formulations}

In this section, we outline formulations recently developed in topological data analysis that utilize alternative homology theories, extending the standard homology for analyzing point-cloud data.

\subsubsection{Persistent sheaf theory}

The aforementioned methods, such as persistent homology and persistent combinatorial Laplacians, have proven to be powerful for analyzing the topological structure of data \cite{qiu2023persistent}. However, they remain limited in their ability to describe the local features along with the relationships among them, which often play a crucial role in practical applications.

Sheaf theory provides a more powerful and general framework for understanding topological spaces. In particular, it enables the study of local-to-global properties of structures defined on spaces. In recent years, it has found various applications in topological data analysis \cite{ghrist2014elementary, yoon2018cellular, hansen2019toward, hansen2020laplacians, cooperband2023cosheaf, cooperband2025unified}, particularly through the notion of cellular sheaves \cite{shepard1985cellular}, which is more accessible to researchers and easier to implement for practical computations. For a simplicial complex, a cellular sheaf can be thought of as an assignment of a vector space to each simplex, together with linear transformations, called restriction maps, defined between the vector spaces by face relations. These restriction maps must satisfy certain consistency conditions, allowing global information to be determined from local features. This makes cellular sheaves valuable for analyzing data with localized properties, from which global features can be inferred.

Due to the direction of the restriction maps, sheaves are naturally associated with cohomology. The consistency conditions lead to a well-defined coboundary operator $\delta$ on sheaf cochain complexes satisfying $\delta\delta=0$, which defines the sheaf cohomology, with the standard cohomology being a special case. A sheaf-theoretical perspective on persistent homology was first investigated in \cite{curry2014sheaves} and later extended to the persistent setting in \cite{yegnesh2016persistence, russold2022persistent}. A notion of persistent sheaves was also proposed in \cite{hang2023correspondence}. 

By introducing an inner product on the cochain groups, the sheaf Laplacians \cite{hansen2019toward} can be defined. Their kernels correspond to the sheaf cohomology groups with dimensions given by the Betti numbers. By incorporating persistence, Wei and Wei \cite{wei2025persistentsheaf} proposed persistent sheaf Laplacians to address the limitation of persistent combinatorial Laplacians in handling data points with labels or weights. The framework generalizes persistent combinatorial Laplacians, as the latter can be seen as a special case of the former when all data points are assigned the same quantity. Given a pair of sheaves in a filtration, the persistent boundary operator is defined analogously to the case of persistent combinatorial Laplacians. The main difference lies in the underlying structure of the data and the use of cohomology instead of homology, along with a coboundary operator defined by scalar multiplication subject to consistency rules. The persistent sheaf Laplacians follow with the same formulas, with kernels isomorphic to the persistent cohomology groups and dimensions given by the persistent Betti numbers. In addition, the matrix representations of persistent sheaf Laplacians can be calculated in the same way as for the persistent combinatorial Laplacians \cite{wei2025persistentsheaf}.

\subsubsection{Persistent Mayer topology}
\label{sec:algebraicTopo.Mayer}

Note that all aforementioned methods in topological data analysis rely on the construction of the standard chain complexes, where the boundary operator satisfies $\partial^2 = 0$. These methods might be limited to revealing more intricate information about the underlying structure of data as the boundary operator captures only the relationship between simplices of consecutive dimensions.

To address this, Shen et al. \cite{shen2024persistent} introduced the persistent Mayer homology and persistent Mayer Laplacians, which generalize the Mayer homology \cite{mayer1942new} to the persistent setting, allowing the study of connections between simplices across multiple dimensions rather than just adjacent ones. The approach is built on an $N$-chain complex with the boundary operator satisfying a more general condition $\partial^N = 0$ with $N\ge 2$.   
This extends the current existing frameworks based on standard homology where $\partial^2 = 0$. For a discussion about the relationship between Mayer homology and the standard homology, see \cite{spanier1949mayer}.

Analogous to the standard case, the kernels of the persistent Mayer Laplacians correspond to the persistent Mayer homology groups, with ranks given by the persistent Mayer Betti numbers. However, in contrast to the standard persistent homology governed solely by the filtration parameter, the persistent Mayer homology and persistent Mayer Laplacians for each degree $k$ are controlled by an additional parameter. The framework thus enables the encoding of richer topological and geometric or combinatorial information about the underlying structure, in particular, the relationship between simplices across different dimensions. A further generalization of persistent Mayer Laplacians, called persistent Mayer Dirac, has been proposed in \cite{suwayyid2024persistent}, which is more sensitive to local features.

Persistent Mayer Laplacians and the Mayer Dirac have been applied to the study of molecular structure analysis \cite{shen2024persistent, suwayyid2024persistent}. Persistent Mayer homology has been used as multichannel features for predicting protein-ligand binding affinities \cite{feng2024mayer}.

\subsubsection{Persistent interaction topology}
In algebraic geometry, intersection theory studies the intersection of subvarieties within a given variety. In topology, intersection (co)homology studies singular spaces. In the discrete setting, interaction involves understanding how components within a complex system influence each other. This concept is particularly important in fields such as complex networks, molecular interactions, and organism interactions in biological systems.
Traditional tools in topological data analysis often lack the ability to capture such intersections, as they typically focus on global topological features instead of focusing on localized interactions between elements in the system. 

Recently, motivated by the notion of interaction cohomology \cite{knill2018cohomology},
Liu et al. \cite{liu2023interaction} developed the interaction topology to gain insight into the interactions between topological spaces and characterize such interactions using topological invariants. The theory is built on the category of interaction spaces. In the paper, they introduced the concept of interaction spaces, studied their homotopy and singular homology, and proved that interaction singular homology is a homotopy invariant. In addition, they proposed the interaction simplicial homology and cohomology to facilitate applications in point-cloud data analysis.

Based on interaction homology, Liu et al. \cite{liu2024persistent} proposed the persistent interaction homology and persistent interaction Laplacians, and demonstrate their stability. The paper also provides an algorithm for computing the interaction homology and interaction Laplacians for finite point sets. The strength of the framework is demonstrated through the analysis of interactions among different types of atoms within molecules.


\section{Differential topology approaches }\label{secDifferentialT}

In this section, we shift our attention to differential topology analysis of data on differential manifolds. The approach uses tools from differential topology, such as vector fields, differential forms, Laplace-Beltrami operators, and Hodge Laplacians, along with their extensions to multiscale analysis, to understand the underlying structure of differential manifolds. These tools are deeply connected to the topology of the underlying manifolds through de Rham-Hodge theory, a landmark theory of the 20th century mathematics, and thus are well-suited for analyzing data naturally represented as differential manifolds. Examples include 2D surfaces, 3D solids, and image data such as 2D slices or 3D volumes from magnetic resonance imaging (MRI) or computed tomography (CT) scans. These data are commonly found in fields such as computer graphics, biological modeling, and medical imaging.

Although in many scenarios, techniques in algebraic topology analysis for point cloud data can be applied to data using point cloud representation, it is more appropriate to model the data as continuous or differential manifolds. These data are naturally continuous in observation, even though they are practically represented as discrete points. Relying solely on sampled points can lead to inconsistencies due to their varying combinatorial representations, resulting in different outcomes using algebraic and combinatorial techniques. A more natural solution is to model the data theoretically as continuous or differential manifolds, and perform discretization only in the implementation. This way ensures that the analysis remains consistent and stable regardless of different discretizations, and the results are convergent.

In addition, the study of differential manifolds naturally incorporates the boundary conditions, such as the normal (Dirichlet) and tangential (Neumann) boundary conditions, which are significantly important in applications where the boundary carries physical or geometric meaning, such as electromagnetism and fluid dynamics. These boundary conditions provide constraints that govern the behavior of a system near the boundary, and ensure the existence and uniqueness of solutions when solving boundary value problems. It also ensures the well-posedness of differential operators and Hodge Laplacians' capability to perform differential calculus on manifolds with boundary. However, there is no intrinsic notion of domain boundary for combinatorial Laplacians.

The implementation of de Rham Hodge theory has been challenging for years due to the complicated boundary conditions. Earlier methods, such as the standard finite element method \cite{ciarlet2002finite}, are insufficient to capture the manifold topology and can lead to entirely wrong eigenvalues for the vector Laplacians under normal or tangential boundary conditions \cite{arnold2018finite}. With recent advances, such as the finite element exterior calculus (FEEC) \cite{arnold2006finite} and the discrete exterior calculus (DEC) \cite{desbrun2006discrete}, the practical realization of the theory has become possible. These techniques ensure that the discretizations are consistent with the continuous theory and thus allow for accurate computation of the cohomologies of the underlying manifolds. Based on these techniques, computational frameworks have been developed for computing Hodge decompositions, Hodge Laplacians and their multiscale extension. These developments have enabled the application of de Rham Hodge theory in topological data analysis, opening new directions for analyzing the topological and geometric features of data on manifolds.

Below, we start with a review of de Rham-Hodge theory for manifolds with boundary, 
including the de Rham cohomology, Hodge Laplacians, and Hodge decomposition, along with their recent implementations. We then provide a short discussion on persistent de Rham cohomology \cite{chen2021evolutionary}. We also focus on the recently developed persistent Hodge Laplacians \cite{chen2021evolutionary, su2024persistent} in capturing the topological and geometric features of 2D and 3D compact domains. We include both the Lagrangian formulation \cite{chen2021evolutionary} and the Eulerian formulation \cite{su2024persistent} of the de Rham-Hodge theory, depending on whether the underlying manifolds are discretized using simplicial meshes or represented as regions in a Cartesian grid defined by level-set functions, respectively.

\subsection{De Rham-Hodge theory for manifolds with boundary}

The de Rham-Hodge theory offers a foundational framework for understanding the topological and geometrical structures of manifolds through differential forms. Its formulation is particularly subtle in the case of manifolds with boundary, where appropriate boundary conditions, such as normal and tangential boundary conditions \cite{shonkwiler2009poincare}, must be imposed. These boundary conditions are essential for the core results of the theory, i.e., the correspondence between harmonic forms and the cohomology.

The notion of differential forms might be abstract at first. However, they are closely related to familiar concepts in vector calculus. For example, in $\RR^3$, a $0$- or $3$-form is just a scalar field, while a $1$- or $2$-form corresponds to a vector field. In general, differential forms can be viewed as a generalization of scalar and vector fields on differential manifolds. These correspondences between differential forms and scalar or vector fields provide an intuitive understanding of the boundary conditions in the case of $1$-forms: a $1$-form is normal if its corresponding vector field is normal to the boundary, or tangential if the field is tangential to the boundary. The formal definitions of normal and tangential boundary conditions, and also a complete correspondence between differential forms and their counterparts in vector calculus in $\RR^3$ can be found in \cite{zhao20193d}.
The two types of boundary conditions are dual to each other under the Hodge star operator.

The differential $d^k$ (i.e., the exterior derivative) and the codifferential operator $\delta^k$, are two fundamental operators in building the de Rham-Hodge theory. We briefly recall them here with their correspondences with classical operators in vector calculus for an intuitive understanding. The differential $d^k$ maps $k$-forms to $(k\!+\!1)$-forms and satisfies $d^kd^{k-1}=0$. It generalizes and unifies the notions of gradient $\nabla$, curl $\nabla\times$, and divergence $\nabla\cdot$ in vector analysis, which correspond to $d$ applied to $0$-, $1$-, and $2$-forms in $\RR^3$, respectively. The codifferential $\delta^k$ can be defined as an adjoint operator of $d^{k-1}$ with respect to the Hodge $L^2$ inner product induced from the Hodge star operator on the space of differential forms. It maps $k$-forms to $(k\!-\!1)$-forms and satisfies $\delta^{k-1}\delta^k = 0$. In $\RR^3$, it corresponds to $-\nabla$, $\nabla\times$, and $-\nabla\cdot$ when applied to $2$-, $1$-, and $0$-forms, respectively.

\subsubsection{De Rham cohomology}

Let $\Omega_n^k$ and $\Omega_t^k$ be the space of normal $k$-forms and the space of tangential $k$-forms over a smooth manifold $M$ of dimension $m$, respectively. For their formal definitions, see \cite{zhao20193d}. These two spaces are dual to each other under the Hodge star operator, i.e., $\Omega_n^k\cong \Omega_t^{m-k}$. The de Rham cohomologies of differential forms are defined on two sequences connected by the differential $d^k$ and the codifferential $\delta^k$:
\begin{align}
\cdots\overset{d^{k-2}}{\longrightarrow} \Omega_n^{k-1} \overset{d^{k-1}}{\longrightarrow} \Omega_n^k \overset{d^{k}}{\longrightarrow} \Omega_n^{k+1} \overset{d^{k+1}}{\longrightarrow} \cdots,\\
\cdots\overset{\delta^{k-1}}{\longleftarrow} \Omega_t^{k-1} \overset{\delta^{k}}{\longleftarrow} \Omega_t^k \overset{\delta^{k+1}}{\longleftarrow} \Omega_n^{k+1} \overset{\delta^{k+2}}{\longleftarrow} \cdots.
\end{align}
These two complexes are well-defined as $d$ preserves the normal boundary conditions and $\delta$ preserves the tangential boundary conditions. The resulting cohomologies correspond to the $k$-th relative de Rham cohomology $H_{dR}^k(M, \partial M)$ and the $k$-th absolute de Rham cohomology $H_{dR}^k(M)$ as follows:
\begin{align}
    H_{dR}^k(M, \partial M) &= \ker d^k/\img d^{k-1},\\
    H_{dR}^k(M) &= \ker \delta^{k}/\img \delta^{k+1},
\end{align}
where $\ker d^k$ (resp. $\ker\delta^{k}$) is called the space of closed (resp. coclosed) $k$-forms and $\img d^{k-1}$ (resp. $\img\delta^{k+1}$) is the space of exact (resp. coexact) $k$-forms. 
The ranks of $H_{dR}^k(M, \partial M)$ and $H_{dR}^k(M)$ are given by the $(m\!-\!k)$-th Betti number $\beta_{m-k}$ and the $k$-th Betti number $\beta_k$, respectively. Here the domains of these operators are restricted under the corresponding boundary conditions. 
For a detailed discussion of the absolute and relative de Rham cohomology and the corresponding isomorphisms, refer to \cite{schwarz2006hodge}.

\subsubsection{Hodge Laplacians}
\label{sec.differentialTopo.HodgeLap}

The $k$-th Hodge Laplacian $\Delta^k$ on the space of $k$-forms is defined as follows:
\begin{align}
    \Delta^k = d^{k-1}\delta^{k} + \delta^{k+1} d^k.
\end{align}
When $k=0$, the Hodge Laplacian reduces to the Laplace-Beltrami operator $\delta^1d^0$ acting on the space of functions defined on the manifold. In the presence of the boundary, the Hodge Laplacian $\Delta^k$ has no direct correspondence with the manifold topology, as its kernel is infinite-dimensional. However, one can define the Hodge Laplacians $\Delta^k_n$ and $\Delta^k_t$ by restricting the Hodge Laplacian $\Delta^k$ to the space of normal forms and to the space of tangential forms, respectively. These two types of Hodge Laplacians lead to finite-dimensional kernels. The space $\ker\Delta^k_n$ corresponds to the relative de Rham cohomology, and $\ker\Delta^k_t$ to the absolute de Rham cohomology \cite{friedrichs1955differential}. In addition, these two kernels are related by the Poincar{\'e} duality, leading to the following identifications:
\begin{align}
    \ker\Delta_n^k &\cong H^k_{dR}(M, \partial M)\cong H^{m-k}_{dR}(M)\cong\ker\Delta_t^{m-k}
\end{align}
The elements in $\ker\Delta_n^k$ are called the normal harmonic $k$-fields, and those in $\ker\Delta_t^k$ are called the tangential harmonic $k$-fields. These identifications enable the use of either type of boundary conditions for the study of the topology of the underlying manifold and its boundary.

Note that similar to the combinatorial case, the harmonic spectra of Hodge Laplacians capture the topological information about the underlying manifolds with the number of zero-eigenvalues given by the Betti numbers. However, the non-harmonic spectra of the Hodge Laplacians encode rich geometric information such as the metric and the curvature of the underlying manifolds. This is fundamentally different from the combinatorial case, where the non-harmonic spectra reflect only the combinatorial or the discrete structural information rather than the smooth geometric features. For a thorough discussion about the similarity and difference between the combinatorial Laplacians and Hodge Laplacians, see \cite{ribando2024graph}.

The implementation of Hodge Laplacians has been realized through FEEC in \cite{arnold2006finite} and through DEC \cite{zhao20193d, ribando2024graph}, with both approaches providing accurate eigenvalue results and preserving the cohomology of the underlying domain. In addition, a notion of boundary-induced graph (BIG) Laplacians has been introduced \cite{ribando2024graph}, defined on Cartesian domains with appropriate boundary conditions. It is a special case of the discretized Hodge Laplacians and is simpler to compute while delivering correct topological dimensions of the data. The Hodge Laplacian of degree $0$, i.e., the Laplace-Beltrami operator, has been widely applied to the study of the shapes of manifolds. For example, ShapeDNA \cite{reuter2006laplace} provides a simple representation of a shape using the spectra of the operator. Global point signature (GPS) \cite{rustamov2007laplace} encodes a shape using a vector of scaled eigenfunctions of the operator. ShapeDNA utilizes standard finite element approaches, while GPS employs tools in DEC.

\subsubsection{Hodge decomposition}

The Hodge decomposition is a central result in Hodge theory, which decomposes the space of differential $ k$-forms into a direct sum of mutually orthogonal subspaces. A $3$-component Hodge decomposition of differential forms, called the Hodge Morrey decomposition \cite{morrey1956variational}, is given as follows:
\begin{align}
    \Omega^k = d\Omega_n^{k-1} \oplus \delta\Omega_t^{k+1} \oplus\mathcal{H}^k
\end{align}
Here $\mathcal{H}^k = \ker d^k\cap\ker\delta^{k}\subset\ker\Delta^k$, known as the space of harmonic $k$-fields. When applied to vector fields on a domain, this decomposition states that any vector field can be orthogonally decomposed into a curl-free, a divergence-free, and a harmonic component, leading to the classical Helmholtz-Hodge decomposition \cite{ladyzhenskaya1969mathematical}. The space $\mathcal{H}^k$ is equal to the kernel of the Hodge Laplacian $\ker\Delta^k$ and isomorphic to the $k$-th de Rham cohomology group. However, in the presence of a boundary, the space $\mathcal{H}^k$ is infinite-dimensional \cite{schwarz2006hodge}, providing no connections to the manifold topology. It has been shown that this space can be further split by enforcing boundary conditions, revealing the relative and absolute cohomology \cite{friedrichs1955differential}.
In particular, there is a $5$-component topology-preserving Hodge decomposition for compact domains in the Euclidean spaces \cite{shonkwiler2009poincare}, given as follows:
\begin{align}\label{eq.hd.5subspaces}
	\Omega^k = d\Omega^{k-1}_n\oplus\delta\Omega^{k+1}_t \oplus \mathcal{H}^{k}_n \oplus \mathcal{H}^{k}_t\oplus (d\Omega^{k-1}\cap\delta\Omega^{k+1}),
\end{align}
where $\mathcal{H}^k_n=\ker\Delta^k_n$ and $\mathcal{H}^k_t=\ker\Delta^k_t$ are the restrictions of $\mathcal{H}^k$ to the space of normal forms and the space of tangential forms, respectively, encoding the relative and absolute cohomology of the underlying manifold as discussed previously.

For the computational frameworks of Hodge decomposition, we focus here only on the complete topology-preserving $5$-component Hodge decomposition. A review of the classical $3$-component Helmholtz-Hodge decomposition is available in \cite{bhatia2012helmholtz}. The discretization of the 5-component Hodge decomposition of vector fields has been proposed using their differential form representations for tetrahedral meshes \cite{zhao20193d} and for 2D and 3D compact domains defined by level set functions in Cartesian grids \cite{su2024hodge, su2024topology}. Both methods are based on DEC. Other approaches, assuming piecewise-constant vector fields, have been developed for surface triangle meshes \cite{poelke2016boundary} and also tetrahedral meshes \cite{poelke2017hodge}. A discussion of this $5$-component Hodge decomposition for compact domains in $\RR^3$ in terms of vector and scalar fields can be found in \cite{cantarella2002vector}.

\subsection{Persistent de Rham cohomology}

Persistent de Rham cohomology was first introduced by Chen et al. in \cite{chen2021evolutionary}. The work provides a detailed discussion about the persistence of harmonic forms under normal and tangential boundary conditions. The persistent de Rham cohomology holds the same importance both in theory and practical applications as persistent homology does for point-cloud data. It shows promise as a foundation tool and opens new directions for analyzing the topological structures of data that can be represented as differential manifolds.

To define the persistence de Rham cohomology, we consider a manifold $M$, defined by a sublevel set of a level set function $f$, i.e., $M = f^{-1}([-\infty, c])$ for $c\in\RR$. By varying the isovalues of $f$, a filtration of differential manifolds with inclusion maps can be obtained as follows: 
\begin{align}
    M_0\hookrightarrow M_1 \cdots \hookrightarrow M_{s-1}\hookrightarrow M_s,
\end{align}
where each $M_i$ is the sublevel set of $f$ corresponding to a different isovalue. The inclusion map for any pair of manifolds $M_i\hookrightarrow M_j$ in the sequence induces a unique extension on the normal differential forms $\Omega^k_n(M_i)\to\Omega^k_n(M_j)$ satisfying certain conditions. This map further induces a map on the cohomology groups $\mathcal{\psi}^k_{i, j}: H^k_{dR, i}(M, \partial M)\to H^k_{dR, j}(M, \partial M)$ that tracks the evolution of topological features, namely, normal harmonic fields as $H^k_{dR, i}(M, \partial M)\cong\mathcal{H}^k_{n, i}$. Here $\mathcal{H}^k_{n, i}$ denotes the space of normal harmonic $k$-fields on $M_i$. One can then define the persistent relative de Rham cohomology as the image of the map:
\begin{align}
    H^k_{dR,i,j}(M, \partial M) = \img \mathcal{\psi}^k_{i,j} = \ker d^k_i/\left(\img d^{k-1}_j\cap\ker d^k_i\right).
\end{align}
Its rank provides the $(m\!-\!k)$-th persistent Betti number $\beta_{m-k}^{i, j}$. The birth and death, as well as the persistence of normal harmonic fields, can then be defined correspondingly by examining the elements in the persistent cohomology groups, analogous to the combinatorial case. Note that only the spaces of normal differential forms are considered here. The dual formulation under the tangential boundary conditions can be constructed in a similar manner \cite{chen2021evolutionary}.

\subsection{Persistent Hodge Laplacians}

In the same paper \cite{chen2021evolutionary}, Chen et al. first theoretically proposed the persistent Hodge Laplacians. The paper also presents an implementation of multiscaled Hodge Laplacians using DEC for manifolds in the Lagrangian formulation, i.e., manifolds discretized using simplicial meshes. Subsequently, a computational framework of persistent Hodge Laplacians was proposed by Su et al. in \cite{su2024persistent} in the Eulerian setting with the adapted DEC from simplicial meshes to regular grids. The manifolds in this context are modeled as domains in a Cartesian grid bounded by isocurves or isosurfaces of level-set functions. These developments are based on combining de Rham cohomology with the evolution of differential manifolds in a filtration, and can be viewed as extensions of earlier work on persistent homology in the cubic setting \cite{wagner2011efficient, wang2016object}. 

Given any pair of manifolds $M_i\hookrightarrow M_j$ in the filtration aforementioned in the previous subsection, the $k$-th persistent Hodge Laplacian on the space of normal $k$-forms $\Omega^k(M_i)$ is given as follows:
\begin{align}
    \Delta^k_{i,j} = d^{k-1}_{i,j}\delta^{k}_{i,j} + \delta^{k+1}_i d^k_i,
\end{align}
where $d^k_i$ and $\delta^k$ is the $k$-th differential and codifferential on normal forms on $M_i$, respectively, and $d^{k-1}_{i,j}$ is the persistent differential operator from normal $(k\!-\!1)$-forms on $M_j$ to normal $k$-forms on $M_i$, with $\delta^{k}_{i,j}$ its adjoint operator. 
The persistent Hodge Laplacians can similarly be defined on the space of tangential forms. However, they are equivalent due to the duality of the normal and tangential boundary conditions. When $i = j$, the persistent Hodge Laplacian reduces to the usual Hodge Laplacian restricted to normal forms. The kernel of $\Delta^k_{i,j}$ defines the persistent harmonic fields. It is isomorphic to the persistent relative de Rham cohomology, and its rank leads to the persistent Betti number.

We now turn to the discretization frameworks and briefly discuss the Lagrangian and Eulerian formulations using DEC, depending on how the underlying manifolds are discretized.

\subsubsection{Lagrangian formulation}

In the Lagrangian formulation, the underlying manifolds are discretized as simplicial meshes, specifically triangular meshes in 2D or tetrahedral meshes in 3D. A discrete $k$-form is just a vector with each entry assigned to a $k$-simplex. For example, a discrete $0$-form is a vector with values assigned to the vertices, while a $1$-form has values corresponding to edges. The differential is a sparse matrix given as the transpose of the signed incidence matrix between simplices of consecutive dimensions. The codifferential is the product of the incident matrix and the discrete Hodge star operators, which takes the form of a diagonal matrix. For further details of these discrete operators, see \cite{zhao20193d, chen2021evolutionary}.

The Lagrangian formulation requires explicit tessellation of the domains, which relies on high-quality meshing tools such as CGAL \cite{fabri2009cgal}, Gmsh \cite{geuzaine2009gmsh}, Pymesh \cite{zhou2019pymesh}, etc. This makes the implementation of persistent Hodge Laplacians challenging in this context, as generating consistent meshes throughout the filtration can be problematic, and the resulting discrete operators may not be sparse matrices. The multiscale version of Hodge Laplacians \cite{chen2021evolutionary} in the Lagrangian formulation, in addition to its dependence on these meshing tools, has several limitations. It is sensitive and numerically inconsistent because of variations in tetrahedral mesh representations of a compact domain. The method is also computationally expensive, as it requires repeated tessellation of the domains to create the evolution of manifolds in a filtration. These limitations make the framework unsuitable for machine learning tasks, which require feature representation computed on the same scale for comparability, and also efficiency, which is essential for computing features for large datasets.

\subsubsection{Eulerian formulation}

The Eulerian formulation of persistent Hodge Laplacians \cite{su2024persistent} has been introduced with the adapted DEC to address the limitations in the Lagrangian case. A computational algorithm for persistent Hodge Laplacians is presented in the paper. The underlying manifold in this setting is modeled as a sublevel set of a level-set function on a standard Cartesian grid, i.e., a domain bounded by an isocurve or an isosurface of a level set function. The adapted DEC from simplicial meshes to regular grids follows a similar structure as in the Lagrangian case. Note that in the Eulerian formulation, adjusting the isovalues of a level-set function naturally leads to a nested sequence of grid complexes in a filtration, making the construction of discrete persistent Hodge Laplacians feasible.

The Eulerian representation offers several advantages: First, it eliminates the need for explicit tessellation, thereby removing the reliance on high-quality meshing tools. Second, the computation of Laplacians for all objects can be performed on a pre-designed Cartesian grid, the extracted features from different objects are thus consistent and comparable for machine learning tasks. In addition, as the vertices, edges, faces, and cells are fixed in a Cartesian grid, the data structure and the discrete operators are significantly simplified. In particular, the evolution of manifolds can be easily obtained by adjusting the isovalues of the level-set functions. These make the computation of Laplacians more efficient in the Eulerian formulation. The framework was demonstrated in the same paper by integrating it into machine learning models for predicting the protein-ligand binding affinities.

\section{Geometric topology approaches }
\label{secGeometricT}

This section focuses on methods grounded in geometric topology, in particular knot theory, for analyzing curves embedded in 3-space, such as knots, links, tangles, and linkoids. Such data structures are often found in chemical compounds \cite{panagiotou2019topological}, DNA \cite{arsuaga2005dna}, and protein structures \cite{sulkowska2012conservation}. Knot theory studies these objects up to ambient isotopy, informally referring to a continuous deformation that does not involve cutting and rejoining. The key tools in Knot theory are knot invariants, which are quantities that remain unchanged under these deformations and play an essential role in distinguishing knots and links. Examples include the crossing number \cite{adams2004knot}, the bridge number \cite{schubert1954numerische}, the knot group \cite{crowell2012introduction}, the Jones polynomial \cite{jones1997polynomial}, the Alexander polynomial \cite{alexander1928topological},  the Khovanov homology \cite{khovanov2000categorification}, and the knot Floer homology \cite{ozsvath2004holomorphic}.

Knot theory has been applied in various fields, such as physics \cite{ohtsuki2001quantum}, chemistry \cite{liang1994knots}, and biology \cite{sumners2020role,schlick2021knot,millett2013identifying}, but its success in practical applications has been somewhat limited. One of the challenges lies in the fact that classical tools in knot theory are inherently global and topological, offering little insight into local geometric and topological features. This makes these tools insufficient for analyzing biological structures, such as DNA and proteins, which often carry intricate local structures and/or local entanglement between segments. In addition, many datasets encountered in real applications consist of open curves rather than closed loops, which makes classical knot invariants inapplicable. These challenges highlight the need to develop more effective tools for analyzing curves in practical settings.

Motivated by the success of topological data analysis for various domains, the concept of knot data analysis (KDA) was introduced in \cite{shen2024knot}. The paper proposes a framework called multiscale Gauss linking integral (mGLI), which is the first to combine multiscale analysis with tools in knot theory. Following the work, several other geometric topological frameworks have been proposed to address the aforementioned limitations of classical methods in analyzing curves in 3-space. These include a localized Jones polynomial \cite{song2025multiscale}, persistent Khovanov homology \cite{shen2024evolutionary}, and Khovanov Laplacians, and Khovanov Dirac operator \cite{jones2025khovanov}. 

In the following, we explore these recently developed geometric topology approaches. The multiscale Gauss linking integral (mGLI) \cite{shen2024knot} and the localized Jones polynomials \cite{song2025multiscale} build upon the Gauss linking integral \cite{Gauss1877} and the Jones polynomial of open-ended knot diagrams \cite{panagiotou2020knot}, respectively, which are currently the only two known continuous measures of single or pairwise topological complexity for open curves in 3-space. These approaches enable the analysis of local and global entanglement properties across multiple spatial scales. The persistent Khovanov homology \cite{shen2024evolutionary} extends the Khovanov homology by introducing filtrations of link diagrams based on smoothing transformations, revealing finer topological and geometric features that are often overlooked by traditional invariants. The Khovanov Laplacians and Dirac operator \cite{jones2025khovanov} extend the Khovanov homology by incorporating spectral information, offering new descriptors that recover known invariants from homology while encoding additional non-harmonic information.

\subsection{Multiscale Gauss linking integral}

 Shen et al. \cite{shen2024knot} proposed the multiscale Gauss linking integral (mGLI), offering a multiscale approach for analyzing the local structure, connectivity, and entanglement of segments in both open and closed curves. These local features are crucial for understanding the physical properties and functions of data, such as chemical compounds and biological molecules. Although local, the method can also recover the global topological features when sufficiently large scales are applied. The framework builds on the Gauss linking integral, which is a double line integral that provides the linking number of two curves and has been successfully applied in understanding polymers \cite{panagiotou2019topological} and proteins \cite{panagiotou2020topological,baldwin2021local,baldwin2022local}.

The method requires curve segmentation to partition curves into small segments for a given pair of curves. The pairwise Gauss linking integrals between these curve segments then result in a segmentation matrix whose entries encode the pairwise linking information. For multiscale analysis, weights on these entries, defined by indicator functions in a specific distance range from each segment, are then incorporated. By varying the distance range, this process generates a collection of segmentation matrices that capture the local entanglement across various spatial scales. The cumulative integrals across scales also provide additional insights into local structures, and reveal the global relationships at sufficiently large scales. For each segment in one curve, the local linking information with all segments of the other curve, calculated within each distance range, can be used to generate the multiscale feature of that segment through statistical methods. The construction can also be applied to a single curve by considering the local linking information between segments of the same curve.

The segmentation of curves and the distance range are crucial in determining the effectiveness of  Gauss linking integrals for feature extraction. Therefore, the segmentation strategy and the distance range should be chosen carefully based on specific applications. The framework has been integrated with machine learning models to address a variety of complex biological problems, including protein $B$-factor prediction, protein-ligand binding affinity prediction, and others. 
The model consistently outperforms the traditional TDA models and demonstrates superior performance in certain tasks \cite{shen2024knot}.

By combining the classical Gauss linking integral with multiscale analysis, mGLI offers a powerful framework for analyzing the local and global structures of curve-type data, opening a new area in data analysis and knot learning. The framework addresses the limitations of traditional TDA techniques by examining the local entanglement of curve segments, and it can be applied to open curves. In addition, it is computationally efficient, requiring only a few minutes on a personal computer to generate the features, which makes it practical for analyzing large datasets. The method shows promise and holds strong potential for practical applications across various fields.

\subsection{Localized Jones polynomials}

The localized Jones polynomials were introduced by Song et al. \cite{song2025multiscale} for analyzing the local and global entanglement properties of collections of disjoint open or closed curves. The paper proposes two models, termed the multiscale Jones polynomial and the persistent Jones polynomial, which extend the Jones polynomial proposed in \cite{barkataki2022jones} by incorporating spatial multiscale analysis. The approach in \cite{barkataki2022jones} generalizes the traditional Jones polynomial \cite{jones1997polynomial} for closed curves and the one in \cite{panagiotou2020knot} for an open curve, offering a continuous measure of linking complexity for collections of open curves in 3-space. 

The stability of both proposed localized Jones polynomials has been analyzed in the same paper, ensuring their robustness for practical use. These models have also been applied to protein structure analysis, providing insight into the entanglement that goes beyond atomic
positions alone. The frameworks demonstrate potential to benefit a range of practical applications.

\paragraph{Multiscale Jones polynomial}
The multiscale Jones polynomial \cite{song2025multiscale} follows a similar strategy used for the multiscale Gauss linking integral \cite{shen2024knot}. Given a collection of disjoint open or closed curves with a segmentation, the method computes, for each curve segment, the Jones polynomial of the set consisting of that curve segment together with all curve segments that lie within a specified distance range from that curve segment. By varying the distance interval, the process generates a vector of Jones polynomials for each segment, which captures the topological entanglement at multiple spacial scales. By considering all curve segments, the method finally leads to a matrix of Jones polynomials, encoding both the local and global entanglement properties for the collection of curves. These Jones polynomials are then evaluated at a fixed parameter for real applications, resulting in a real characteristic matrix suitable for downstream analysis. Analogous to the multiscale Gauss linking integral \cite{shen2024knot}, the featurization depends on the choice of segmentation of the collections of curves, and the selection of distance ranges, both of which can be tuned appropriately for specific applications.

\paragraph{Persistent Jones polynomial}

The persistent Jones polynomial was proposed by considering the pairwise distance matrix of curve segments for a collection of disjoint open or closed curves with a segmentation, which then induces a filtration of Vietoris-Rips complexes (other types of simplicial complexes may also be used), with the vertices corresponding to the curve segments \cite{song2025multiscale}. Each facet in the filtration is assigned a weight given by the Jones polynomial of the set of curve segments forming that facet, leading to a weighted filtration of polynomials of the segmentation. This filtration gives rise to weighted barcodes, which encode the birth and death of facets along with their associated polynomial weights. For practical applications, the Jones polynomials are evaluated at a fixed parameter value. By introducing a filtration, the persistent Jones polynomial encodes the local and global interactions among curve segments within a segmentation. In contrast to the persistent homology for point cloud data, which focuses on the behavior of holes, the persistent Jones polynomial captures the entanglement of curves in $3$-space.

\subsection{Persistent Khovanov homology}
\label{sec:geometricTopo.Khovanov}

Khovanov homology \cite{khovanov2000categorification} is one of the most significant advances in knot theory, providing a powerful link invariant that categorifies the Jones polynomial \cite{jones1997polynomial}. It reveals richer information such as torsion. It can also detect the unknot, which is a very hard problem in general, while it remains an open question whether the Jones polynomial can do so. To incorporate the persistence into the analysis of curve-type data, Shen et al. \cite{shen2024evolutionary} proposed the persistent Khovanov homology (or evolutionary Khovanov homology) of knots and links by defining filtrations of links obtained through the smoothing transformations of crossings in a link diagram. This work addresses the limitations that classical knot theory primarily focuses on global topological properties and may overlook the important local features of curve-type data. It enables a more refined characterization of the topological structures and geometric shapes of knots and links across various scales in a filtration. 

Conceptually, Khovanov homology for knots and links shares similarities with the classical homology for simplicial complexes, although they are developed in different contexts. Both are constructed from chain complexes, which consist of sequences of abelian groups connected by differentials or boundary maps that vanish when applied twice. A key difference lies in the grading, where Khovanov homology is bigraded, with both a homological grading and a quantum grading, whereas the simplicial homology has only a single homological grading. The graded Betti numbers in Khovanov homology reflect the graded dimensions, compared to the classical Betti numbers in simplicial homology, which count the number of holes in each dimension.

The paper presents two approaches for constructing a filtration of link diagrams: distance-based filtration and unzipping filtration. The idea is to generate an order of crossings, which are gradually resolved by a specific smoothing, finally leading to a simple link diagram consisting of disjoint circles. This process provides a filtration. The distance-based filtration considers disks centered at all crossings determined by one single radius parameter, and a crossing is isolated if its disk does not intersect with any other disks. Starting with a large initial radius, the order is then obtained by labeling isolated crossings as the radius gradually decreases. For the unzipping filtration, the order is chosen by starting with an initial crossing and labeling the crossings of the link by moving in a specific direction along the curve.

The persistent Khovanov can capture nontrivial evolutionary information even for trivial knots and unknotted links. The effectiveness of the methods has been demonstrated in the analysis of the SARS-CoV-2 frameshifting pseudo-knot in the same paper. Later, Liu et al. \cite{liu2024persistentTangles} proposed theoretically the persistent Khovanov homology of tangles, focusing on characterizing local topological characteristics in curve-type data. By enabling a persistence analysis of topological invariants across multiple scales, these frameworks provide powerful tools capable of capturing finer topological details or local structures of curve-type data, offering promising geometric topology approaches in data science.

\subsection{Khovanov Laplacians and Khovanov Dirac}

Motivated by the success of combinatorial and Hodge Laplacians and Dirac in data science, and their connections to the (co)homology, Jones and Wei \cite{jones2025khovanov} introduced the Khovanov Laplacians and Khovanov Dirac for knots and links, extending the Khovanov homology by incorporating the spectral analysis in the study of knots and links. The harmonic spectra of these operators can recover the topological invariants from Khovanov homology, while the non-harmonic spectra reveal additional information not accessible through homology alone. The framework can also be seen as an extension of the Laplacians on directed graphs associated with link diagrams \cite{silver2019knot}, which establishes a connection between graph Laplacians and key properties of knots. These graph Laplacians naturally lead to the Seifert matrix, the Alexander matrix, and the Goeritz matrix of a link, along with the invariants derived from them, including the Alexander polynomial and $\omega$-signature. 

The construction of the Khovanov Laplacians and Khovanov Dirac parallels that of the combinatorial Laplacians and the Dirac operator on simplicial complexes. As the differential in Khovanov homology increases only the homological grading by $1$ while preserving the quantum grading, one can just restrict to chain complexes with fixed quantum grading. The paper defines an inner product on a chain group by assuming that its basis elements are orthonormal, which then induces a unique adjoint of the differential. The Khovanov Laplacians are subsequently defined using the same formula as for the combinatorial Laplacians, and the Dirac operator follows the same matrix structure.
Analogous to the combinatorial case, the kernels of Khovanov Laplacians correspond to the Khovanov homology groups.

As noted in the paper \cite{jones2025khovanov}, the construction of Khovanov Laplacians and the Dirac operator may also be extended to other knot homology theories, such as the knot Floer homology \cite{ozsvath2004holomorphic}, which categorifies the Alexander polynomial. This raises possible connections between the Laplacians derived from the Floer knot homology and the Laplacians mentioned above on directed graphs \cite{silver2019knot}. In addition, in analogy with the combinatorial and differential manifold settings, one may define the persistent Khovanov Laplacians and the persistent Dirac, further generalizing the framework of persistent Khovanov homology \cite{shen2024evolutionary}.

\section{Machine learning featurization}
\label{sec:MachineLearningFeaturization}
Simplicial complexes are the most widely used tools for modeling point cloud data in topological data analysis, owing to their conceptual simplicity and computational efficiency. However, they are not universally optimal for all data types. For example, in network data, a coauthorship network is more appropriately modeled as a hypergraph rather than a simplicial complex. This is because the presence of a joint publication among three authors does not necessarily imply that each pair of authors has collaborated independently, a condition inherent in the simplicial complex construction. 
Similarly, data on manifolds such as images often arise from smooth structures and therefore possess intrinsic differential properties and boundary information. These characteristics cannot be fully captured by the purely combinatorial structure of simplicial complexes. Instead, modeling the data as differential manifolds provides a more suitable framework for incorporating smoothness and boundary information. 
Thus, the selection of a mathematical structure for modeling data in TDA should be guided by the intrinsic nature of the data and the objectives of the specific analysis task. Moreover, different modeling choices naturally lead to different classes of topological tools and representations.

In this section, we begin by categorizing various data types commonly encountered in TDA, including point-cloud data, distance matrix data, network data, data on manifolds, sequential data, knots and links data, and data with additional non-geometric information. For each category, we discuss appropriate mathematical models and their implications. We then present methods for extracting topological representations from these data structures. Finally, we describe approaches for transforming these topological representations into features suitable for machine learning applications. We mainly focus on vectorization techniques, which convert topological representations into fixed-size feature vectors. Other strategies, such as kernel-based methods, are also briefly discussed.

\subsection{Input data types and the choice of topological modeling}

\subsubsection{Point-cloud data}
Point-cloud data is one of the most commonly encountered data formats in a wide range of fields. Examples include atoms in biomolecules such as proteins, DNA, and RNA, as well as nodes in various networks such as social networks, transportation networks, and molecular interaction networks.
Modeling point cloud data using mathematical objects such as digraphs, simplicial complexes, and hypergraphs enables the application of advanced topological tools, including persistent homology and persistent combinatorial Laplacians, for comprehensive data analysis. The selection of an appropriate mathematical model for point cloud data depends on the nature and available additional information of the data. 

The Vietoris-Rips complex and Alpha complex are suitable when data points are given without any additional information or when they need to be treated equally, while the weighted simplicial complex or sheaf is better suited for data points that carry additional label information, such as type, mass, or color. The digraphs or hyperdigraphs are appropriate for representing data with asymmetry or directed relations, such as gene expression datasets, where directional information encodes gene regulation. Hypergraphs and super-hypergraphs can be effective for handling data points with complex relationships or incomplete structures.

\subsubsection{Distance matrix data}
For distance matrix data, it is important to emphasize the flexibility and diversity in the construction of distance metrics. While the Euclidean distance is the most widely used, it is not universally optimal for all applications. In many cases, modifying the distance metric or adopting alternative definitions can lead to a more accurate extraction of topological features relevant to the underlying problem.

The interactive distance \cite{cang2018integration} is designed to capture the interactions between two distinct types of objects. Taking protein-ligand interactions as an example, directly applying Euclidean distance to construct simplicial complexes from the entire protein-ligand structure can obscure the biologically meaningful interactions due to overwhelming contributions from intra-protein or intra-ligand interactions. To address this, the interactive distance assigns the original Euclidean distance to point pairs that originate from different structures (i.e., protein and ligand), while assigning an infinite distance to point pairs within the same structure. This strategy ensures that the resulting simplicial complex focuses specifically on cross-structure interactions, thereby offering a more targeted characterization of the protein-ligand interactions.

The multi-level distance \cite{cang2018representability} is proposed to highlight noncovalent atomic interactions, such as hydrogen bonds and van der Waals forces. In this approach, bonded atom pairs are assigned infinite distances since bonded atoms stay closer than non-bonded ones in most cases, allowing us to focus on non-bonded, noncovalent interactions. Although this metric may violate the triangle inequality, it remains valid for constructing Vietoris-Rips complexes, where only pairwise distances are required. This strategy leads to not only zero-dimensional topological features encoding noncovalent interactions but also higher-dimensional features that capture small structural fluctuations among different conformations of the same molecule.

The kernel function-based distance \cite{cang2018representability} was introduced to model atomic interaction strengths, as they often do not align linearly with the Euclidean distance. The frameworks uses the negative of a radial basis function to define a correlation function-based filtration matrix.
Stronger interactions correspond to smaller distances and appear earlier in the filtration process. Similarly, in the flexibility rigidity index-based distance, a Lorentz-type kernel with a tunable scale parameter is employed. By adjusting this scale parameter, one can capture structural patterns at different spatial resolutions.

The KNN-based distance strategy \cite{le2025persistent} assigns distances between pairs of points based on their mutual nearest-neighbor relationships. Specifically, for a pair of points $(a,b)$, if $a$ is the $s$-th nearest neighbor of $b$, and $b$ is the $t$-th nearest neighbor of $a$, then the distance between $a$ and $b$ is defined as the smaller of the two values, i.e., $\min(s,t)$. The standard Vietoris-Rips or Alpha filtration can then be constructed from the resulting distance matrix of the point cloud. This strategy has been used to analyze single-cell RNA sequencing data \cite{cottrell2024k}.

\subsubsection{Network data}
Network data consists of nodes and edges that represent interactions or relationships among entities. Such data arises in a wide range of domains. For example, in molecular networks, nodes correspond to atoms and edges to chemical bonds, in social networks, nodes represent individuals and edges denote social connections, and in trade networks, nodes are countries while edges indicate trading relationships. Mathematically, these networks can be modeled as graphs, enabling the application of tools from graph theory for pairwise interaction analysis. 
The choice of topological tools depends on the nature of the network and the specific application context.

The (directed) undirected network can be naturally modeled as a one-dimensional (ordered) simplicial complex, where nodes represent 0-simplices and edges represent 1-simplices. 
To analyze high-order interactions beyond pairwise connections, one can extend the graph to a higher-dimensional simplicial complex. Several constructions are available: the flag complex (or clique complex) \cite{jones2025persistent}, which includes a simplex for each complete subgraph, the neighborhood complex \cite{liu2021neighborhood}, which encodes shared neighborhood structures, and the Hom complex \cite{liu2022hom}, which can be interpreted analogously to convolutional kernels in convolutional neural networks. In addition, hypergraphs and other generalized mathematical objects can be employed to model high-order interactions directly \cite{grbic2022aspects,liu2024intcomplex}.

\subsubsection{Data on manifolds}
Common examples for data on manifolds include images, such as 3D volumes or 2D slices from MRI and CT scans, as well as isosurfaces, such as isothermal surfaces representing a uniform heat distribution in fluid dynamics and electron density isosurfaces modeling molecular orbitals in molecular visualization. 
Although they are discrete in their raw form, it is natural to treat them as data on manifolds, as they are often sampled from continuous or smooth structures. In previous studies, topological analysis of such data, particularly image data, has often been conducted using point cloud-based techniques. For example, a common approach in medical imaging is to select points from the region of interest in the image and then construct a simplicial complex from these selected points to perform topological analysis \cite{singh2023topological}. Additionally, cubical complex structures have also been employed to model images for topological analysis \cite{wang2016object}. Mathematically, both simplicial complexes and cubical complexes are combinatorial objects, so they are capable of capturing the topological structures arising from the data on manifolds. However, due to their inherently combinatorial nature, these representations neglect the differential structure and boundary conditions of the underlying manifolds, thereby limiting their ability to preserve the differential properties intrinsic to such data.

Therefore, modeling these data as differential manifolds offers a more suitable framework for capturing their intrinsic differential properties. In this setting, tools from differential topology or differential geometry, such as differential forms, vector fields, and differential operators, can be employed to facilitate quantitative analysis of the data. In addition, these data, such as images, inherently have boundaries that separate objects from the background, making the boundary conditions essential. These boundary conditions not only impose constraints that govern the behavior of the system near the boundary but also ensure the existence and uniqueness of solutions in boundary value problems. Moreover, the boundary conditions modify the differential operators, such as the Laplacian operators and the Dirac operator, making them different from their counterparts defined on closed manifolds (without boundary).
It is worth noting that the Hodge Laplacians on differential manifolds differ fundamentally from the combinatorial Laplacians defined on simplicial complexes \cite{ribando2024graph}. 
The effectiveness of differential topology-based methods has been demonstrated in various applications, such as Hodge decomposition in single-cell RNA velocity analysis \cite{su2024hodge}, PHL in drug design and discovery \cite{su2024persistent}, and MTDL \cite{liu2025manifold} in medical image analysis.

Building on this perspective, recent work has integrated topological and geometric priors into deep learning architectures. The Topological Convolutional Neural Network (TCNN) introduced by Love et al. \cite{love2023topological} replaces standard convolutional layers with manifold-parameterized filters that reflect the structure of image and video data. This approach preserves topological relationships and leverages geometric localization across layers, leading to faster learning, better generalization, and improved interpretability. TCNNs also extend naturally to 3D data, highlighting the broader potential of topology in deep learning.

\subsubsection{Sequential data}
Sequential data arise in a wide range of fields, such as protein and DNA sequences in biology, word sequences in natural language processing, temporal stock price fluctuations in time series analysis, and frame sequences in video processing. Many approaches have been developed for topological sequence analysis (TSA). One strategy for TSA is to directly compute a pairwise distance matrix \cite{chan2013topology}, which can then be used to construct Vietoris-Rips filtration by varying the distance parameter. For time series data, an approach involves dividing curves into short segments and building a Vietoris-Rips filtration for each segment \cite{zheng2024towards}. There are also TSA approaches to transform the sequence into a point cloud, thereby enabling the direct use of established techniques for point cloud data. These strategies include Takens' delay embedding, the sliding window method, and the $k$-mer topology approach.

In Takens' delay embedding \cite{takens2006detecting}, the sequence ${x_1,x_2,\cdots}$ is transformed into points using two parameters: the embedding dimension $d$ and the delay parameter $\tau$. For instance, with $d = 2$, each point in the resulting point cloud takes the form $(x_i, x_{i+\tau})$. 
The sliding window method similarly employs two parameters: the point dimension parameter $M$ and the window parameter $\tau$. For example, if $M = 3$, a point in the embedded space is constructed as $(x_i, x_{i+\tau}, x_{i+2\tau})$. There are also other variations for the sliding window method, such as SW1PerS \cite{perea2015sliding} and toroidal sliding window embeddings \cite{perea2016persistent}.

In the $k$-mer topology method \cite{hozumi2024revealing}, the parameter $k$ specifies the window size used to slide along the sequence, generating a collection of $k$-mers (subsequences of length $k$), which are treated as a point cloud. The element-specific strategy is then applied to emphasize particular types of points. The pairwise distances between these points are computed to construct a distance matrix, which serves as the foundation for subsequent topological analysis using methods such as persistent homology or persistent Laplacians. 

 Delta complex approaches have been proposed for TSA \cite{liu2025topological}. The essential idea is to consider the segments of a sequence as simplices, which serve as the building blocks of a $\Delta$-complex. Since $\Delta$-complexes allow repeated elements, they align with the nature of biological sequences. Unlike $k$-mer topology, $\Delta$-complex-based TSA methods do not systematically analyze k-mers in a sequence. These approaches eventually utilize persistent homology or persistent Laplacians and are faster than the early $k$-mer topology method for large genome sequences.  


Recent developments in TDL introduce novel methods for integrating temporal and topological structures in sequential data, particularly in the context of neural spike decoding \cite{mitchell2024topological}. These approaches enable end-to-end learning where topology is a learnable model component \cite{papamarkou2024position,karan2021time}. In many cases, the data consists of temporally ordered events, such as neural spikes, where each observation contributes to the construction of a simplicial complex that captures multi-way relationships among entities (e.g., neurons spiking together within a short time window). These complexes encode higher-order interactions beyond simple pairwise connections, providing a richer structural representation.

Rather than focusing on individual time points or fixed-size windows, some TDL methods track the temporal evolution of these simplicial complexes, producing a sequence of topological states. In event-based settings like spike decoding, this evolution is typically modeled through discrete-time or state-space processes. Other approaches, \cite{einizade2025cosmos,montagna2024topological}, model the evolution of topological features over complexes using continuous-time dynamics and partial differential equations. In these cases, each simplicial complex can be embedded as a learnable vector using neural architectures designed for topological spaces, such as simplicial message passing networks \cite{pmlr-v139-bodnar21a} or cell complex encoders.

The result is a unified, data-driven framework that captures both temporal dynamics and topological structure. In applications like neural decoding, this enables the model to learn how patterns of co-activation change over time and relate to behavior, stimulus, or latent neural states. By representing sequences as evolving topological objects, TDL provides a flexible and expressive way to analyze structured time-series data, enabling the discovery of complex patterns embedded in both time and topology.

\subsubsection{Knots, links, tangles, and braids}

Knots, links, tangles, and braids are fundamental topological objects that can be naturally modeled as curves embedded in 3-space
\cite{shen2024knot,song2025multiscale, shen2024evolutionary}. In practical applications, such data structures are commonly found in a wide range of fields. Examples include ropes, shoelaces, polymers, DNA, RNA, nucleosomes, and chromosomes.

The representation of these objects as curves in 3-space has advantages, as it allows the use of tools from knot theory for analyzing complex structure relationships. The approach has proven effective in various studies. For example, the entanglement of polymer chains can be studied through their knotting and linking properties, with the Gauss linking integral correlating with the physical properties of polymeric material \cite{panagiotou2019topological}. In addition, studies have shown that the Gauss linking integral is related to protein entanglements, particularly in understanding protein folding kinetics \cite{panagiotou2020topological}. It has also been applied to the study of SARS-CoV-2 spike protein \cite{baldwin2022local}. The effectiveness of the knot-based topological frameworks has also been demonstrated in protein flexibility analysis and protein-ligand interactions \cite{shen2024knot}, protein structure and functions \cite{song2025multiscale}, and SARS-CoV-2 frame-shifting pseudo-knot \cite{shen2024evolutionary}.

Databases such as KnotProt \cite{jamroz2015knotprot} and LinkProt \cite{dabrowski2016linkprot} have been developed to support the study of knots and links of proteins. The Python package Topoly \cite{dabrowski2021topoly} provides tools for calculating topological polynomial invariants to distinguish knots, slipknots, links, and spatial graphs.

\subsubsection{Data with additional non-geometric information}

For data with additional non-geometric information, we refer to all the aforementioned data types that are further enriched with background-related attributes beyond geometric structures. We set this as a distinct data type because such data is common in real-world applications, and the supplementary information often plays a crucial role in addressing domain-specific problems. For example, the protein structures not only contain the 3D coordinates of atoms (geometric information) but also include the residue types, partial charges, and solvent accessibility (non-geometric biological information), which are critical for protein stability and protein-ligand interactions. Similarly, in network data, beyond the topological structure of nodes and edges, one may have node-associated attributes such as the user interests in social networks, which are critical for the clustering and link prediction tasks.

The element-specific method \cite{cang2017topologynet,cang2018integration,cang2018representability} is a powerful framework for systematically modeling different atomic and molecular interactions into a unified topological system. For example, given a point cloud of molecular atoms, one can construct element-specific simplicial complexes by grouping atoms based on their element types. This enables us to focus on particular atomic interactions, such as covalent bonds, hydrogen bonds, van der Waals forces, and other noncovalent interactions. The integration of element-specific strategy with standard persistent homology has demonstrated great promise in various biological data analysis tasks, such as protein-ligand binding affinity prediction, protein-protein interaction analysis, and protein stability analysis.

Another approach for incorporating non-geometric information is to encode it as weights and integrate these weights into the underlying mathematical frameworks. For example, in the enriched barcode method \cite{cang2020persistent}, the persistent cohomology is computed from the underlying simplicial complex modeling of the protein atoms, and the additional biochemical information is embedded as weights into the cohomology representatives. This results in an enriched, or colored, persistence barcode that captures both topological and chemical properties. In the persistent sheaf Laplacian method \cite{wei2025persistent}, the atomic properties are encoded into the sheaf construction, particularly, the construction of the boundary operators. Thus the derived persistent sheaf cohomology and sheaf Laplacians will automatically carry this additional information.

\subsection{Approaches for topological representations}

\subsubsection{Representations of persistent homology/cohomology}

Persistent homology and persistent cohomology produce identical barcodes when computed with coefficients in a field \cite{de2011dualities}, which is the most common choice in practical applications. Therefore, here we focus only on the case of persistent homology for simplicity. 
Once an appropriate topological model is established for a given data set, the corresponding persistent homology can be computed to extract meaningful topological representations, such as persistent simplicial homology for simplicial complexes, persistent hypergraph homology for hypergraphs, and persistent sheaf homology for sheaves. Note that persistent homology can be visualized and encoded in various forms, each offering a distinct perspective on the topological features. Next, we summarize several common representation formats of persistent homology.

\paragraph{Persistence barcode/diagram}
The two standard representations of persistent homology, persistence barcodes \cite{carlsson2004persistence} and persistence diagrams \cite{edelsbrunner2002topological}, both offer intuitive ways to represent the birth and death of topological features during the filtration process. Persistence barcodes summarize the lifespans of topological features, such as connected components, holes, and cavities, across different dimensions. Each barcode consists of horizontal line segments, where the left endpoint of each line segment represents the birth time and the right endpoint indicates the death time of a topological feature. A persistence diagram is equivalent to a barcode but encodes the same topological information as a 2D scatter plot, where the $x$-coordinate represents the birth time and the $y$-coordinate represents the death time. As the death time of a feature always occurs after its birth time, all points in a persistence diagram are positioned above the diagonal $y=x$. Long segments in a barcode typically correspond to significant topological features, while short segments are often considered as noise. Similarly, in a persistence diagram, features far from the diagonal are usually considered significant, and those closer to the diagonal are often treated as noise. However, this assumption might not hold when considering physical, chemical and biological data, where short segments or features near the diagonal might carry important physical information such as bond lengths, benzene ring size, and other properties \cite{xia2014persistent}. There is also a Grassmannian persistence diagram construction \cite{gulen2023orthogonal,gulen2025grassmannian} for multiparameter persistence. Persistence barcodes/diagrams can effectively capture the topological features of the data, and have been used in the study of image patches \cite{carlsson2008local, adams2009nonlinear}, hand-drawn letters \cite{collins2004barcode}, and unusual events in viral evolution \cite{chan2013topology}.

\paragraph{Persistent Betti numbers/ Betti curves}

Persistent Betti numbers \cite{edelsbrunner2002topological} and Betti curves are multiscale analogs of the Betti numbers that count the number of topological features that persist in the filtration. Persistent Betti numbers are defined as the ranks of persistent homology groups and are often visualized as Betti curves, which offer a step function-based representation of persistent homology. The representation transforms persistence barcodes/diagrams into functional representations. The vector space structure of functions then enables basic operations like addition and scalar multiplication, making Betti curves suitable for statistical tasks such as computing the means and variances. These operators are not naturally defined for persistence barcodes/diagrams.

Betti curves are simple and easy to compute. In addition, as small perturbations often introduce short lifetimes that have minimal impact on the overall Betti number distribution, they are also robust to noise. However, one drawback of Betti curves is their limited expression power, as they only contain counts of topological features. In addition, the mapping from a diagram to a curve is not injective \cite{chevyrev2018persistence}, which means that different diagrams can lead to the same curve.

\paragraph{Persistence landscape}
A stable representation of persistent homology, named persistence landscapes \cite{bubenik2015statistical,bubenik2017persistence}, offers a more structured and expressive representation of persistence diagrams while maintaining compatibility with statistical and machine learning frameworks.
They were introduced to convert persistence diagrams into sequences of piecewise linear functions, which can be viewed as a horizontal version of persistence diagrams. Compared to Betti curves, which track only the count of topological features, persistence landscapes encode both the birth-and-death pairs and their persistence and maintain a bijective correspondence with persistence diagrams, making the latter preferable when detailed topological analysis is needed. Recently, it has been extended to the multiparameter persistence landscape (MPPL) \cite{vipond2020multiparameter} for the multiparameter case.

A key feature of persistence landscapes is their hierarchical structure, where different levels of functions capture topological features with varying significance. 
The function with the lowest lifespan subscript captures the most persistent topological features, while those with higher lifespan-subscript capture less persistent features. Persistence landscapes have been used to study conformational changes in protein binding\cite{kovacev2016using}, Microstructure Analysis \cite{dlotko2016topological}, time-dependent functional networks constructed from time series \cite{stolz2017persistent}, and image classification tasks and orbit classifications \cite{kim2020pllay}. 

\paragraph{Persistence image}
Persistence images \cite{adams2017persistence} were introduced as an alternative stable representation that transforms persistence diagrams into fixed-size, grid-based functional representations, making them particularly well suited for convolutional neural
networks (CNN). Persistence images are constructed by computing the integral of a persistence surface on a grid obtained by discretizing a subdomain of $\RR^2$. Here, a persistence surface is defined as a weighted sum of Gaussian functions centered at each point in the persistence diagram, with the weighting function controlling the significance of each point. Persistence images have been proven to be  stable with respect to the 1-Wasserstein distance between persistence diagrams in the same paper. 

A key advantage of persistence images is that they are easy to compute and offer flexibility in three choices: the resolution, the distribution (and its associated parameters), and the weighting function. However, these choices are noncanonical, meaning that there are no guidelines for setting these parameters to achieve the best performance in machine learning tasks. PI-Net \cite{som2020pi}, a simple one-step differentiable architecture, is the first to propose the use of deep learning for computing persistence images directly from data. The persistence image has also been extended to a multiparameter setting, leading to the Multiparameter Persistence Image (MPPI) \cite{carriere2020multiparameter}.

\paragraph{Persistent Betti number image}
The notion of multidimensional persistence was proposed in \cite{xia2015multidimensional} for transforming dynamical biomolecular data into persistent Betti number images by encoding the topological evolution of molecular systems over time. Specifically, given a sequence of data frames representing a dynamical process, persistent homology is applied to each frame, and the resulting persistent Betti numbers are encoded as one-dimensional vectors. By stacking these vectors across all time steps, a two-dimensional image is constructed that reflects the temporal topological changes of the system.
The multiscale multidimensional persistence proposed in the same paper further encodes molecular structures at various spatial scales within a single image. In particular, the flexibility-rigidity function with tunable scale parameters is used to generate a volumetric density representation of the molecule. For specific chosen scale parameters, the persistent Betti numbers are computed and recorded as one-dimensional vectors, then a 2D, 3D, or high-dimensional image is formed by assembling these vectors across all scales.

\paragraph{Colored persistence barcode }
Several methods have been proposed to enrich persistence barcodes with additional information, resulting in the colored persistence barcode. For example, in the persistent enriched barcode approach \cite{cang2020persistent}, persistent cohomology is first computed, and supplementary weight information associated with simplices is incorporated into the cohomology representatives through specific functions to assign colors to the corresponding bars. In the persistent Jones polynomial method \cite{song2025multiscale}, each simplex is assigned a Jones polynomial, and the polynomial associated with the positive simplex corresponding to a bar is evaluated at a selected value to determine its color. These enhanced barcodes effectively encode background-related information, thereby providing a more informative topological characterization.

\paragraph{Graded persistence barcode}
The most commonly used persistence barcode is derived from single-graded persistent homology. However, there are also methods that provide multi-graded persistent homology representations. For example, in the persistent Mayer homology method \cite{shen2024persistent}, the additional grading in the $N$-chain complex determines a specific type of boundary operator. By varying the grading index, a series of chain complexes can be constructed, resulting in a sequence of persistent homologies, i.e., a bigraded persistent homology. Similarly, in the persistent Khovanov homology method \cite{liu2024persistentTangles}, an additional quantum grading that is determined by the smoothing type of the crossings, is introduced into the cochain complex, resulting in a bigraded persistent homology that captures richer structural information. These graded persistent homologies naturally give rise to graded persistence barcode representations.

\subsubsection{Representations of persistent Laplacians}
The representation of Laplacians primarily involves their spectral properties, that is, their eigenvalues and eigenvectors, which encode essential information about the topological, and the combinatorial or geometric structures of the underlying datasets. Another important aspect is that the Hodge Laplacian enables the decomposition of vector fields into three or five orthogonal components through the Hodge decomposition, providing a powerful tool for analyzing flows, gradients, and cycles in data. In the following, we describe these representations in detail.

\paragraph{Eigenvalue representation}
A fundamental property of topological Laplacians is the isomorphism between their kernel and the corresponding cohomology groups. In practice, computations are typically performed with a field coefficient, under which cohomology is isomorphic to homology. This implies that the multiplicity of zero eigenvalues of the topological Laplacian is equal to the Betti number in the corresponding dimension. 
This result extends to the persistent setting: the persistent Laplacian exhibits a similar isomorphism, where the multiplicity of its zero eigenvalues corresponds to the persistent Betti numbers \cite{memoli2022persistent}. Consequently, the zero eigenvalues of the persistent Laplacian fully encode the persistent homology information. In contrast, the nonzero eigenvalues provide additional geometric and combinatorial insights of the underlying data structure.
Laplacian matrices are symmetric and semi-positive definite, which guarantees that all nonzero eigenvalues are positive. For Hodge Laplacians on manifolds, the larger eigenvalues are associated with higher-frequency components of the corresponding eigenfunctions. These high-frequency modes reflect fine-scale structural variations or local irregularities and may represent noise, especially when the manifold exhibits non-smooth features such as sharp corners. 
The integration of the eigenvalues of Laplacians with machine learning techniques has been successfully applied in a variety of domains, including drug design and discovery, computational biology and chemistry, and material science.

\paragraph{Eigenvector representation}

It is known that the eigenvectors associated with topological Laplacians capture rich topological and geometric information about the underlying objects. For example, in shape analysis, the global point signature (GPS) \cite{rustamov2007laplace}, a representation of surfaces defined using scaled eigenfunctions of the Laplace-Beltrami operator, has been shown to carry sufficient information of surfaces for various shape processing tasks, including shape classification, segmentation, and correspondence. The eigenvectors of Hodge Laplacians, including the harmonic eigenfields, the gradient eigenfields, and the curl eigenfields under certain boundary conditions, have been applied to the analysis of biological macromolecules \cite{zhao2020rham}. These eigenvectors can be used either individually as inputs in machine learning models, or as a supplement to eigenvalue features.

For persistent topological Laplacians, while the eigenvalues are known to capture the topological and geometric changes during the filtration process, the interpretation of the eigenvectors, however, remains unclear and requires further investigation.

\paragraph{Hodge decomposition representation}
Hodge decomposition orthogonally splits vector fields into gradient, curl, and harmonic components. This decomposition enables a more detailed analysis of a vector field by extracting different dynamic features, making it particularly useful in understanding the underlying mechanisms of a field. It has numerous applications in electromagnetism \cite{miller1984interpretations}, fluids and dynamic systems \cite{yang2021clebsch, Yin2023FluidCohomology}, geometric modeling~\cite{wang2021computing}, and spectral data analysis~\cite{keros2023spectral}. 

The most well-known Hodge decomposition is the 3-component Helmholtz-Hodge decomposition \cite{ladyzhenskaya1969mathematical}. For compact domains in Euclidean spaces, there is a refined 5-component Hodge decomposition of vector fields \cite{shonkwiler2009poincare}, which preserves the topology of the underlying domains. The discrete 5-component Hodge decomposition has been proposed for surface triangle meshes \cite{poelke2016boundary}, for tetrahedral meshes \cite{poelke2017hodge, zhao20193d}, and for domains defined by level set functions in Cartesian grids \cite{su2024hodge, su2024topology}.
Hodge decomposition has been successfully applied to the study of single-cell RNA velocities for extracting different dynamic information of cells by utilizing their gene expressions \cite{su2024hodgerna}. In addition, manifold topological deep learning (MTDL) \cite{liu2025manifold} has been recently proposed based on the Hodge decomposition technique for medical image analysis.

\subsubsection{Representations of persistent Dirac operator}

In the combinatorial case, the persistent Dirac operator can be seen as a generalization of the persistent combinatorial Laplacians \cite{ameneyro2024quantum, wee2023persistent}. It captures the same topological information as the persistent combinatorial Laplacians, with its eigenvalues corresponding to the square root of the eigenvalues of the persistent combinatorial Laplacians in all dimensions. However, as a first-order operator, it is more sensitive to local features compared to the second-order Laplacians, and its matrix couples different topological dimensions. 

One significant advantage of the persistent Dirac operator is that it enables efficient quantum computation, providing exponential speed-up over classical TDA algorithms. This operator not only recovers the topological outputs from persistent homology but also leads to spectral features of the point cloud data. The persistent Dirac operator can also be formulated similarly for differential manifolds, extending the Hodge-Dirac operator in the differential topology setting, thus facilitating the analysis of data on manifolds.

\subsubsection{Representations of multiscale Gauss linking integral}
The multiscale Gauss linking integral (mGLI) method \cite{shen2024knot} is developed for the analysis of one-dimensional curve data, including knots, links, and tangles. The method begins by decomposing the curve into a collection of short segments, after which the pairwise Gauss linking integrals are computed between all segment pairs. These integrals form a segmentation matrix whose entries quantify the topological linking information between curve segments. Then pairwise distances between segments are incorporated as weights into the entries, yielding a refined segmentation matrix that encodes linking structures across multiple spatial scales. By systematically varying the scale parameter, a sequence of such matrices can be constructed, each capturing localized topological linking information at a distinct resolution. Collectively, these matrices provide a comprehensive multiscale topological representation of the underlying curve structure.

\subsubsection{Representations of multiscale/persistent Jones polynomial}

The multiscale and persistent Jones polynomials \cite{song2025multiscale} extend the classical Jones polynomial \cite{barkataki2022jones} to a multiscale setting, enabling the analysis of local and global entanglements for collections of disjoint open and closed curves. Given a specific curve segmentation, the multiscale Jones polynomial framework offers a matrix of Jones polynomials, where each row corresponds to a feature vector of a curve segment determined by the Jones polynomials computed across various distance ranges. Each entry of the feature vector is the Jones polynomial of the set of curves consisting of that curve segment and all other segments that fall within a specific distance range from that segment. The persistent Jones polynomial builds on a Vietoris-Rips complex constructed from the pairwise distance matrix of curve segments, which leads to the birth, death, and persistence information of facets in the complex. By incorporating the Jones polynomials of facets, this framework leads to a weighted barcode, with weights given by the Jones polynomials of the associated facets.

\subsubsection{Software and code resources}
Table \ref{tab:software} lists several widely used software packages for computing persistent homology, persistent Laplacians, and Hodge decompositions. In addition to these, another popular tool is Mapper \cite{singh2007topological}, which provides a simplified representation of high-dimensional datasets in the form of simplicial complexes. Taking point cloud data as an example, the data is first divided into clusters using a clustering algorithm. A simplicial complex is then constructed from these clusters: each cluster corresponds to a vertex, and a set of 
$n$ clusters forms an $(n\!-\!1)$-simplex if their intersection is nonempty. The one-skeleton of the resulting simplicial complex is commonly referred to as the Mapper graph, which provides a simple description of the data and captures important information about its topological structures. Several software tools have implemented the Mapper algorithm, including KeplerMapper \cite{van2019kepler}, giotto-tda \cite{tauzin2021giotto}, GUDHI \cite{maria2014gudhi}, and Mapper Interactive \cite{zhou2021mapper}.

\begin{table}[h]
\centering
\caption{Software packages for computing persistent homology and persistent Laplacians.  }
\label{tab:software}
\resizebox{\textwidth}{!}{
\begin{tabular}{|l|c|c|c|c|c|c| }
    \hline
    Software & Simplicial & Cubical & Persistent & Persistent & Persistent & Hodge  \\
    & Complex & Complex & Homology & Cohomology & Laplacian & Decom.\\
    \hline
    \href{JavaPlex}{https://appliedtopology.github.io/javaplex/} \cite{adams2014javaplex} & Yes & & Yes & & & \\
    \hline
    \href{Perseus}{https://people.maths.ox.ac.uk/nanda/perseus/index.html} \cite{mischaikow2013morse} & Yes & Yes & Yes & & & \\
    \hline
    \href{Dinoysus}{https://mrzv.org/software/dionysus2/} \cite{morozov2007dionysus} & Yes & & Yes & Yes & & \\
    \hline
    \href{GUDHI}{https://gudhi.inria.fr/} \cite{maria2014gudhi} & Yes & Yes & Yes & Yes & & \\
    \hline
    \href{Ripser}{https://ripser.scikit-tda.org/en/latest/} \cite{bauer2021ripser} & Yes & & Yes & Yes & & \\
    \hline
    \href{PHAT}{https://github.com/xoltar/phat} \cite{bauer2017phat} & Yes & Yes & Yes & & &\\
    \hline
    \href{DIPHA}{https://github.com/DIPHA/dipha/} \cite{bauer2014dipha} & Yes & Yes & Yes & & &\\
    \hline
    \href{R-TDA}{https://cran.r-project.org/web/packages/TDA/index.html} \cite{fasy2014introduction} & Yes & Yes & Yes & Yes & & \\
    \hline
    \href{Giotto-tda}{https://giotto-ai.github.io/gtda-docs/0.5.1/library.html} \cite{tauzin2021giotto} & Yes & Yes & Yes & & & \\
    \hline
    \href{OpenPH}{https://github.com/rodrgo/OpenPH} \cite{mendoza2017parallel} & Yes & & Yes & & &  \\
    \hline
    \href{CTL}{https://github.com/appliedtopology/ctl} & Yes & Yes & Yes & & & \\
    \hline
    \href{Eirene}{https://github.com/henselman-petrusek/Eirene.jl} \cite{henselmanghristl6} & Yes & Yes & Yes & & & \\
    \hline
    \href{Cubicle}{https://bitbucket.org/hubwag/cubicle/src/master/} \cite{wagner2023slice} &  & Yes & Yes & & &\\
    \hline
    \href{CubicalRips}{https://github.com/shizuo-kaji/CubicalRipser_3dim} \cite{kaji2020cubical} & & Yes & Yes & & & \\
    \hline
    \href{TTK}{https://topology-tool-kit.github.io/} \cite{tierny2017topology} & Yes & Yes & Yes & & & \\
    \hline
    \href{libstick}{https://www.sthu.org/code/libstick/} \cite{edelsbrunner2010computational} & Yes & & Yes & & &  \\
    \hline
    \href{HERMES}{https://weilab.math.msu.edu/HERMES/} \cite{wang2021hermes} & Yes & & & &Yes &  \\
    \hline
    \href{PersistLap}{https://github.com/ndag/Persistent-Laplacian} \cite{memoli2022persistent} & Yes& & & & Yes & \\
    \hline
    \href{HHD}{https://github.com/zhesu1/HHD} \cite{su2024hodgerna} & &Yes & & & & Yes\\
    \hline
    \href{5ComponentHD}{https://github.com/zhesu1/5ComponentHD} \cite{su2024hodge} & &Yes & & & & Yes \\
    \hline
    \href{3DHodgeDecom}{https://github.com/rdzhao/3DHodgeDecomposition} \cite{zhao20193d} &Yes & & & & & Yes \\
    \hline
\end{tabular}
}
\end{table}

\subsection{Featurization methods}
Persistent homology and persistent Jones polynomial can be encoded as persistence barcodes, which consist of a collection of (birth, death) pairs representing the lifespan of topological features across the filtration. In contrast, the multiscale Gauss linking integral and multiscale Jones polynomial are typically represented by segmentation matrices, where the entries qualify the local linking information between segments of the underlying data. For persistent Laplacians and persistent Dirac operator, the spectral information (eigenvalue and eigenvector) is usually used for feature extraction. These representations give rise to three types of topological descriptors: persistence barcode, persistence spectral, and segmentation matrix. In addition to these, the Hodge Laplacians enable three- and five-component orthogonal decompositions of vector fields, which naturally provide a multi-channel image representation of the original vector field in the Eulerian setting. In the following, we describe methods for transforming these topological descriptors into feature vectors suitable for machine learning algorithms.

\subsubsection{Functional formulations}
Persistence barcodes/diagrams can not be directly applied to machine learning tasks, as they are not naturally structured as fixed-size feature vectors. Persistence barcodes consist of varying numbers of intervals with different lengths, while persistent diagrams contain varying numbers of unordered points. One straightforward approach for vectorizing persistent barcodes/diagrams is to use their functional formulations, which can be directly discretized as machine learning inputs by sampling the resulting functions at selected points. Commonly used functional formulations include Betti curves, persistence landscape \cite{bubenik2015statistical}, persistence surfaces and persistence images \cite{adams2017persistence}, etc.

\subsubsection{Counting in bins}
Another method for vectorizing persistence barcode is the binning-based approach \cite{cang2018representability}. In this method, the uniform-length filtration range is divided into a collection of bins, and the numbers of births, deaths, and persistences within each bin are counted. This process yields three corresponding feature vectors: birth count, death count, and persistence count. This approach effectively captures structural properties such as bond lengths, including those of noncovalent interactions in biomolecules, and has been widely adopted in various topological data analysis studies \cite{cang2018integration,cang2017topologynet,nguyen2019mathematical}.

\subsubsection{Statistical property}
For persistence barcode, a basic and widely used approach for vectorization is to extract their statistical properties \cite{ali2023survey,cang2015topological, cang2018representability,asaad2022persistent}. For example, one can compute descriptive statistics of the bar lengths, birth times, death times, and midpoints, including the maximum, minimum, sum, median, mean, standard deviation, variance, and entropy. This yields a feature vector of dimension $4k$, where $k$ denotes the number of statistical measures for each of the four quantities. A more fine-grained representation can be obtained by partitioning the filtration range into $s$ fixed bins and computing statistical properties within each bin. For instance, counting the number of bars whose birth-death intervals cover each bin produces an $s$-dimensional vector, while recording the maximum and minimum birth values in each bin results in a $2s$-dimensional vector.

For persistence spectra, a common practice is to divide the filtration range into uniform intervals and compute the corresponding Laplacian matrices at each interval, the eigenvalues and eigenvectors of these matrices are then extracted to form feature representations \cite{meng2021persistent, liu2021persistent}. Basic statistical measures, such as maximum, minimum, sum, median, mean, standard deviation, variance, and entropy, can be applied to the eigenvalues of each Laplacian matrix, yielding feature vectors that capture topological, geometric, and combinatorial characteristics across scales. However, when the Laplacian matrices are large, computing the full spectrum becomes computationally expensive. In such cases, it is common to retain only the smallest $k$ nonzero eigenvalues, as well as the Betti numbers, which capture the most global and stable spectral components. This truncated spectral representation is not only computationally efficient but also robust to noise.

The segmentation matrix encodes the pairwise linking information between segments of a given pair of curves. A straightforward vectorization strategy is to compute global statistical properties, such as the maximum, minimum, sum, median, mean, standard deviation, variance, and entropy, over all entries in the matrix to construct feature vectors. Alternatively, localized linking patterns can be captured by summing the entries of each row or column, thereby quantifying the overall linking contribution of individual segments. Statistical measures can then be applied to these localized values across all segments to generate feature vectors that preserve local structural information.

\subsubsection{Algebraic formulations}
For persistence barcode, algebraic vectorization is primarily achieved through polynomial maps constructed from the barcode. For example, in the Adcock-Carlsson coordinate method \cite{adcock2013ring}, given a barcode with $n$ intervals, a ring of algebraic functions, polynomials in $2n$ variables that satisfy certain properties, is constructed. A $k$-dimensional feature vector is then obtained by selecting $k$ algebraic functions from this ring and evaluating them on the $n$ (birth, death) pairs. In the tropical coordinate method \cite{kalivsnik2019tropical}, for a barcode with $n$ intervals, a function involving $2n$ variables is constructed using only the operations of maximum, minimum, addition, and subtraction. This function is invariant under permutation of the variables. The feature vector is obtained by evaluating the function on all $n$ (birth, death) pairs or (birth, persistence) pairs. There are also methods that utilize complex polynomials defined on barcode intervals to generate feature vectors \cite{ferri1999representing,di2015comparing}.

For persistence spectra, vectorization primarily involves applying functions to the eigenvalues of the associated Laplacian matrices. For example, in the persistent attribute method \cite{meng2021persistent}, the filtration range is divided into uniform bins, and spectral indices of the eigenvalues are computed within each bin. These spectral indices \cite{puzyn2010recent} include the Laplacian graph energy, generalized graph energy, generalized average graph energy, spectral radius, spectral diameter, Ivanciuc matrix spectrum operators, spectral moments, Lov{\'a}sz–Pelik{\'a}n index, and quasi-Wiener index. In addition to these spectral descriptors, other spectral properties, such as algebraic connectivity, modularity, Cheeger constant, vertex/edge expansion, and quantities derived from flow, random walks, and heat kernels, can also be incorporated into the vectorization process. Beyond these classical spectral properties, more advanced functions such as the Riemann zeta function \cite{liu2022dowker} and symmetry functions \cite{behler2007generalized} have also been applied to eigenvalues for feature extraction.

\subsubsection{Multidimensional formulations}

Several formulations incorporate multiple parameters or indices into the feature generation process. For example, Mayer homology includes an additional index in the filtration \cite{shen2024persistent}, resulting in 3D persistence barcodes. In addition, multidimensional persistent homology provides an effective framework for transforming dynamic molecular data into image representations, and the multiscale multidimensional persistence further encodes molecular structures at various spatial scales within a single image \cite{xia2015multidimensional}.

For persistence spectra, one approach to generate a multidimensional feature vector is to use the filtration values as the $x$-axis and plot either the first $k$ eigenvalues or the statistical properties of the entire eigenvalue spectrum on the $y$-axis. An alternative strategy is to extend the idea of the persistence image to spectral features. Specifically, for each filtration value, one computes the first $k$ eigenvalues or $k$ statistical measures of the full spectrum, resulting in $k$ (filtration, eigen-property) pairs. These pairs can be treated analogously to (birth, death) pairs in persistence barcodes, allowing for the construction of a persistence image. If all such pairs are collected across the entire filtration process, a 2D image representation is obtained. Alternatively, if a separate image is constructed at each filtration value and these are concatenated along the filtration axis, a 3D image can be derived. 

Additionally, when applying Hodge Laplacians to Cartesian grids, image representations can be obtained directly via Hodge decomposition. For example, applying the three-component decomposition on a grid yields curl-free, divergence-free, and harmonic components, each of which has the same spatial dimensions as the original one-form. A multi-channel image can then be constructed by concatenating these three components along the channel axis.

\subsubsection{Kernel-based models}

Kernel-based machine learning models rely on kernel functions to implicitly map data into high-dimensional feature spaces, enabling the modeling of complex nonlinear relationships without explicitly computing the coordinates in those spaces. The core component of these methods is the kernel function, which acts as a similarity measure between all pairs of data points. 

To integrate topological representations such as persistence diagrams into kernel-based learning frameworks, it is essential to define a notion of comparability between such diagrams. A natural way is to define metric structures on the space of persistence diagrams. These metric structures induce distances that can quantify the differences between persistence barcodes/diagrams. Several metrics, such as Hausdorff distance, bottleneck distance \cite{cohen2005stability}, and Wasserstein distance \cite{mileyko2011probability,cang2018representability}, have been studied and shown to be robust to small perturbations. However, the computation of these distances, especially for large datasets, can be computationally expensive \cite{di2015comparing}. It is also challenging to define an appropriate averaging process, as the local mean may not be unique \cite{mileyko2011probability, munch2015probabilistic}. In particular, these metrics are primarily used for distance-based learning methods such as k-nearest neighbors (KNN) and some clustering algorithms. Therefore, they are unsuitable for machine learning techniques that require a Hilbert space structure, or fixed-size feature inputs, such as support vector machines (SVM), Principal Component Analysis (PCA), decision tree classification, and deep neural networks.

A kernel-based model addresses the issue by defining a kernel that maps persistence diagrams into a Hilbert space, allowing them to be directly integrated with kernel-based machine learning techniques, such as SVM and PCA. The kernel encodes the similarity between persistence diagrams by computing an inner product in the Hilbert space. Based on this construction, various persistent homology-based kernels, such as persistence scale space kernel \cite{reininghaus2015stable}, persistence weighted Gaussian kernel \cite{kusano2016persistence}, sliced Wasserstein kernel \cite{carriere2017sliced}, and persistence Fisher kernel \cite{le2018persistence}, have been proposed and demonstrated success in applications as reported in their papers, including shape classification, shape segmentation, texture recognition, etc.

\subsubsection{Statistical Frameworks for Persistent Homology}\label{stat-framework}

The study of persistence diagrams through statistical methods has advanced significantly in the past decade, addressing the challenges posed by their unordered nature and variable cardinality, especially in the presence of small data perturbations. The metric structure of diagram space, particularly under Wasserstein or bottleneck distances, naturally leads to statistical notions of central tendency, such as Fr\'echet means and medians \cite{turner2014frechet, mileyko2011probability,cohen2005stability}.

To facilitate the application of classical statistical principles within topological data analysis, persistence diagrams have been increasingly modeled as random variables. The development of functional representations that embed persistence diagrams into more tractable mathematical spaces has played a key role in driving this progress. For example, persistence landscapes \cite{bubenik2015statistical} embed persistence diagrams into Banach spaces, allowing for the direct computation of means, variances, and confidence bands, thus integrating topological summaries more directly into established statistical methodologies. Similarly, persistence images \cite{adams2017persistence} offer a complementary approach, transforming diagrams into stable, finite-dimensional vector representations that readily support standard statistical analyses. These embeddings have facilitated the establishment of Central Limit Theorems (CLT) for persistence diagrams under Wasserstein distances and analogues of the Law of Large Numbers (LLN) for various functional representations \cite{divol2021understanding,bubenik2015statistical}, strengthening TDA as a more robust and statistically grounded field, and enhancing its applicability across a wider range of scientific and engineering domains. Furthermore, bootstrap methods have been developed to assess population-level differences in topological summaries \cite{chazal2014stochastic,chazal2018robust}, enabling formal comparisons between groups by constructing confidence sets and analyzing the asymptotic behavior of topological features \cite{fasy2014confidence}.

Sufficient statistics encapsulate all the information in a sample that is relevant to the underlying probability distribution. A tropical geometry-based framework constructs sufficient statistics for persistent homology \cite{monod2020tropical}. The developed model maps barcodes into tropical coordinates that are injective and stable. This embedding into a Euclidean space via tropical geometry retains the complete topological information required for parametric modeling while enabling a seamless integration with classical statistical techniques.

Kernel Density Estimation (KDE)-based frameworks work directly on the space of persistence diagrams. Some approaches define expected diagrams as mean objects in a space of Radon measures, drawing on point process theory \cite{divol2019density}, while others utilize Wasserstein and bottleneck distances to capture the intrinsic geometry of topological features \cite{maroulas2019nonparametricPDF}. These methods estimate the distribution of topological features without converting diagrams into vector representations, offering a more intrinsic and geometrically faithful characterization of their variability.

Maroulas et al. \cite{maroulas2019nonparametricPDF}  introduced nonparametric density estimation techniques tailored for persistence diagrams using the framework of random finite sets. Using Probability Hypothesis Density (PHD), they estimate the distribution of topological features without restrictive assumptions, enabling a robust separation of meaningful topological features from noise, which is essential to analyze complex or noisy data. They provide theoretical convergence guarantees, including bounds on the mean absolute bottleneck deviation, ensuring that density estimates converge to the true distribution.

In a complementary direction, a Bayesian framework for persistent homology has been developed in which persistence diagrams are modeled as realizations from  Poisson point processes with prior intensities, enabling posterior inference over topological feature distributions and supporting a probabilistic interpretation of topological summaries \cite{maroulas2020bayesian}. 
In a related development focused on generative modeling, Papamarkou et al. \cite{papamarkou2022random} introduced a method to generate synthetic persistence diagrams based on pairwise interacting point processes and RJ-MCMC sampling, enabling realistic simulation of topological feature distributions while capturing dependencies among features. This helps perform statistical analysis in topological data analysis, especially when data is limited, such as in materials science.
A related model extends this probabilistic framework to the classification of biological networks by modeling persistence diagrams as independent and identically distributed cluster point processes, using Gaussian mixtures to capture spatial distribution and binomial distributions to model variability in the number of topological features \cite{maroulas2022bayesian}. Additionally, Bayesian topological signal processing further advances this line of work to handle time-series data, supporting uncertainty-aware classification of dynamic topological features in noisy or evolving signals \cite{oballe2022bayesian}.

\section{Concluding remarks}
\label{sec:conclusion}
Aided by artificial intelligence (AI), topological data analysis (TDA) has become a rapidly growing field in applied mathematics and data science. A key technique of TDA is persistent homology, a multiscale algebraic topology tool that has demonstrated impressive success across diverse disciplines in science, engineering, medicine, industry, and defense. However, persistent homology with simplicial complex faces several limitations, including: inability to handle non-geometric information, reliance on qualitative descriptions of high-order cycles, inability to track non-topological changes, incapability in handling directed networks, and insensitivity to structured data such as hypergraphs. Additionally, while algebraic topology techniques typically analyze point-cloud data, they are not directly applicable to data on differentiable manifolds or one-dimensional curves embedded in 3-space. These challenges in TDA have motivated significant advancements in the past decade.

This paper provides a comprehensive review of TDA techniques developed to address persistent homology’s limitations, emphasizing their mathematical foundations and strengths. Some techniques extend persistent homology from simplicial complexes to other topological domains, such as cell complexes, path complexes, directed flag complexes, cellular sheaves, and hyper(di)graphs. Others introduce alternative topological formulations, including Laplacians \cite{wang2020persistent}, Dirac operators \cite{ameneyro2024quantum}, sheaf theory, Mayer topology \cite{shen2024persistent}, and interaction topology  \cite{liu2023interaction}. For data on differentiable manifolds, methods rooted in differential topology, such as persistent de Rham cohomology, persistent Hodge Laplacians, and Hodge decomposition, have been developed \cite{chen2021evolutionary,su2024persistent}. For data on one-dimensional curves in 3-space, tools from geometric topology, such as the Gauss linking integral \cite{shen2024knot}, Jones polynomial \cite{song2025multiscale}, and Khovanov homology \cite{shen2024evolutionary,liu2024persistentTangles}, have been proposed.

The rapid growth of TDA is driven by machine learning, such as topological deep learning (TDL) \cite{cang2017topologynet} and mathematical AI. However, topological representations from TDA cannot be directly implemented in machine learning without appropriate vectorization. This paper also reviews the selection of topological algorithms for different input data formats, their associated topological representations, and methods for topological vectorization in machine learning contexts.

TDA and TDL are dynamic and rapidly evolving fields. The future development of TDA and TDL will be driven by both theoretical curiosity and practical needs, including the demand for interpretable deep neural networks. We envision future advancements in the following directions.

Currently, the local geometry and/or topology surrounding individual data points can be described by persistent cohomology \cite{cang2020persistent}, evolutionary homology \cite{cang2020evolutionary}, interaction topology \cite{liu2023interaction}, and persistent sheaf Laplacian \cite{wei2025persistentsheaf}. Further efforts are needed to refine local topology methods for practical applications.

Most current TDA approaches were developed for point-cloud data. Sequential data, however, are widely used in applications such as DNA sequencing, natural language processing, and electroencephalogram (EEG) analysis. These areas require novel TDA methodologies.

In differential geometry, the Atiyah–Singer index theory establishes a topological index. A persistent index theory for families of elliptic differential operators could be developed for data on compact manifolds.

The advancement of low-dimensional topology methods is critical for analyzing 1D curves embedded in 3-space (e.g., entangled polymers and neural networks). Persistent Floer homology may offer an effective framework for modeling protein-protein binding interactions. Additionally, robust computational algorithms for persistent Khovanov homology and persistent Khovanov Laplacians remain an urgent priority.

Finally, large language models (LLMs) herald a new era in AI, creating opportunities across disciplines. Chatbots like DeepSeek and ChatGPT can transform abstract mathematical theories—including computational topology-into practical tools \cite{liu2024chatgpt}. Rapid growth is anticipated in AI-assisted topology and topology-enhanced AI in the coming years.

TDA is a vast interdisciplinary field, spanning algebraic, differential, and geometric topology; computational algorithms; AI-driven methodologies; and real-world applications. While this paper aims to comprehensively review TDA and TDL beyond persistent homology, we inevitably omit many significant contributions due to the field’s breadth. 
 We note that tremendous efforts have been made to develop topological neural networks since the introduction of TDL in 2017 \cite{cang2017topologynet}. 
 To focus on TDA, we have restricted ourselves from elaborating on this development. The interested reader is referred to a recent review on various aspects of TDL\cite{papamarkou2024position}.  
We hope this work serves as a guide for researchers exploring topological techniques and helps them select optimal approaches for analyzing specific problems across diverse fields.


\vspace*{1cm}

\section*{Acknowledgments}
This work was supported in part by NIH grants R01AI164266, and R35GM148196, NSF grant DMS-2052983, MSU Research Foundation, and Bristol-Myers Squibb 65109. VM's has been partially funded by the National Science Foundation under grants DMS-2012609, DMR-2309083, DRL-2314155, DMS-1821241, MCB-1715794, and DGE-2152168; the U.S. Army Research Laboratory (ARL) under contracts W911NF-22-2-0143, W911NF-21-2-0186, and W911NF-19-2-0328; the U.S. Army Research Office (ARO) under contracts W911NF-21-1-0094 and W911NF-17-1-0313; and Thor Industries/ARL under contract W911NF-17-2-0141.

\vspace*{1cm}
\bibliography{refs}


\begin{thebibliography}{334}
\ifx \bisbn   \undefined \def \bisbn  #1{ISBN #1}\fi
\ifx \binits  \undefined \def \binits#1{#1}\fi
\ifx \bauthor  \undefined \def \bauthor#1{#1}\fi
\ifx \batitle  \undefined \def \batitle#1{#1}\fi
\ifx \bjtitle  \undefined \def \bjtitle#1{#1}\fi
\ifx \bvolume  \undefined \def \bvolume#1{\textbf{#1}}\fi
\ifx \byear  \undefined \def \byear#1{#1}\fi
\ifx \bissue  \undefined \def \bissue#1{#1}\fi
\ifx \bfpage  \undefined \def \bfpage#1{#1}\fi
\ifx \blpage  \undefined \def \blpage #1{#1}\fi
\ifx \burl  \undefined \def \burl#1{\textsf{#1}}\fi
\ifx \doiurl  \undefined \def \doiurl#1{\url{https://doi.org/#1}}\fi
\ifx \betal  \undefined \def \betal{\textit{et al.}}\fi
\ifx \binstitute  \undefined \def \binstitute#1{#1}\fi
\ifx \binstitutionaled  \undefined \def \binstitutionaled#1{#1}\fi
\ifx \bctitle  \undefined \def \bctitle#1{#1}\fi
\ifx \beditor  \undefined \def \beditor#1{#1}\fi
\ifx \bpublisher  \undefined \def \bpublisher#1{#1}\fi
\ifx \bbtitle  \undefined \def \bbtitle#1{#1}\fi
\ifx \bedition  \undefined \def \bedition#1{#1}\fi
\ifx \bseriesno  \undefined \def \bseriesno#1{#1}\fi
\ifx \blocation  \undefined \def \blocation#1{#1}\fi
\ifx \bsertitle  \undefined \def \bsertitle#1{#1}\fi
\ifx \bsnm \undefined \def \bsnm#1{#1}\fi
\ifx \bsuffix \undefined \def \bsuffix#1{#1}\fi
\ifx \bparticle \undefined \def \bparticle#1{#1}\fi
\ifx \barticle \undefined \def \barticle#1{#1}\fi
\bibcommenthead
\ifx \bconfdate \undefined \def \bconfdate #1{#1}\fi
\ifx \botherref \undefined \def \botherref #1{#1}\fi
\ifx \url \undefined \def \url#1{\textsf{#1}}\fi
\ifx \bchapter \undefined \def \bchapter#1{#1}\fi
\ifx \bbook \undefined \def \bbook#1{#1}\fi
\ifx \bcomment \undefined \def \bcomment#1{#1}\fi
\ifx \oauthor \undefined \def \oauthor#1{#1}\fi
\ifx \citeauthoryear \undefined \def \citeauthoryear#1{#1}\fi
\ifx \endbibitem  \undefined \def \endbibitem {}\fi
\ifx \bconflocation  \undefined \def \bconflocation#1{#1}\fi
\ifx \arxivurl  \undefined \def \arxivurl#1{\textsf{#1}}\fi
\csname PreBibitemsHook\endcsname

\bibitem{adams2004knot}
\begin{bbook}
\bauthor{\bsnm{Adams}, \binits{C.C.}}:
\bbtitle{The Knot Book: An Elementary Introduction to the Mathematical Theory
  of Knots}.
\bpublisher{American Mathematical Society},
\blocation{\hspace{0pt}}
(\byear{2004})
\end{bbook}
\endbibitem

\bibitem{adams2009nonlinear}
\begin{barticle}
\bauthor{\bsnm{Adams}, \binits{H.}},
\bauthor{\bsnm{Carlsson}, \binits{G.}}:
\batitle{On the nonlinear statistics of range image patches}.
\bjtitle{SIAM Journal on Imaging Sciences}
\bvolume{2}(\bissue{1}),
\bfpage{110}--\blpage{117}
(\byear{2009})
\end{barticle}
\endbibitem

\bibitem{adams2014javaplex}
\begin{bchapter}
\bauthor{\bsnm{Adams}, \binits{H.}},
\bauthor{\bsnm{Tausz}, \binits{A.}},
\bauthor{\bsnm{Vejdemo-Johansson}, \binits{M.}}:
\bctitle{Javaplex: a research software package for persistent (co)homology}.
In: \bbtitle{Mathematical Software--ICMS 2014: 4th International Congress,
  Seoul, South Korea, August 5-9, 2014. Proceedings 4},
pp. \bfpage{129}--\blpage{136}
(\byear{2014}).
\bcomment{Springer}
\end{bchapter}
\endbibitem

\bibitem{adams2017persistence}
\begin{barticle}
\bauthor{\bsnm{Adams}, \binits{H.}},
\bauthor{\bsnm{Emerson}, \binits{T.}},
\bauthor{\bsnm{Kirby}, \binits{M.}},
\bauthor{\bsnm{Neville}, \binits{R.}},
\bauthor{\bsnm{Peterson}, \binits{C.}},
\bauthor{\bsnm{Shipman}, \binits{P.}},
\bauthor{\bsnm{Chepushtanova}, \binits{S.}},
\bauthor{\bsnm{Hanson}, \binits{E.}},
\bauthor{\bsnm{Motta}, \binits{F.}},
\bauthor{\bsnm{Ziegelmeier}, \binits{L.}}:
\batitle{Persistence images: a stable vector representation of persistent
  homology}.
\bjtitle{Journal of Machine Learning Research}
\bvolume{18}(\bissue{8}),
\bfpage{1}--\blpage{35}
(\byear{2017})
\end{barticle}
\endbibitem

\bibitem{adams2024persistent}
\begin{botherref}
\oauthor{\bsnm{Adams}, \binits{H.}},
\oauthor{\bsnm{Lagoda}, \binits{E.}},
\oauthor{\bsnm{Moy}, \binits{M.}},
\oauthor{\bsnm{Sadovek}, \binits{N.}},
\oauthor{\bsnm{De~Saha}, \binits{A.}}:
Persistent equivariant cohomology.
arXiv preprint arXiv:2408.17331
(2024)
\end{botherref}
\endbibitem

\bibitem{adcock2013ring}
\begin{barticle}
\bauthor{\bsnm{Adcock}, \binits{A.}},
\bauthor{\bsnm{Carlsson}, \binits{E.}},
\bauthor{\bsnm{Carlsson}, \binits{G.}}:
\batitle{The ring of algebraic functions on persistence bar codes}.
\bjtitle{Homology, Homotopy and Applications}
\bvolume{18}(\bissue{1}),
\bfpage{381}--\blpage{402}
(\byear{2016})
\end{barticle}
\endbibitem

\bibitem{alexander1928topological}
\begin{barticle}
\bauthor{\bsnm{Alexander}, \binits{J.W.}}:
\batitle{Topological invariants of knots and links}.
\bjtitle{Transactions of the American Mathematical Society}
\bvolume{30}(\bissue{2}),
\bfpage{275}--\blpage{306}
(\byear{1928})
\end{barticle}
\endbibitem

\bibitem{ali2023survey}
\begin{barticle}
\bauthor{\bsnm{Ali}, \binits{D.}},
\bauthor{\bsnm{Asaad}, \binits{A.}},
\bauthor{\bsnm{Jimenez}, \binits{M.-J.}},
\bauthor{\bsnm{Nanda}, \binits{V.}},
\bauthor{\bsnm{Paluzo-Hidalgo}, \binits{E.}},
\bauthor{\bsnm{Soriano-Trigueros}, \binits{M.}}:
\batitle{A survey of vectorization methods in topological data analysis}.
\bjtitle{IEEE Transactions on Pattern Analysis and Machine Intelligence}
\bvolume{45}(\bissue{12}),
\bfpage{14069}--\blpage{14080}
(\byear{2023})
\end{barticle}
\endbibitem

\bibitem{ameneyro2024quantum}
\begin{barticle}
\bauthor{\bsnm{Ameneyro}, \binits{B.}},
\bauthor{\bsnm{Maroulas}, \binits{V.}},
\bauthor{\bsnm{Siopsis}, \binits{G.}}:
\batitle{Quantum persistent homology}.
\bjtitle{Journal of Applied and Computational Topology}
\bvolume{8}(\bissue{7}),
\bfpage{1961}--\blpage{1980}
(\byear{2024})
\end{barticle}
\endbibitem

\bibitem{ameneyro2023quantum}
\begin{bchapter}
\bauthor{\bsnm{Ameneyro}, \binits{B.}},
\bauthor{\bsnm{Siopsis}, \binits{G.}},
\bauthor{\bsnm{Maroulas}, \binits{V.}}:
\bctitle{Quantum persistent homology for time series}.
In: \bbtitle{APS March Meeting Abstracts},
vol. \bseriesno{2023},
pp. \bfpage{73}--\blpage{003}
(\byear{2023})
\end{bchapter}
\endbibitem

\bibitem{arnold2018finite}
\begin{bbook}
\bauthor{\bsnm{Arnold}, \binits{D.N.}}:
\bbtitle{Finite Element Exterior Calculus}.
\bpublisher{SIAM},
\blocation{\hspace{0pt}}
(\byear{2018})
\end{bbook}
\endbibitem

\bibitem{arnold2006finite}
\begin{barticle}
\bauthor{\bsnm{Arnold}, \binits{D.N.}},
\bauthor{\bsnm{Falk}, \binits{R.S.}},
\bauthor{\bsnm{Winther}, \binits{R.}}:
\batitle{Finite element exterior calculus, homological techniques, and
  applications}.
\bjtitle{Acta numerica}
\bvolume{15},
\bfpage{1}--\blpage{155}
(\byear{2006})
\end{barticle}
\endbibitem

\bibitem{arsuaga2005dna}
\begin{barticle}
\bauthor{\bsnm{Arsuaga}, \binits{J.}},
\bauthor{\bsnm{Vazquez}, \binits{M.}},
\bauthor{\bsnm{McGuirk}, \binits{P.}},
\bauthor{\bsnm{Trigueros}, \binits{S.}},
\bauthor{\bsnm{Sumners}, \binits{D.W.}},
\bauthor{\bsnm{Roca}, \binits{J.}}:
\batitle{{DNA} knots reveal a chiral organization of {DNA} in phage capsids}.
\bjtitle{Proceedings of the National Academy of Sciences}
\bvolume{102}(\bissue{26}),
\bfpage{9165}--\blpage{9169}
(\byear{2005})
\end{barticle}
\endbibitem

\bibitem{asaad2022persistent}
\begin{barticle}
\bauthor{\bsnm{Asaad}, \binits{A.}},
\bauthor{\bsnm{Ali}, \binits{D.}},
\bauthor{\bsnm{Majeed}, \binits{T.}},
\bauthor{\bsnm{Rashid}, \binits{R.}}:
\batitle{Persistent homology for breast tumor classification using mammogram
  scans}.
\bjtitle{Mathematics}
\bvolume{10}(\bissue{21}),
\bfpage{4039}
(\byear{2022})
\end{barticle}
\endbibitem

\bibitem{atienza2019persistent}
\begin{barticle}
\bauthor{\bsnm{Atienza}, \binits{N.}},
\bauthor{\bsnm{Gonzalez-Diaz}, \binits{R.}},
\bauthor{\bsnm{Rucco}, \binits{M.}}:
\batitle{Persistent entropy for separating topological features from noise in
  {V}ietoris-{R}ips complexes}.
\bjtitle{Journal of Intelligent Information Systems}
\bvolume{52},
\bfpage{637}--\blpage{655}
(\byear{2019})
\end{barticle}
\endbibitem

\bibitem{baccini2022weighted}
\begin{barticle}
\bauthor{\bsnm{Baccini}, \binits{F.}},
\bauthor{\bsnm{Geraci}, \binits{F.}},
\bauthor{\bsnm{Bianconi}, \binits{G.}}:
\batitle{Weighted simplicial complexes and their representation power of
  higher-order network data and topology}.
\bjtitle{Physical Review E}
\bvolume{106}(\bissue{3}),
\bfpage{034319}
(\byear{2022})
\end{barticle}
\endbibitem

\bibitem{bae2017beyond}
\begin{bchapter}
\bauthor{\bsnm{Bae}, \binits{W.}},
\bauthor{\bsnm{Yoo}, \binits{J.}},
\bauthor{\bsnm{Chul~Ye}, \binits{J.}}:
\bctitle{Beyond deep residual learning for image restoration: Persistent
  homology-guided manifold simplification}.
In: \bbtitle{Proceedings of the IEEE Conference on Computer Vision and Pattern
  Recognition Workshops},
pp. \bfpage{145}--\blpage{153}
(\byear{2017})
\end{bchapter}
\endbibitem

\bibitem{baldwin2021local}
\begin{barticle}
\bauthor{\bsnm{Baldwin}, \binits{Q.}},
\bauthor{\bsnm{Panagiotou}, \binits{E.}}:
\batitle{The local topological free energy of proteins}.
\bjtitle{Journal of Theoretical Biology}
\bvolume{529},
\bfpage{110854}
(\byear{2021})
\end{barticle}
\endbibitem

\bibitem{baldwin2022local}
\begin{barticle}
\bauthor{\bsnm{Baldwin}, \binits{Q.}},
\bauthor{\bsnm{Sumpter}, \binits{B.}},
\bauthor{\bsnm{Panagiotou}, \binits{E.}}:
\batitle{The local topological free energy of the {SARS-CoV-2} spike protein}.
\bjtitle{Polymers}
\bvolume{14}(\bissue{15}),
\bfpage{3014}
(\byear{2022})
\end{barticle}
\endbibitem

\bibitem{barkataki2022jones}
\begin{barticle}
\bauthor{\bsnm{Barkataki}, \binits{K.}},
\bauthor{\bsnm{Panagiotou}, \binits{E.}}:
\batitle{The {J}ones polynomial of collections of open curves in 3-space}.
\bjtitle{Proceedings of the Royal Society A}
\bvolume{478}(\bissue{2267}),
\bfpage{20220302}
(\byear{2022})
\end{barticle}
\endbibitem

\bibitem{barnes2021comparative}
\begin{barticle}
\bauthor{\bsnm{Barnes}, \binits{D.}},
\bauthor{\bsnm{Polanco}, \binits{L.}},
\bauthor{\bsnm{Perea}, \binits{J.A.}}:
\batitle{A comparative study of machine learning methods for persistence
  diagrams}.
\bjtitle{Frontiers in Artificial Intelligence}
\bvolume{4},
\bfpage{681174}
(\byear{2021})
\end{barticle}
\endbibitem

\bibitem{bauer2021ripser}
\begin{barticle}
\bauthor{\bsnm{Bauer}, \binits{U.}}:
\batitle{Ripser: efficient computation of {Vietoris--Rips} persistence
  barcodes}.
\bjtitle{Journal of Applied and Computational Topology}
\bvolume{5}(\bissue{3}),
\bfpage{391}--\blpage{423}
(\byear{2021})
\end{barticle}
\endbibitem

\bibitem{bauer2014dipha}
\begin{botherref}
\oauthor{\bsnm{Bauer}, \binits{U.}},
\oauthor{\bsnm{Kerber}, \binits{M.}},
\oauthor{\bsnm{Reininghaus}, \binits{J.}}:
Dipha (a distributed persistent homology algorithm).
Software available at https://github. com/DIPHA/dipha
(2014)
\end{botherref}
\endbibitem

\bibitem{bauer2017phat}
\begin{barticle}
\bauthor{\bsnm{Bauer}, \binits{U.}},
\bauthor{\bsnm{Kerber}, \binits{M.}},
\bauthor{\bsnm{Reininghaus}, \binits{J.}},
\bauthor{\bsnm{Wagner}, \binits{H.}}:
\batitle{Phat--persistent homology algorithms toolbox}.
\bjtitle{Journal of symbolic computation}
\bvolume{78},
\bfpage{76}--\blpage{90}
(\byear{2017})
\end{barticle}
\endbibitem

\bibitem{behler2007generalized}
\begin{barticle}
\bauthor{\bsnm{Behler}, \binits{J.}},
\bauthor{\bsnm{Parrinello}, \binits{M.}}:
\batitle{Generalized neural-network representation of high-dimensional
  potential-energy surfaces}.
\bjtitle{Physical review letters}
\bvolume{98}(\bissue{14}),
\bfpage{146401}
(\byear{2007})
\end{barticle}
\endbibitem

\bibitem{bhatia2012helmholtz}
\begin{barticle}
\bauthor{\bsnm{Bhatia}, \binits{H.}},
\bauthor{\bsnm{Norgard}, \binits{G.}},
\bauthor{\bsnm{Pascucci}, \binits{V.}},
\bauthor{\bsnm{Bremer}, \binits{P.-T.}}:
\batitle{The {Helmholtz-Hodge} decomposition - a survey}.
\bjtitle{IEEE Transactions on visualization and computer graphics}
\bvolume{19}(\bissue{8}),
\bfpage{1386}--\blpage{1404}
(\byear{2012})
\end{barticle}
\endbibitem

\bibitem{bi2022cayley}
\begin{botherref}
\oauthor{\bsnm{Bi}, \binits{W.}},
\oauthor{\bsnm{Li}, \binits{J.}},
\oauthor{\bsnm{Liu}, \binits{J.}},
\oauthor{\bsnm{Wu}, \binits{J.}}:
On the {Cayley}-persistence algebra.
arXiv preprint arXiv:2205.10796
(2022)
\end{botherref}
\endbibitem

\bibitem{bianconi2021topological}
\begin{barticle}
\bauthor{\bsnm{Bianconi}, \binits{G.}}:
\batitle{The topological {Dirac} equation of networks and simplicial
  complexes}.
\bjtitle{Journal of Physics: Complexity}
\bvolume{2}(\bissue{3}),
\bfpage{035022}
(\byear{2021})
\end{barticle}
\endbibitem

\bibitem{biasotti2008reeb}
\begin{barticle}
\bauthor{\bsnm{Biasotti}, \binits{S.}},
\bauthor{\bsnm{Giorgi}, \binits{D.}},
\bauthor{\bsnm{Spagnuolo}, \binits{M.}},
\bauthor{\bsnm{Falcidieno}, \binits{B.}}:
\batitle{Reeb graphs for shape analysis and applications}.
\bjtitle{Theoretical computer science}
\bvolume{392}(\bissue{1-3}),
\bfpage{5}--\blpage{22}
(\byear{2008})
\end{barticle}
\endbibitem

\bibitem{pmlr-v139-bodnar21a}
\begin{bchapter}
\bauthor{\bsnm{Bodnar}, \binits{C.}},
\bauthor{\bsnm{Frasca}, \binits{F.}},
\bauthor{\bsnm{Wang}, \binits{Y.}},
\bauthor{\bsnm{Otter}, \binits{N.}},
\bauthor{\bsnm{Mont{\'u}far}, \binits{G.F.}},
\bauthor{\bsnm{Li{\`o}}, \binits{P.}},
\bauthor{\bsnm{Bronstein}, \binits{M.}}:
\bctitle{Weisfeiler and lehman go topological: Message passing simplicial
  networks}.
In: \bbtitle{Proceedings of the 38th International Conference on Machine
  Learning},
pp. \bfpage{1026}--\blpage{1037}
(\byear{2021})
\end{bchapter}
\endbibitem

\bibitem{botnan2023introduction}
\begin{bchapter}
\bauthor{\bsnm{Botnan}, \binits{M.}},
\bauthor{\bsnm{Lesnick}, \binits{M.}}:
\bctitle{An introduction to multiparameter persistence}.
In: \bbtitle{Representations of Algebras and Related Structures},
pp. \bfpage{77}--\blpage{150}
(\byear{2023})
\end{bchapter}
\endbibitem

\bibitem{botnan2020rectangle}
\begin{bchapter}
\bauthor{\bsnm{Botnan}, \binits{M.B.}},
\bauthor{\bsnm{Lebovici}, \binits{V.}},
\bauthor{\bsnm{Oudot}, \binits{S.}}:
\bctitle{On rectangle-decomposable 2-parameter persistence modules}.
In: \bbtitle{Leibniz International Proceedings in Informatics},
vol. \bseriesno{164},
pp. \bfpage{22}--\blpage{1}
(\byear{2020}).
\bcomment{Leibniz-Zentrum f{\"u}r Informatik}
\end{bchapter}
\endbibitem

\bibitem{bressan2019embedded}
\begin{barticle}
\bauthor{\bsnm{Bressan}, \binits{S.}},
\bauthor{\bsnm{Li}, \binits{J.}},
\bauthor{\bsnm{Ren}, \binits{S.}},
\bauthor{\bsnm{Wu}, \binits{J.}}:
\batitle{The embedded homology of hypergraphs and applications}.
\bjtitle{Asian Journal of Mathematics}
\bvolume{23}(\bissue{3}),
\bfpage{479}--\blpage{500}
(\byear{2019})
\end{barticle}
\endbibitem

\bibitem{bubenik2015statistical}
\begin{barticle}
\bauthor{\bsnm{Bubenik}, \binits{P.}}:
\batitle{Statistical topological data analysis using persistence landscapes.}
\bjtitle{Journal of Machine Learning Research}
\bvolume{16}(\bissue{1}),
\bfpage{77}--\blpage{102}
(\byear{2015})
\end{barticle}
\endbibitem

\bibitem{bubenik2017persistence}
\begin{barticle}
\bauthor{\bsnm{Bubenik}, \binits{P.}},
\bauthor{\bsnm{D{\l}otko}, \binits{P.}}:
\batitle{A persistence landscapes toolbox for topological statistics}.
\bjtitle{Journal of Symbolic Computation}
\bvolume{78},
\bfpage{91}--\blpage{114}
(\byear{2017})
\end{barticle}
\endbibitem

\bibitem{calmon2023dirac}
\begin{barticle}
\bauthor{\bsnm{Calmon}, \binits{L.}},
\bauthor{\bsnm{Schaub}, \binits{M.T.}},
\bauthor{\bsnm{Bianconi}, \binits{G.}}:
\batitle{Dirac signal processing of higher-order topological signals}.
\bjtitle{New Journal of Physics}
\bvolume{25}(\bissue{9}),
\bfpage{093013}
(\byear{2023})
\end{barticle}
\endbibitem

\bibitem{cang2017topologynet}
\begin{barticle}
\bauthor{\bsnm{Cang}, \binits{Z.}},
\bauthor{\bsnm{Wei}, \binits{G.-W.}}:
\batitle{Topologynet: Topology based deep convolutional and multi-task neural
  networks for biomolecular property predictions}.
\bjtitle{PLoS computational biology}
\bvolume{13}(\bissue{7}),
\bfpage{1005690}
(\byear{2017})
\end{barticle}
\endbibitem

\bibitem{cang2018integration}
\begin{barticle}
\bauthor{\bsnm{Cang}, \binits{Z.}},
\bauthor{\bsnm{Wei}, \binits{G.-W.}}:
\batitle{Integration of element specific persistent homology and machine
  learning for protein-ligand binding affinity prediction}.
\bjtitle{International journal for numerical methods in biomedical engineering}
\bvolume{34}(\bissue{2}),
\bfpage{2914}
(\byear{2018})
\end{barticle}
\endbibitem

\bibitem{cang2020persistent}
\begin{barticle}
\bauthor{\bsnm{Cang}, \binits{Z.}},
\bauthor{\bsnm{Wei}, \binits{G.-W.}}:
\batitle{Persistent cohomology for data with multicomponent heterogeneous
  information}.
\bjtitle{SIAM journal on mathematics of data science}
\bvolume{2}(\bissue{2}),
\bfpage{396}--\blpage{418}
(\byear{2020})
\end{barticle}
\endbibitem

\bibitem{cang2018representability}
\begin{barticle}
\bauthor{\bsnm{Cang}, \binits{Z.}},
\bauthor{\bsnm{Mu}, \binits{L.}},
\bauthor{\bsnm{Wei}, \binits{G.-W.}}:
\batitle{Representability of algebraic topology for biomolecules in machine
  learning based scoring and virtual screening}.
\bjtitle{PLoS computational biology}
\bvolume{14}(\bissue{1}),
\bfpage{1005929}
(\byear{2018})
\end{barticle}
\endbibitem

\bibitem{cang2020evolutionary}
\begin{barticle}
\bauthor{\bsnm{Cang}, \binits{Z.}},
\bauthor{\bsnm{Munch}, \binits{E.}},
\bauthor{\bsnm{Wei}, \binits{G.-W.}}:
\batitle{Evolutionary homology on coupled dynamical systems with applications
  to protein flexibility analysis}.
\bjtitle{Journal of applied and computational topology}
\bvolume{4}(\bissue{4}),
\bfpage{481}--\blpage{507}
(\byear{2020})
\end{barticle}
\endbibitem

\bibitem{cang2015topological}
\begin{botherref}
\oauthor{\bsnm{Cang}, \binits{Z.}},
\oauthor{\bsnm{Mu}, \binits{L.}},
\oauthor{\bsnm{Wu}, \binits{K.}},
\oauthor{\bsnm{Opron}, \binits{K.}},
\oauthor{\bsnm{Xia}, \binits{K.}},
\oauthor{\bsnm{Wei}, \binits{G.-W.}}:
A topological approach for protein classification.
Computational and Mathematical Biophysics
\textbf{3}(1)
(2015)
\end{botherref}
\endbibitem

\bibitem{cantarella2002vector}
\begin{barticle}
\bauthor{\bsnm{Cantarella}, \binits{J.}},
\bauthor{\bsnm{DeTurck}, \binits{D.}},
\bauthor{\bsnm{Gluck}, \binits{H.}}:
\batitle{Vector calculus and the topology of domains in 3-space}.
\bjtitle{The American Mathematical Monthly}
\bvolume{109}(\bissue{5}),
\bfpage{409}--\blpage{442}
(\byear{2002})
\end{barticle}
\endbibitem

\bibitem{carlsson2009topology}
\begin{barticle}
\bauthor{\bsnm{Carlsson}, \binits{G.}}:
\batitle{Topology and data}.
\bjtitle{Bulletin of the American Mathematical Society}
\bvolume{46}(\bissue{2}),
\bfpage{255}--\blpage{308}
(\byear{2009})
\end{barticle}
\endbibitem

\bibitem{carlsson2010zigzag}
\begin{barticle}
\bauthor{\bsnm{Carlsson}, \binits{G.}},
\bauthor{\bsnm{De~Silva}, \binits{V.}}:
\batitle{Zigzag persistence}.
\bjtitle{Foundations of computational mathematics}
\bvolume{10},
\bfpage{367}--\blpage{405}
(\byear{2010})
\end{barticle}
\endbibitem

\bibitem{carlsson2007theory}
\begin{bchapter}
\bauthor{\bsnm{Carlsson}, \binits{G.}},
\bauthor{\bsnm{Zomorodian}, \binits{A.}}:
\bctitle{The theory of multidimensional persistence}.
In: \bbtitle{Proceedings of the Twenty-third Annual Symposium on Computational
  Geometry},
pp. \bfpage{184}--\blpage{193}
(\byear{2007})
\end{bchapter}
\endbibitem

\bibitem{carlsson2009zigzag}
\begin{bchapter}
\bauthor{\bsnm{Carlsson}, \binits{G.}},
\bauthor{\bsnm{De~Silva}, \binits{V.}},
\bauthor{\bsnm{Morozov}, \binits{D.}}:
\bctitle{Zigzag persistent homology and real-valued functions}.
In: \bbtitle{Proceedings of the Twenty-fifth Annual Symposium on Computational
  Geometry},
pp. \bfpage{247}--\blpage{256}
(\byear{2009})
\end{bchapter}
\endbibitem

\bibitem{carlsson2004persistence}
\begin{bchapter}
\bauthor{\bsnm{Carlsson}, \binits{G.}},
\bauthor{\bsnm{Zomorodian}, \binits{A.}},
\bauthor{\bsnm{Collins}, \binits{A.}},
\bauthor{\bsnm{Guibas}, \binits{L.}}:
\bctitle{Persistence barcodes for shapes}.
In: \bbtitle{Proceedings of the 2004 Eurographics/ACM SIGGRAPH Symposium on
  Geometry Processing},
pp. \bfpage{124}--\blpage{135}
(\byear{2004})
\end{bchapter}
\endbibitem

\bibitem{carlsson2008local}
\begin{barticle}
\bauthor{\bsnm{Carlsson}, \binits{G.}},
\bauthor{\bsnm{Ishkhanov}, \binits{T.}},
\bauthor{\bsnm{De~Silva}, \binits{V.}},
\bauthor{\bsnm{Zomorodian}, \binits{A.}}:
\batitle{On the local behavior of spaces of natural images}.
\bjtitle{International journal of computer vision}
\bvolume{76},
\bfpage{1}--\blpage{12}
(\byear{2008})
\end{barticle}
\endbibitem

\bibitem{carlsson2019parametrized}
\begin{barticle}
\bauthor{\bsnm{Carlsson}, \binits{G.}},
\bauthor{\bsnm{De~Silva}, \binits{V.}},
\bauthor{\bsnm{Kali{\v{s}}nik}, \binits{S.}},
\bauthor{\bsnm{Morozov}, \binits{D.}}:
\batitle{Parametrized homology via zigzag persistence}.
\bjtitle{Algebraic \& Geometric Topology}
\bvolume{19}(\bissue{2}),
\bfpage{657}--\blpage{700}
(\byear{2019})
\end{barticle}
\endbibitem

\bibitem{gudhi:CoverComplex}
\begin{bchapter}
\bauthor{\bsnm{Carri{\`{e}}re}, \binits{M.}}:
\bctitle{Cover complex}.
In: \bbtitle{GUDHI User and Reference Manual},
(\byear{2025})
\end{bchapter}
\endbibitem

\bibitem{carriere2020multiparameter}
\begin{barticle}
\bauthor{\bsnm{Carriere}, \binits{M.}},
\bauthor{\bsnm{Blumberg}, \binits{A.}}:
\batitle{Multiparameter persistence image for topological machine learning}.
\bjtitle{Advances in Neural Information Processing Systems}
\bvolume{33},
\bfpage{22432}--\blpage{22444}
(\byear{2020})
\end{barticle}
\endbibitem

\bibitem{carriere2017sliced}
\begin{bchapter}
\bauthor{\bsnm{Carriere}, \binits{M.}},
\bauthor{\bsnm{Cuturi}, \binits{M.}},
\bauthor{\bsnm{Oudot}, \binits{S.}}:
\bctitle{Sliced {W}asserstein kernel for persistence diagrams}.
In: \bbtitle{International Conference on Machine Learning},
pp. \bfpage{664}--\blpage{673}
(\byear{2017}).
\bcomment{PMLR}
\end{bchapter}
\endbibitem

\bibitem{chan2013topology}
\begin{barticle}
\bauthor{\bsnm{Chan}, \binits{J.M.}},
\bauthor{\bsnm{Carlsson}, \binits{G.}},
\bauthor{\bsnm{Rabadan}, \binits{R.}}:
\batitle{Topology of viral evolution}.
\bjtitle{Proceedings of the National Academy of Sciences}
\bvolume{110}(\bissue{46}),
\bfpage{18566}--\blpage{18571}
(\byear{2013})
\end{barticle}
\endbibitem

\bibitem{chazal2014stochastic}
\begin{bchapter}
\bauthor{\bsnm{Chazal}, \binits{F.}},
\bauthor{\bsnm{Fasy}, \binits{B.T.}},
\bauthor{\bsnm{Lecci}, \binits{F.}},
\bauthor{\bsnm{Rinaldo}, \binits{A.}},
\bauthor{\bsnm{Wasserman}, \binits{L.}}:
\bctitle{Stochastic convergence of persistence landscapes and silhouettes}.
In: \bbtitle{Proceedings of the Thirtieth Annual Symposium on Computational
  Geometry},
pp. \bfpage{474}--\blpage{483}
(\byear{2014})
\end{bchapter}
\endbibitem

\bibitem{chazal2016structure}
\begin{bbook}
\bauthor{\bsnm{Chazal}, \binits{F.}},
\bauthor{\bsnm{De~Silva}, \binits{V.}},
\bauthor{\bsnm{Glisse}, \binits{M.}},
\bauthor{\bsnm{Oudot}, \binits{S.}}:
\bbtitle{The Structure and Stability of Persistence Modules}
vol. \bseriesno{10}.
\bpublisher{Springer},
\blocation{\hspace{0pt}}
(\byear{2016})
\end{bbook}
\endbibitem

\bibitem{chazal2018robust}
\begin{barticle}
\bauthor{\bsnm{Chazal}, \binits{F.}},
\bauthor{\bsnm{Fasy}, \binits{B.T.}},
\bauthor{\bsnm{Lecci}, \binits{F.}},
\bauthor{\bsnm{Michel}, \binits{B.}},
\bauthor{\bsnm{Rinaldo}, \binits{A.}},
\bauthor{\bsnm{Wasserman}, \binits{L.}}:
\batitle{Robust topological inference: Distance to a measure and kernel
  distance}.
\bjtitle{Journal of Machine Learning Research}
\bvolume{18}(\bissue{159}),
\bfpage{1}--\blpage{40}
(\byear{2018})
\end{barticle}
\endbibitem

\bibitem{chen2024multiscale}
\begin{barticle}
\bauthor{\bsnm{Chen}, \binits{D.}},
\bauthor{\bsnm{Liu}, \binits{J.}},
\bauthor{\bsnm{Wei}, \binits{G.-W.}}:
\batitle{Multiscale topology-enabled structure-to-sequence transformer for
  protein--ligand interaction predictions}.
\bjtitle{Nature Machine Intelligence}
\bvolume{6}(\bissue{7}),
\bfpage{799}--\blpage{810}
(\byear{2024})
\end{barticle}
\endbibitem

\bibitem{chen2023path}
\begin{barticle}
\bauthor{\bsnm{Chen}, \binits{D.}},
\bauthor{\bsnm{Liu}, \binits{J.}},
\bauthor{\bsnm{Wu}, \binits{J.}},
\bauthor{\bsnm{Wei}, \binits{G.-W.}},
\bauthor{\bsnm{Pan}, \binits{F.}},
\bauthor{\bsnm{Yau}, \binits{S.-T.}}:
\batitle{Path topology in molecular and materials sciences}.
\bjtitle{The journal of physical chemistry letters}
\bvolume{14}(\bissue{4}),
\bfpage{954}--\blpage{964}
(\byear{2023})
\end{barticle}
\endbibitem

\bibitem{chen2023persistent}
\begin{barticle}
\bauthor{\bsnm{Chen}, \binits{D.}},
\bauthor{\bsnm{Liu}, \binits{J.}},
\bauthor{\bsnm{Wu}, \binits{J.}},
\bauthor{\bsnm{Wei}, \binits{G.-W.}}:
\batitle{Persistent hyperdigraph homology and persistent hyperdigraph
  {L}aplacians}.
\bjtitle{Foundations of Data Science}
\bvolume{5}(\bissue{4}),
\bfpage{558}--\blpage{588}
(\byear{2023})
\end{barticle}
\endbibitem

\bibitem{chen2008efficient}
\begin{barticle}
\bauthor{\bsnm{Chen}, \binits{G.}},
\bauthor{\bsnm{Mischaikow}, \binits{K.}},
\bauthor{\bsnm{Laramee}, \binits{R.S.}},
\bauthor{\bsnm{Zhang}, \binits{E.}}:
\batitle{Efficient {M}orse decompositions of vector fields}.
\bjtitle{IEEE Transactions on Visualization and Computer Graphics}
\bvolume{14}(\bissue{4}),
\bfpage{848}--\blpage{862}
(\byear{2008})
\end{barticle}
\endbibitem

\bibitem{chen2022omicronBA2}
\begin{barticle}
\bauthor{\bsnm{Chen}, \binits{J.}},
\bauthor{\bsnm{Wei}, \binits{G.-W.}}:
\batitle{Omicron {BA}. 2 ({B}. 1.1. 529.2): high potential for becoming the
  next dominant variant}.
\bjtitle{The journal of physical chemistry letters}
\bvolume{13}(\bissue{17}),
\bfpage{3840}--\blpage{3849}
(\byear{2022})
\end{barticle}
\endbibitem

\bibitem{chen2020mutations}
\begin{barticle}
\bauthor{\bsnm{Chen}, \binits{J.}},
\bauthor{\bsnm{Wang}, \binits{R.}},
\bauthor{\bsnm{Wang}, \binits{M.}},
\bauthor{\bsnm{Wei}, \binits{G.-W.}}:
\batitle{Mutations strengthened {SARS-CoV-2} infectivity}.
\bjtitle{Journal of molecular biology}
\bvolume{432}(\bissue{19}),
\bfpage{5212}--\blpage{5226}
(\byear{2020})
\end{barticle}
\endbibitem

\bibitem{chen2021evolutionary}
\begin{barticle}
\bauthor{\bsnm{Chen}, \binits{J.}},
\bauthor{\bsnm{Zhao}, \binits{R.}},
\bauthor{\bsnm{Tong}, \binits{Y.}},
\bauthor{\bsnm{Wei}, \binits{G.-W.}}:
\batitle{Evolutionary de {R}ham-{H}odge method}.
\bjtitle{Discrete and continuous dynamical systems. Series B}
\bvolume{26}(\bissue{7}),
\bfpage{3785}
(\byear{2021})
\end{barticle}
\endbibitem

\bibitem{chen2022persistent}
\begin{barticle}
\bauthor{\bsnm{Chen}, \binits{J.}},
\bauthor{\bsnm{Qiu}, \binits{Y.}},
\bauthor{\bsnm{Wang}, \binits{R.}},
\bauthor{\bsnm{Wei}, \binits{G.-W.}}:
\batitle{Persistent {L}aplacian projected {Omicron BA. 4 and BA. 5} to become
  new dominating variants}.
\bjtitle{Computers in Biology and Medicine}
\bvolume{151},
\bfpage{106262}
(\byear{2022})
\end{barticle}
\endbibitem

\bibitem{chen2023topological}
\begin{barticle}
\bauthor{\bsnm{Chen}, \binits{J.}},
\bauthor{\bsnm{Woldring}, \binits{D.R.}},
\bauthor{\bsnm{Huang}, \binits{F.}},
\bauthor{\bsnm{Huang}, \binits{X.}},
\bauthor{\bsnm{Wei}, \binits{G.-W.}}:
\batitle{Topological deep learning based deep mutational scanning}.
\bjtitle{Computers in biology and medicine}
\bvolume{164},
\bfpage{107258}
(\byear{2023})
\end{barticle}
\endbibitem

\bibitem{chen2021z}
\begin{bchapter}
\bauthor{\bsnm{Chen}, \binits{Y.}},
\bauthor{\bsnm{Segovia}, \binits{I.}},
\bauthor{\bsnm{Gel}, \binits{Y.R.}}:
\bctitle{{Z-GCNETs}: time zigzags at graph convolutional networks for time
  series forecasting}.
In: \bbtitle{International Conference on Machine Learning},
pp. \bfpage{1684}--\blpage{1694}
(\byear{2021}).
\bcomment{PMLR}
\end{bchapter}
\endbibitem

\bibitem{chevyrev2018persistence}
\begin{barticle}
\bauthor{\bsnm{Chevyrev}, \binits{I.}},
\bauthor{\bsnm{Nanda}, \binits{V.}},
\bauthor{\bsnm{Oberhauser}, \binits{H.}}:
\batitle{Persistence paths and signature features in topological data
  analysis}.
\bjtitle{IEEE transactions on pattern analysis and machine intelligence}
\bvolume{42}(\bissue{1}),
\bfpage{192}--\blpage{202}
(\byear{2018})
\end{barticle}
\endbibitem

\bibitem{chowdhury2018persistent}
\begin{bchapter}
\bauthor{\bsnm{Chowdhury}, \binits{S.}},
\bauthor{\bsnm{M{\'e}moli}, \binits{F.}}:
\bctitle{Persistent path homology of directed networks}.
In: \bbtitle{Proceedings of the Twenty-Ninth Annual ACM-SIAM Symposium on
  Discrete Algorithms},
pp. \bfpage{1152}--\blpage{1169}
(\byear{2018}).
\bcomment{SIAM}
\end{bchapter}
\endbibitem

\bibitem{chung2022persistence}
\begin{barticle}
\bauthor{\bsnm{Chung}, \binits{Y.-M.}},
\bauthor{\bsnm{Lawson}, \binits{A.}}:
\batitle{Persistence curves: A canonical framework for summarizing persistence
  diagrams}.
\bjtitle{Advances in Computational Mathematics}
\bvolume{48}(\bissue{1}),
\bfpage{6}
(\byear{2022})
\end{barticle}
\endbibitem

\bibitem{ciarlet2002finite}
\begin{bbook}
\bauthor{\bsnm{Ciarlet}, \binits{P.G.}}:
\bbtitle{The Finite Element Method for Elliptic Problems}.
\bpublisher{SIAM},
\blocation{\hspace{0pt}}
(\byear{2002})
\end{bbook}
\endbibitem

\bibitem{clough2020topological}
\begin{barticle}
\bauthor{\bsnm{Clough}, \binits{J.R.}},
\bauthor{\bsnm{Byrne}, \binits{N.}},
\bauthor{\bsnm{Oksuz}, \binits{I.}},
\bauthor{\bsnm{Zimmer}, \binits{V.A.}},
\bauthor{\bsnm{Schnabel}, \binits{J.A.}},
\bauthor{\bsnm{King}, \binits{A.P.}}:
\batitle{A topological loss function for deep-learning based image segmentation
  using persistent homology}.
\bjtitle{IEEE transactions on pattern analysis and machine intelligence}
\bvolume{44}(\bissue{12}),
\bfpage{8766}--\blpage{8778}
(\byear{2020})
\end{barticle}
\endbibitem

\bibitem{cohen2005stability}
\begin{bchapter}
\bauthor{\bsnm{Cohen-Steiner}, \binits{D.}},
\bauthor{\bsnm{Edelsbrunner}, \binits{H.}},
\bauthor{\bsnm{Harer}, \binits{J.}}:
\bctitle{Stability of persistence diagrams}.
In: \bbtitle{Proceedings of the Twenty-first Annual Symposium on Computational
  Geometry},
pp. \bfpage{263}--\blpage{271}
(\byear{2005})
\end{bchapter}
\endbibitem

\bibitem{cohen2009extending}
\begin{barticle}
\bauthor{\bsnm{Cohen-Steiner}, \binits{D.}},
\bauthor{\bsnm{Edelsbrunner}, \binits{H.}},
\bauthor{\bsnm{Harer}, \binits{J.}}:
\batitle{Extending persistence using {P}oincar{\'e} and {L}efschetz duality}.
\bjtitle{Foundations of Computational Mathematics}
\bvolume{9}(\bissue{1}),
\bfpage{79}--\blpage{103}
(\byear{2009})
\end{barticle}
\endbibitem

\bibitem{collins2004barcode}
\begin{barticle}
\bauthor{\bsnm{Collins}, \binits{A.}},
\bauthor{\bsnm{Zomorodian}, \binits{A.}},
\bauthor{\bsnm{Carlsson}, \binits{G.}},
\bauthor{\bsnm{Guibas}, \binits{L.J.}}:
\batitle{A barcode shape descriptor for curve point cloud data}.
\bjtitle{Computers \& Graphics}
\bvolume{28}(\bissue{6}),
\bfpage{881}--\blpage{894}
(\byear{2004})
\end{barticle}
\endbibitem

\bibitem{conley1978isolated}
\begin{bbook}
\bauthor{\bsnm{Conley}, \binits{C.C.}}:
\bbtitle{Isolated Invariant Sets and the {M}orse Index}
vol. \bseriesno{38}.
\bpublisher{American Mathematical Soc.},
\blocation{\hspace{0pt}}
(\byear{1978})
\end{bbook}
\endbibitem

\bibitem{conti2022topological}
\begin{barticle}
\bauthor{\bsnm{Conti}, \binits{F.}},
\bauthor{\bsnm{Moroni}, \binits{D.}},
\bauthor{\bsnm{Pascali}, \binits{M.A.}}:
\batitle{A topological machine learning pipeline for classification}.
\bjtitle{Mathematics}
\bvolume{10}(\bissue{17}),
\bfpage{3086}
(\byear{2022})
\end{barticle}
\endbibitem

\bibitem{cooperband2025unified}
\begin{botherref}
\oauthor{\bsnm{Cooperband}, \binits{Z.}},
\oauthor{\bsnm{Ghrist}, \binits{R.}}:
Unified origami kinematics via cosheaf homology.
arXiv preprint arXiv:2501.02581
(2025)
\end{botherref}
\endbibitem

\bibitem{cooperband2023cosheaf}
\begin{botherref}
\oauthor{\bsnm{Cooperband}, \binits{Z.}},
\oauthor{\bsnm{Ghrist}, \binits{R.}},
\oauthor{\bsnm{Hansen}, \binits{J.}}:
A cosheaf theory of reciprocal figures: Planar and higher genus graphic
  statics.
arXiv preprint arXiv:2311.12946
(2023)
\end{botherref}
\endbibitem

\bibitem{cottrell2024k}
\begin{barticle}
\bauthor{\bsnm{Cottrell}, \binits{S.}},
\bauthor{\bsnm{Hozumi}, \binits{Y.}},
\bauthor{\bsnm{Wei}, \binits{G.-W.}}:
\batitle{K-nearest-neighbors induced topological {PCA} for single cell
  {RNA}-sequence data analysis}.
\bjtitle{Computers in biology and medicine}
\bvolume{175},
\bfpage{108497}
(\byear{2024})
\end{barticle}
\endbibitem

\bibitem{cottrell2023plpca}
\begin{barticle}
\bauthor{\bsnm{Cottrell}, \binits{S.}},
\bauthor{\bsnm{Wang}, \binits{R.}},
\bauthor{\bsnm{Wei}, \binits{G.-W.}}:
\batitle{{PLPCA}: persistent {L}aplacian-enhanced {PCA} for microarray data
  analysis}.
\bjtitle{Journal of chemical information and modeling}
\bvolume{64}(\bissue{7}),
\bfpage{2405}--\blpage{2420}
(\byear{2023})
\end{barticle}
\endbibitem

\bibitem{crowell2012introduction}
\begin{bbook}
\bauthor{\bsnm{Crowell}, \binits{R.H.}},
\bauthor{\bsnm{Fox}, \binits{R.H.}}:
\bbtitle{Introduction to Knot Theory}
vol. \bseriesno{57}.
\bpublisher{Springer},
\blocation{\hspace{0pt}}
(\byear{2012})
\end{bbook}
\endbibitem

\bibitem{curry2014sheaves}
\begin{bbook}
\bauthor{\bsnm{Curry}, \binits{J.M.}}:
\bbtitle{Sheaves, Cosheaves and Applications}.
\bpublisher{University of Pennsylvania},
\blocation{\hspace{0pt}}
(\byear{2014})
\end{bbook}
\endbibitem

\bibitem{dabaghian2012topological}
\begin{botherref}
\oauthor{\bsnm{Dabaghian}, \binits{Y.}},
\oauthor{\bsnm{M{\'e}moli}, \binits{F.}},
\oauthor{\bsnm{Frank}, \binits{L.}},
\oauthor{\bsnm{Carlsson}, \binits{G.}}:
A topological paradigm for hippocampal spatial map formation using persistent
  homology
(2012)
\end{botherref}
\endbibitem

\bibitem{dabrowski2016linkprot}
\begin{botherref}
\oauthor{\bsnm{Dabrowski-Tumanski}, \binits{P.}},
\oauthor{\bsnm{Jarmolinska}, \binits{A.I.}},
\oauthor{\bsnm{Niemyska}, \binits{W.}},
\oauthor{\bsnm{Rawdon}, \binits{E.J.}},
\oauthor{\bsnm{Millett}, \binits{K.C.}},
\oauthor{\bsnm{Sulkowska}, \binits{J.I.}}:
Linkprot: a database collecting information about biological links.
Nucleic acids research,
976
(2016)
\end{botherref}
\endbibitem

\bibitem{dabrowski2021topoly}
\begin{barticle}
\bauthor{\bsnm{Dabrowski-Tumanski}, \binits{P.}},
\bauthor{\bsnm{Rubach}, \binits{P.}},
\bauthor{\bsnm{Niemyska}, \binits{W.}},
\bauthor{\bsnm{Gren}, \binits{B.A.}},
\bauthor{\bsnm{Sulkowska}, \binits{J.I.}}:
\batitle{Topoly: {P}ython package to analyze topology of polymers}.
\bjtitle{Briefings in Bioinformatics}
\bvolume{22}(\bissue{3}),
\bfpage{196}
(\byear{2021})
\end{barticle}
\endbibitem

\bibitem{dawson1990homology}
\begin{barticle}
\bauthor{\bsnm{Dawson}, \binits{R.J.M.}}:
\batitle{Homology of weighted simplicial complexes}.
\bjtitle{Cahiers de Topologie et G{\'e}om{\'e}trie Diff{\'e}rentielle
  Cat{\'e}goriques}
\bvolume{31}(\bissue{3}),
\bfpage{229}--\blpage{243}
(\byear{1990})
\end{barticle}
\endbibitem

\bibitem{de2004topological}
\begin{bchapter}
\bauthor{\bsnm{De~Silva}, \binits{V.}},
\bauthor{\bsnm{Carlsson}, \binits{G.E.}}:
\bctitle{Topological estimation using witness complexes.}
In: \bbtitle{PBG},
pp. \bfpage{157}--\blpage{166}
(\byear{2004})
\end{bchapter}
\endbibitem

\bibitem{de2009persistent}
\begin{bchapter}
\bauthor{\bsnm{De~Silva}, \binits{V.}},
\bauthor{\bsnm{Vejdemo-Johansson}, \binits{M.}}:
\bctitle{Persistent cohomology and circular coordinates}.
In: \bbtitle{Proceedings of the Twenty-fifth Annual Symposium on Computational
  Geometry},
pp. \bfpage{227}--\blpage{236}
(\byear{2009})
\end{bchapter}
\endbibitem

\bibitem{de2011dualities}
\begin{barticle}
\bauthor{\bsnm{De~Silva}, \binits{V.}},
\bauthor{\bsnm{Morozov}, \binits{D.}},
\bauthor{\bsnm{Vejdemo-Johansson}, \binits{M.}}:
\batitle{Dualities in persistent (co)homology}.
\bjtitle{Inverse Problems}
\bvolume{27}(\bissue{12}),
\bfpage{124003}
(\byear{2011})
\end{barticle}
\endbibitem

\bibitem{delaunay1934sphere}
\begin{barticle}
\bauthor{\bsnm{Delaunay}, \binits{B.}}:
\batitle{Sur la sphere vide}.
\bjtitle{Izvestia Akademii Nauk SSSR}
\bvolume{7},
\bfpage{793}--\blpage{800}
(\byear{1934})
\end{barticle}
\endbibitem

\bibitem{deligne1974conjecture}
\begin{barticle}
\bauthor{\bsnm{Deligne}, \binits{P.}}:
\batitle{La conjecture de weil. i}.
\bjtitle{Publications Math{\'e}matiques de l'Institut des Hautes {\'E}tudes
  Scientifiques}
\bvolume{43},
\bfpage{273}--\blpage{307}
(\byear{1974})
\end{barticle}
\endbibitem

\bibitem{deligne1980conjecture}
\begin{barticle}
\bauthor{\bsnm{Deligne}, \binits{P.}}:
\batitle{La conjecture de weil: Ii}.
\bjtitle{Publications Math{\'e}matiques de l'IH{\'E}S}
\bvolume{52},
\bfpage{137}--\blpage{252}
(\byear{1980})
\end{barticle}
\endbibitem

\bibitem{desbrun2006discrete}
\begin{bchapter}
\bauthor{\bsnm{Desbrun}, \binits{M.}},
\bauthor{\bsnm{Kanso}, \binits{E.}},
\bauthor{\bsnm{Tong}, \binits{Y.}}:
\bctitle{Discrete differential forms for computational modeling}.
In: \bbtitle{ACM SIGGRAPH 2006 Courses},
pp. \bfpage{39}--\blpage{54}
(\byear{2006})
\end{bchapter}
\endbibitem

\bibitem{dey2021computing}
\begin{bchapter}
\bauthor{\bsnm{Dey}, \binits{T.K.}},
\bauthor{\bsnm{Hou}, \binits{T.}}:
\bctitle{Computing zigzag persistence on graphs in near-linear time}.
In: \bbtitle{37th International Symposium on Computational Geometry}
(\byear{2021})
\end{bchapter}
\endbibitem

\bibitem{dey2021updating}
\begin{botherref}
\oauthor{\bsnm{Dey}, \binits{T.K.}},
\oauthor{\bsnm{Hou}, \binits{T.}}:
Updating zigzag persistence and maintaining representatives over changing
  filtrations.
arXiv preprint arXiv:2112.02352
(2021)
\end{botherref}
\endbibitem

\bibitem{dey2022fast}
\begin{bchapter}
\bauthor{\bsnm{Dey}, \binits{T.K.}},
\bauthor{\bsnm{Hou}, \binits{T.}}:
\bctitle{Fast computation of zigzag persistence}.
In: \bbtitle{30th Annual European Symposium on Algorithms (ESA 2022)}.
\bsertitle{Leibniz International Proceedings in Informatics (LIPIcs)},
vol. \bseriesno{244}
(\byear{2022}).
\bcomment{pp.~43:1--43:15}
\end{bchapter}
\endbibitem

\bibitem{dey2024computing}
\begin{barticle}
\bauthor{\bsnm{Dey}, \binits{T.K.}},
\bauthor{\bsnm{Kim}, \binits{W.}},
\bauthor{\bsnm{M{\'e}moli}, \binits{F.}}:
\batitle{Computing generalized rank invariant for 2-parameter persistence
  modules via zigzag persistence and its applications}.
\bjtitle{Discrete \& Computational Geometry}
\bvolume{71}(\bissue{1}),
\bfpage{67}--\blpage{94}
(\byear{2024})
\end{barticle}
\endbibitem

\bibitem{dey2022efficient}
\begin{barticle}
\bauthor{\bsnm{Dey}, \binits{T.K.}},
\bauthor{\bsnm{Li}, \binits{T.}},
\bauthor{\bsnm{Wang}, \binits{Y.}}:
\batitle{An efficient algorithm for 1-dimensional (persistent) path homology}.
\bjtitle{Discrete \& Computational Geometry}
\bvolume{68}(\bissue{4}),
\bfpage{1102}--\blpage{1132}
(\byear{2022})
\end{barticle}
\endbibitem

\bibitem{dey2020persistence}
\begin{bchapter}
\bauthor{\bsnm{Dey}, \binits{T.K.}},
\bauthor{\bsnm{Mrozek}, \binits{M.}},
\bauthor{\bsnm{Slechta}, \binits{R.}}:
\bctitle{Persistence of the {Conley} index in combinatorial dynamical systems}.
In: \bbtitle{36th International Symposium on Computational Geometry (SoCG
  2020)}.
\bsertitle{Leibniz International Proceedings in Informatics (LIPIcs)},
vol. \bseriesno{164}
(\byear{2020}).
\bcomment{pp.~37:1--37:17}
\end{bchapter}
\endbibitem

\bibitem{dey2022persistence}
\begin{barticle}
\bauthor{\bsnm{Dey}, \binits{T.K.}},
\bauthor{\bsnm{Mrozek}, \binits{M.}},
\bauthor{\bsnm{Slechta}, \binits{R.}}:
\batitle{Persistence of {C}onley-{M}orse graphs in combinatorial dynamical
  systems}.
\bjtitle{SIAM Journal on Applied Dynamical Systems}
\bvolume{21}(\bissue{2}),
\bfpage{817}--\blpage{839}
(\byear{2022})
\end{barticle}
\endbibitem

\bibitem{dey2019persistent}
\begin{barticle}
\bauthor{\bsnm{Dey}, \binits{T.K.}},
\bauthor{\bsnm{Juda}, \binits{M.}},
\bauthor{\bsnm{Kapela}, \binits{T.}},
\bauthor{\bsnm{Kubica}, \binits{J.}},
\bauthor{\bsnm{Lipi\'nski}, \binits{M.}},
\bauthor{\bsnm{Mrozek}, \binits{M.}}:
\batitle{Persistent homology of {M}orse decompositions in combinatorial
  dynamics}.
\bjtitle{SIAM Journal on Applied Dynamical Systems}
\bvolume{18}(\bissue{1}),
\bfpage{510}--\blpage{530}
(\byear{2019})
\end{barticle}
\endbibitem

\bibitem{dey2022computational}
\begin{bbook}
\bauthor{\bsnm{Dey}, \binits{T.K.}},
\bauthor{\bsnm{Wang}, \binits{Y.}}:
\bbtitle{Computational Topology for Data Analysis}.
\bpublisher{Cambridge University Press},
\blocation{\hspace{0pt}}
(\byear{2022})
\end{bbook}
\endbibitem

\bibitem{di2024path}
\begin{barticle}
\bauthor{\bsnm{Di}, \binits{S.}},
\bauthor{\bsnm{Ivanov}, \binits{S.O.}},
\bauthor{\bsnm{Mukoseev}, \binits{L.}},
\bauthor{\bsnm{Zhang}, \binits{M.}}:
\batitle{On the path homology of cayley digraphs and covering digraphs}.
\bjtitle{Journal of Algebra}
\bvolume{653},
\bfpage{156}--\blpage{199}
(\byear{2024})
\end{barticle}
\endbibitem

\bibitem{di2015comparing}
\begin{bchapter}
\bauthor{\bsnm{Di~Fabio}, \binits{B.}},
\bauthor{\bsnm{Ferri}, \binits{M.}}:
\bctitle{Comparing persistence diagrams through complex vectors}.
In: \bbtitle{Image Analysis and Processing-ICIAP 2015: 18th International
  Conference, Genoa, Italy, September 7-11, 2015, Proceedings, Part I 18},
pp. \bfpage{294}--\blpage{305}
(\byear{2015}).
\bcomment{Springer}
\end{bchapter}
\endbibitem

\bibitem{dirac1928quantum}
\begin{barticle}
\bauthor{\bsnm{Dirac}, \binits{P.A.M.}}:
\batitle{The quantum theory of the electron}.
\bjtitle{Proceedings of the Royal Society of London. Series A, Containing
  Papers of a Mathematical and Physical Character}
\bvolume{117}(\bissue{778}),
\bfpage{610}--\blpage{624}
(\byear{1928})
\end{barticle}
\endbibitem

\bibitem{divol2019density}
\begin{barticle}
\bauthor{\bsnm{Divol}, \binits{V.}},
\bauthor{\bsnm{Chazal}, \binits{F.}}:
\batitle{The density of expected persistence diagrams and its kernel based
  estimation}.
\bjtitle{Journal of Computational Geometry}
\bvolume{10}(\bissue{2}),
\bfpage{1}--\blpage{31}
(\byear{2019}).
\doiurl{10.20382/jocg.v10i2a7}
\end{barticle}
\endbibitem

\bibitem{divol2021understanding}
\begin{barticle}
\bauthor{\bsnm{Divol}, \binits{V.}},
\bauthor{\bsnm{Lacombe}, \binits{T.}}:
\batitle{Understanding the topology and the geometry of the space of
  persistence diagrams via optimal partial transport}.
\bjtitle{Journal of Applied and Computational Topology}
\bvolume{5}(\bissue{1}),
\bfpage{1}--\blpage{53}
(\byear{2021}).
\doiurl{10.1007/s41468-020-00061-z}
\end{barticle}
\endbibitem

\bibitem{dlotko2016topological}
\begin{barticle}
\bauthor{\bsnm{D{\l}otko}, \binits{P.}},
\bauthor{\bsnm{Wanner}, \binits{T.}}:
\batitle{Topological microstructure analysis using persistence landscapes}.
\bjtitle{Physica D: Nonlinear Phenomena}
\bvolume{334},
\bfpage{60}--\blpage{81}
(\byear{2016})
\end{barticle}
\endbibitem

\bibitem{du2024multiscale}
\begin{barticle}
\bauthor{\bsnm{Du}, \binits{H.}},
\bauthor{\bsnm{Wei}, \binits{G.-W.}},
\bauthor{\bsnm{Hou}, \binits{T.}}:
\batitle{Multiscale topology in interactomic network: from transcriptome to
  antiaddiction drug repurposing}.
\bjtitle{Briefings in Bioinformatics}
\bvolume{25}(\bissue{2}),
\bfpage{054}
(\byear{2024})
\end{barticle}
\endbibitem

\bibitem{eckmann1944harmonische}
\begin{barticle}
\bauthor{\bsnm{Eckmann}, \binits{B.}}:
\batitle{Harmonische funktionen und randwertaufgaben in einem komplex}.
\bjtitle{Commentarii Mathematici Helvetici}
\bvolume{17}(\bissue{1}),
\bfpage{240}--\blpage{255}
(\byear{1944})
\end{barticle}
\endbibitem

\bibitem{edelsbrunner2002topological}
\begin{barticle}
\bauthor{\bsnm{Edelsbrunner}},
\bauthor{\bsnm{Letscher}},
\bauthor{\bsnm{Zomorodian}}:
\batitle{Topological persistence and simplification}.
\bjtitle{Discrete \& computational geometry}
\bvolume{28},
\bfpage{511}--\blpage{533}
(\byear{2002})
\end{barticle}
\endbibitem

\bibitem{edelsbrunner2011alpha}
\begin{bchapter}
\bauthor{\bsnm{Edelsbrunner}, \binits{H.}}:
\bctitle{Alpha shapes - a survey}.
In: \bbtitle{Tessellations in the Sciences: Virtues, Techniques and
  Applications of Geometric Tilings},
(\byear{2011})
\end{bchapter}
\endbibitem

\bibitem{edelsbrunner2010computational}
\begin{bbook}
\bauthor{\bsnm{Edelsbrunner}, \binits{H.}},
\bauthor{\bsnm{Harer}, \binits{J.}}:
\bbtitle{Computational Topology: An Introduction}.
\bpublisher{American Mathematical Soc.},
\blocation{\hspace{0pt}}
(\byear{2010})
\end{bbook}
\endbibitem

\bibitem{edelsbrunner2013persistent}
\begin{bbook}
\bauthor{\bsnm{Edelsbrunner}, \binits{H.}},
\bauthor{\bsnm{Morozov}, \binits{D.}}:
\bbtitle{Persistent Homology: Theory and Practice}.
\bpublisher{eScholarship, University of California},
\blocation{\hspace{0pt}}
(\byear{2013})
\end{bbook}
\endbibitem

\bibitem{edelsbrunner2008persistent}
\begin{barticle}
\bauthor{\bsnm{Edelsbrunner}, \binits{H.}},
\bauthor{\bsnm{Harer}, \binits{J.}}, \betal:
\batitle{Persistent homology - a survey}.
\bjtitle{Contemporary mathematics}
\bvolume{453}(\bissue{26}),
\bfpage{257}--\blpage{282}
(\byear{2008})
\end{barticle}
\endbibitem

\bibitem{einizade2025cosmos}
\begin{botherref}
\oauthor{\bsnm{Einzade}, \binits{A.}},
\oauthor{\bsnm{Thanou}, \binits{D.}},
\oauthor{\bsnm{Malliaros}, \binits{F.D.}},
\oauthor{\bsnm{Giraldo}, \binits{J.H.}}:
Cosmos: Continuous simplicial neural networks.
arXiv preprint arXiv:2503.12919
(2025)
\end{botherref}
\endbibitem

\bibitem{estrada2012path}
\begin{barticle}
\bauthor{\bsnm{Estrada}, \binits{E.}}:
\batitle{Path {Laplacian} matrices: introduction and application to the
  analysis of consensus in networks}.
\bjtitle{Linear algebra and its applications}
\bvolume{436}(\bissue{9}),
\bfpage{3373}--\blpage{3391}
(\byear{2012})
\end{barticle}
\endbibitem

\bibitem{fabri2009cgal}
\begin{bchapter}
\bauthor{\bsnm{Fabri}, \binits{A.}},
\bauthor{\bsnm{Pion}, \binits{S.}}:
\bctitle{{CGAL}: the computational geometry algorithms library}.
In: \bbtitle{Proceedings of the 17th ACM SIGSPATIAL International Conference on
  Advances in Geographic Information Systems},
pp. \bfpage{538}--\blpage{539}
(\byear{2009})
\end{bchapter}
\endbibitem

\bibitem{fasy2014confidence}
\begin{barticle}
\bauthor{\bsnm{Fasy}, \binits{B.T.}},
\bauthor{\bsnm{Lecci}, \binits{F.}},
\bauthor{\bsnm{Rinaldo}, \binits{A.}},
\bauthor{\bsnm{Wasserman}, \binits{L.}},
\bauthor{\bsnm{Balakrishnan}, \binits{S.}},
\bauthor{\bsnm{Singh}, \binits{A.}}:
\batitle{Confidence sets for persistence diagrams}.
\bjtitle{The Annals of Statistics}
\bvolume{42}(\bissue{6}),
\bfpage{2301}--\blpage{2339}
(\byear{2014}).
\doiurl{10.1214/14-AOS1252}
\end{barticle}
\endbibitem

\bibitem{fasy2014introduction}
\begin{botherref}
\oauthor{\bsnm{Fasy}, \binits{B.T.}},
\oauthor{\bsnm{Kim}, \binits{J.}},
\oauthor{\bsnm{Lecci}, \binits{F.}},
\oauthor{\bsnm{Maria}, \binits{C.}}:
Introduction to the {R} package {TDA}.
arXiv preprint arXiv:1411.1830
(2014)
\end{botherref}
\endbibitem

\bibitem{feng2024mayer}
\begin{botherref}
\oauthor{\bsnm{Feng}, \binits{H.}},
\oauthor{\bsnm{Shen}, \binits{L.}},
\oauthor{\bsnm{Liu}, \binits{J.}},
\oauthor{\bsnm{Wei}, \binits{G.-W.}}:
Mayer-homology learning prediction of protein-ligand binding affinities.
Journal of Computational Biophysics and Chemistry,
1--14
(2024)
\end{botherref}
\endbibitem

\bibitem{feng2024hypernetwork}
\begin{barticle}
\bauthor{\bsnm{Feng}, \binits{L.}},
\bauthor{\bsnm{Gong}, \binits{H.}},
\bauthor{\bsnm{Zhang}, \binits{S.}},
\bauthor{\bsnm{Liu}, \binits{X.}},
\bauthor{\bsnm{Wang}, \binits{Y.}},
\bauthor{\bsnm{Che}, \binits{J.}},
\bauthor{\bsnm{Dong}, \binits{A.}},
\bauthor{\bsnm{Griffin}, \binits{C.H.}},
\bauthor{\bsnm{Gragnoli}, \binits{C.}},
\bauthor{\bsnm{Wu}, \binits{J.}},
\bauthor{\bsnm{Yau}, \binits{S.-T.}},
\bauthor{\bsnm{Wu}, \binits{R.}}:
\batitle{Hypernetwork modeling and topology of high-order interactions for
  complex systems}.
\bjtitle{Proceedings of the National Academy of Sciences}
\bvolume{121}(\bissue{40}),
\bfpage{2412220121}
(\byear{2024})
\end{barticle}
\endbibitem

\bibitem{feng2025network}
\begin{barticle}
\bauthor{\bsnm{Feng}, \binits{L.}},
\bauthor{\bsnm{Yang}, \binits{D.}},
\bauthor{\bsnm{Wu}, \binits{S.}},
\bauthor{\bsnm{Xue}, \binits{C.}},
\bauthor{\bsnm{Sang}, \binits{M.}},
\bauthor{\bsnm{Liu}, \binits{X.}},
\bauthor{\bsnm{Che}, \binits{J.}},
\bauthor{\bsnm{Wu}, \binits{J.}},
\bauthor{\bsnm{Gragnoli}, \binits{C.}},
\bauthor{\bsnm{Griffin}, \binits{C.}},
\bauthor{\bsnm{Wang}, \binits{C.}},
\bauthor{\bsnm{Yau}, \binits{S.-T.}},
\bauthor{\bsnm{Wu}, \binits{R.}}:
\batitle{Network modeling and topology of aging}.
\bjtitle{Physics reports}
\bvolume{1101},
\bfpage{1}--\blpage{65}
(\byear{2025})
\end{barticle}
\endbibitem

\bibitem{ferri1999representing}
\begin{barticle}
\bauthor{\bsnm{Ferri}, \binits{M.}},
\bauthor{\bsnm{Landi}, \binits{C.}}:
\batitle{Representing size functions by complex polynomials}.
\bjtitle{Proc. Math. Met. in Pattern Recognition}
\bvolume{9},
\bfpage{16}--\blpage{19}
(\byear{1999})
\end{barticle}
\endbibitem

\bibitem{friedrichs1955differential}
\begin{barticle}
\bauthor{\bsnm{Friedrichs}, \binits{K.O.}}:
\batitle{Differential forms on {R}iemannian manifolds}.
\bjtitle{Communications on Pure and Applied Mathematics}
\bvolume{8}(\bissue{4}),
\bfpage{551}--\blpage{590}
(\byear{1955})
\end{barticle}
\endbibitem

\bibitem{frosini1999size}
\begin{barticle}
\bauthor{\bsnm{Frosini}, \binits{P.}},
\bauthor{\bsnm{Mulazzani}, \binits{M.}}:
\batitle{Size homotopy groups for computation of natural size distances}.
\bjtitle{Bulletin of the Belgian Mathematical Society-Simon Stevin}
\bvolume{6}(\bissue{3}),
\bfpage{455}--\blpage{464}
(\byear{1999})
\end{barticle}
\endbibitem

\bibitem{gameiro2015topological}
\begin{barticle}
\bauthor{\bsnm{Gameiro}, \binits{M.}},
\bauthor{\bsnm{Hiraoka}, \binits{Y.}},
\bauthor{\bsnm{Izumi}, \binits{S.}},
\bauthor{\bsnm{Kramar}, \binits{M.}},
\bauthor{\bsnm{Mischaikow}, \binits{K.}},
\bauthor{\bsnm{Nanda}, \binits{V.}}:
\batitle{A topological measurement of protein compressibility}.
\bjtitle{Japan Journal of Industrial and Applied Mathematics}
\bvolume{32},
\bfpage{1}--\blpage{17}
(\byear{2015})
\end{barticle}
\endbibitem

\bibitem{Gauss1877}
\begin{bbook}
\bauthor{\bsnm{Gauss}, \binits{C.F.}}:
\bbtitle{Zur Mathematischen Theorie der Electrodynamischen Wirkungen},
pp. \bfpage{601}--\blpage{630}.
\bpublisher{Springer},
\blocation{Berlin, Heidelberg}
(\byear{1877})
\end{bbook}
\endbibitem

\bibitem{ge2011data}
\begin{botherref}
\oauthor{\bsnm{Ge}, \binits{X.}},
\oauthor{\bsnm{Safa}, \binits{I.}},
\oauthor{\bsnm{Belkin}, \binits{M.}},
\oauthor{\bsnm{Wang}, \binits{Y.}}:
Data skeletonization via {Reeb} graphs.
Advances in neural information processing systems
\textbf{24}
(2011)
\end{botherref}
\endbibitem

\bibitem{geuzaine2009gmsh}
\begin{barticle}
\bauthor{\bsnm{Geuzaine}, \binits{C.}},
\bauthor{\bsnm{Remacle}, \binits{J.-F.}}:
\batitle{Gmsh: a 3-{D} finite element mesh generator with built-in pre- and
  post-processing facilities}.
\bjtitle{International journal for numerical methods in engineering}
\bvolume{79}(\bissue{11}),
\bfpage{1309}--\blpage{1331}
(\byear{2009})
\end{barticle}
\endbibitem

\bibitem{ghrist2008barcodes}
\begin{barticle}
\bauthor{\bsnm{Ghrist}, \binits{R.}}:
\batitle{Barcodes: the persistent topology of data}.
\bjtitle{Bulletin of the American Mathematical Society}
\bvolume{45}(\bissue{1}),
\bfpage{61}--\blpage{75}
(\byear{2008})
\end{barticle}
\endbibitem

\bibitem{ghrist2014elementary}
\begin{bbook}
\bauthor{\bsnm{Ghrist}, \binits{R.W.}}:
\bbtitle{Elementary Applied Topology}
vol. \bseriesno{1}.
\bpublisher{Createspace Seattle},
\blocation{\hspace{0pt}}
(\byear{2014})
\end{bbook}
\endbibitem

\bibitem{goldberg2002combinatorial}
\begin{botherref}
\oauthor{\bsnm{Goldberg}, \binits{T.E.}}:
Combinatorial {L}aplacians of simplicial complexes.
PhD thesis,
Citeseer
(2002)
\end{botherref}
\endbibitem

\bibitem{gong2024topological}
\begin{botherref}
\oauthor{\bsnm{Gong}, \binits{H.}},
\oauthor{\bsnm{Wang}, \binits{H.}},
\oauthor{\bsnm{Wang}, \binits{Y.}},
\oauthor{\bsnm{Zhang}, \binits{S.}},
\oauthor{\bsnm{Liu}, \binits{X.}},
\oauthor{\bsnm{Che}, \binits{J.}},
\oauthor{\bsnm{Wu}, \binits{S.}},
\oauthor{\bsnm{Wu}, \binits{J.}},
\oauthor{\bsnm{Sun}, \binits{X.}},
\oauthor{\bsnm{Zhang}, \binits{S.}},
\oauthor{\bsnm{Yau}, \binits{S.-T.}},
\oauthor{\bsnm{Wu}, \binits{R.}}:
Topological change of soil microbiota networks for forest resilience under
  global warming.
Physics of Life Reviews
(2024)
\end{botherref}
\endbibitem

\bibitem{grbic2022aspects}
\begin{barticle}
\bauthor{\bsnm{Grbi{\'c}}, \binits{J.}},
\bauthor{\bsnm{Wu}, \binits{J.}},
\bauthor{\bsnm{Xia}, \binits{K.}},
\bauthor{\bsnm{Wei}, \binits{G.-W.}}:
\batitle{Aspects of topological approaches for data science}.
\bjtitle{Foundations of data science (Springfield, Mo.)}
\bvolume{4}(\bissue{2}),
\bfpage{165}
(\byear{2022})
\end{barticle}
\endbibitem

\bibitem{grigor2012homologies}
\begin{botherref}
\oauthor{\bsnm{Grigor'yan}, \binits{A.}},
\oauthor{\bsnm{Lin}, \binits{Y.}},
\oauthor{\bsnm{Muranov}, \binits{Y.}},
\oauthor{\bsnm{Yau}, \binits{S.-T.}}:
Homologies of path complexes and digraphs.
arXiv preprint arXiv:1207.2834
(2012)
\end{botherref}
\endbibitem

\bibitem{grigor2022advances}
\begin{barticle}
\bauthor{\bsnm{Grigor’yan}, \binits{A.}}:
\batitle{Advances in path homology theory of digraphs}.
\bjtitle{Notices of the International Consortium of Chinese Mathematicians}
\bvolume{10}(\bissue{2}),
\bfpage{61}--\blpage{124}
(\byear{2022})
\end{barticle}
\endbibitem

\bibitem{grigor2016cohomology}
\begin{barticle}
\bauthor{\bsnm{Grigor’yan}, \binits{A.}},
\bauthor{\bsnm{Muranov}, \binits{Y.}},
\bauthor{\bsnm{Yau}, \binits{S.-T.}}:
\batitle{On a cohomology of digraphs and hochschild cohomology}.
\bjtitle{Journal of Homotopy and Related Structures}
\bvolume{11}(\bissue{2}),
\bfpage{209}--\blpage{230}
(\byear{2016})
\end{barticle}
\endbibitem

\bibitem{grigor2023homotopy}
\begin{barticle}
\bauthor{\bsnm{Grigor’yan}, \binits{A.}},
\bauthor{\bsnm{Lin}, \binits{Y.}},
\bauthor{\bsnm{Muranov}, \binits{Y.}},
\bauthor{\bsnm{Yau}, \binits{S.-T.}}:
\batitle{Homotopy theory for digraphs}.
\bjtitle{Pure and Applied Mathematics Quarterly}
\bvolume{10}(\bissue{4}),
\bfpage{619}--\blpage{674}
(\byear{2023})
\end{barticle}
\endbibitem

\bibitem{gulen2023orthogonal}
\begin{botherref}
\oauthor{\bsnm{G{\"u}len}, \binits{A.B.}},
\oauthor{\bsnm{M{\'e}moli}, \binits{F.}},
\oauthor{\bsnm{Wan}, \binits{Z.}}:
Orthogonal m\"obius inversion and {Grassmannian} persistence diagrams.
arXiv preprint arXiv:2311.06870
(2023)
\end{botherref}
\endbibitem

\bibitem{gulen2025grassmannian}
\begin{botherref}
\oauthor{\bsnm{G{\"u}len}, \binits{A.B.}},
\oauthor{\bsnm{M{\'e}moli}, \binits{F.}},
\oauthor{\bsnm{Wan}, \binits{Z.}}:
Grassmannian persistence diagrams: Special properties in the 1-parameter
  setting.
arXiv preprint arXiv:2504.06077
(2025)
\end{botherref}
\endbibitem

\bibitem{gulen2023generalization}
\begin{bchapter}
\bauthor{\bsnm{G\"{u}len}, \binits{A.B.}},
\bauthor{\bsnm{M\'{e}moli}, \binits{F.}},
\bauthor{\bsnm{Wan}, \binits{Z.}},
\bauthor{\bsnm{Wang}, \binits{Y.}}:
\bctitle{{A Generalization of the Persistent {Laplacian} to Simplicial Maps}}.
In: \bbtitle{39th International Symposium on Computational Geometry (SoCG
  2023)}.
\bsertitle{Leibniz International Proceedings in Informatics (LIPIcs)},
vol. \bseriesno{258}
(\byear{2023}).
\bcomment{pp.~37:1--37:17}
\end{bchapter}
\endbibitem

\bibitem{gundert2014higher}
\begin{bchapter}
\bauthor{\bsnm{Gundert}, \binits{A.}},
\bauthor{\bsnm{Szedl{\'a}k}, \binits{M.}}:
\bctitle{Higher dimensional {C}heeger inequalities}.
In: \bbtitle{Proceedings of the Thirtieth Annual Symposium on Computational
  Geometry},
pp. \bfpage{181}--\blpage{188}
(\byear{2014})
\end{bchapter}
\endbibitem

\bibitem{gunther2012efficient}
\begin{barticle}
\bauthor{\bsnm{G{\"u}nther}, \binits{D.}},
\bauthor{\bsnm{Reininghaus}, \binits{J.}},
\bauthor{\bsnm{Wagner}, \binits{H.}},
\bauthor{\bsnm{Hotz}, \binits{I.}}:
\batitle{Efficient computation of {3D} {M}orse--{S}male complexes and
  persistent homology using discrete {M}orse theory}.
\bjtitle{The Visual Computer}
\bvolume{28},
\bfpage{959}--\blpage{969}
(\byear{2012})
\end{barticle}
\endbibitem

\bibitem{hang2023correspondence}
\begin{barticle}
\bauthor{\bsnm{Hang}, \binits{H.}},
\bauthor{\bsnm{Mio}, \binits{W.}}:
\batitle{Correspondence modules and persistence sheaves: a unifying perspective
  on one-parameter persistent homology}.
\bjtitle{Japan Journal of Industrial and Applied Mathematics}
\bvolume{40}(\bissue{1}),
\bfpage{41}--\blpage{93}
(\byear{2023})
\end{barticle}
\endbibitem

\bibitem{hansen2020laplacians}
\begin{botherref}
\oauthor{\bsnm{Hansen}, \binits{J.}}:
Laplacians of cellular sheaves: Theory and applications.
PhD thesis,
University of Pennsylvania
(2020)
\end{botherref}
\endbibitem

\bibitem{hansen2019toward}
\begin{barticle}
\bauthor{\bsnm{Hansen}, \binits{J.}},
\bauthor{\bsnm{Ghrist}, \binits{R.}}:
\batitle{Toward a spectral theory of cellular sheaves}.
\bjtitle{Journal of Applied and Computational Topology}
\bvolume{3}(\bissue{4}),
\bfpage{315}--\blpage{358}
(\byear{2019})
\end{barticle}
\endbibitem

\bibitem{hayakawa2022quantum}
\begin{barticle}
\bauthor{\bsnm{Hayakawa}, \binits{R.}}:
\batitle{Quantum algorithm for persistent {Betti} numbers and topological data
  analysis}.
\bjtitle{Quantum}
\bvolume{6},
\bfpage{873}
(\byear{2022})
\end{barticle}
\endbibitem

\bibitem{he2025multi}
\begin{barticle}
\bauthor{\bsnm{He}, \binits{Y.}},
\bauthor{\bsnm{Liu}, \binits{J.}}:
\batitle{Multi-scale hochschild spectral analysis on graph data}.
\bjtitle{AIMS Mathematics}
\bvolume{10}(\bissue{1}),
\bfpage{1384}--\blpage{1406}
(\byear{2025})
\end{barticle}
\endbibitem

\bibitem{hensel2021survey}
\begin{barticle}
\bauthor{\bsnm{Hensel}, \binits{F.}},
\bauthor{\bsnm{Moor}, \binits{M.}},
\bauthor{\bsnm{Rieck}, \binits{B.}}:
\batitle{A survey of topological machine learning methods}.
\bjtitle{Frontiers in Artificial Intelligence}
\bvolume{4},
\bfpage{681108}
(\byear{2021})
\end{barticle}
\endbibitem

\bibitem{henselmanghristl6}
\begin{botherref}
\oauthor{\bsnm{Henselman}, \binits{G.}},
\oauthor{\bsnm{Ghrist}, \binits{R.}}:
Matroid filtrations and computational persistent homology.
arXiv preprint arXiv:1606.00199
(2016)
\end{botherref}
\endbibitem

\bibitem{hernandez2025persistence}
\begin{botherref}
\oauthor{\bsnm{Hern{\'a}ndez-Garc{\'\i}a}, \binits{P.}},
\oauthor{\bsnm{Serrano}, \binits{D.H.}},
\oauthor{\bsnm{G{\'o}mez}, \binits{D.S.}}:
From persistence to resilience: New betti numbers for analyzing robustness in
  simplicial complex networks.
arXiv preprint arXiv:2505.10467
(2025)
\end{botherref}
\endbibitem

\bibitem{holme2012temporal}
\begin{barticle}
\bauthor{\bsnm{Holme}, \binits{P.}},
\bauthor{\bsnm{Saram{\"a}ki}, \binits{J.}}:
\batitle{Temporal networks}.
\bjtitle{Physics reports}
\bvolume{519}(\bissue{3}),
\bfpage{97}--\blpage{125}
(\byear{2012})
\end{barticle}
\endbibitem

\bibitem{horak2013spectra}
\begin{barticle}
\bauthor{\bsnm{Horak}, \binits{D.}},
\bauthor{\bsnm{Jost}, \binits{J.}}:
\batitle{Spectra of combinatorial {L}aplace operators on simplicial complexes}.
\bjtitle{Advances in Mathematics}
\bvolume{244},
\bfpage{303}--\blpage{336}
(\byear{2013})
\end{barticle}
\endbibitem

\bibitem{horak2009persistent}
\begin{barticle}
\bauthor{\bsnm{Horak}, \binits{D.}},
\bauthor{\bsnm{Maleti{\'c}}, \binits{S.}},
\bauthor{\bsnm{Rajkovi{\'c}}, \binits{M.}}:
\batitle{Persistent homology of complex networks}.
\bjtitle{Journal of Statistical Mechanics: Theory and Experiment}
\bvolume{2009}(\bissue{03}),
\bfpage{03034}
(\byear{2009})
\end{barticle}
\endbibitem

\bibitem{hozumi2024revealing}
\begin{botherref}
\oauthor{\bsnm{Hozumi}, \binits{Y.}},
\oauthor{\bsnm{Wei}, \binits{G.-W.}}:
Revealing the shape of genome space via k-mer topology.
arXiv preprint arXiv:2412.20202
(2024)
\end{botherref}
\endbibitem

\bibitem{ivanov2024simplicial}
\begin{barticle}
\bauthor{\bsnm{Ivanov}, \binits{S.O.}},
\bauthor{\bsnm{Pavutnitskiy}, \binits{F.}}:
\batitle{Simplicial approach to path homology of quivers, marked categories,
  groups and algebras}.
\bjtitle{Journal of the London Mathematical Society}
\bvolume{109}(\bissue{1}),
\bfpage{12812}
(\byear{2024})
\end{barticle}
\endbibitem

\bibitem{gudhi:TangentialComplex}
\begin{bchapter}
\bauthor{\bsnm{Jamin}, \binits{C.}}:
\bctitle{Tangential complex}.
In: \bbtitle{GUDHI User and Reference Manual},
(\byear{2025})
\end{bchapter}
\endbibitem

\bibitem{jamroz2015knotprot}
\begin{barticle}
\bauthor{\bsnm{Jamroz}, \binits{M.}},
\bauthor{\bsnm{Niemyska}, \binits{W.}},
\bauthor{\bsnm{Rawdon}, \binits{E.J.}},
\bauthor{\bsnm{Stasiak}, \binits{A.}},
\bauthor{\bsnm{Millett}, \binits{K.C.}},
\bauthor{\bsnm{Su{\l}kowski}, \binits{P.}},
\bauthor{\bsnm{Sulkowska}, \binits{J.I.}}:
\batitle{Knotprot: a database of proteins with knots and slipknots}.
\bjtitle{Nucleic acids research}
\bvolume{43}(\bissue{D1}),
\bfpage{306}--\blpage{314}
(\byear{2015})
\end{barticle}
\endbibitem

\bibitem{jiang2021topological}
\begin{barticle}
\bauthor{\bsnm{Jiang}, \binits{Y.}},
\bauthor{\bsnm{Chen}, \binits{D.}},
\bauthor{\bsnm{Chen}, \binits{X.}},
\bauthor{\bsnm{Li}, \binits{T.}},
\bauthor{\bsnm{Wei}, \binits{G.-W.}},
\bauthor{\bsnm{Pan}, \binits{F.}}:
\batitle{Topological representations of crystalline compounds for the
  machine-learning prediction of materials properties}.
\bjtitle{npj computational materials}
\bvolume{7}(\bissue{1}),
\bfpage{28}
(\byear{2021})
\end{barticle}
\endbibitem

\bibitem{jones2025khovanov}
\begin{barticle}
\bauthor{\bsnm{Jones}, \binits{B.}},
\bauthor{\bsnm{Wei}, \binits{G.-W.}}:
\batitle{Khovanov {Laplacian} and {Khovanov} {Dirac} for knots and links}.
\bjtitle{Journal of Physics: Complexity}
\bvolume{6}(\bissue{2}),
\bfpage{025014}
(\byear{2025})
\end{barticle}
\endbibitem

\bibitem{jones2025persistent}
\begin{barticle}
\bauthor{\bsnm{Jones}, \binits{B.}},
\bauthor{\bsnm{Wei}, \binits{G.-W.}}:
\batitle{Persistent directed flag {L}aplacian}.
\bjtitle{Foundations of Data Science}
\bvolume{7}(\bissue{3}),
\bfpage{737}--\blpage{758}
(\byear{2025})
\end{barticle}
\endbibitem

\bibitem{jones1997polynomial}
\begin{bchapter}
\bauthor{\bsnm{Jones}, \binits{V.F.}}:
\bctitle{A polynomial invariant for knots via von {N}eumann algebras}.
In: \bbtitle{Fields Medallists' Lectures},
pp. \bfpage{448}--\blpage{458}.
\bpublisher{World Scientific},
\blocation{\hspace{0pt}}
(\byear{1997})
\end{bchapter}
\endbibitem

\bibitem{jonsson2008simplicial}
\begin{bbook}
\bauthor{\bsnm{Jonsson}, \binits{J.}}:
\bbtitle{Simplicial Complexes of Graphs}
vol. \bseriesno{1928}.
\bpublisher{Springer},
\blocation{\hspace{0pt}}
(\byear{2008})
\end{bbook}
\endbibitem

\bibitem{kaji2020cubical}
\begin{botherref}
\oauthor{\bsnm{Kaji}, \binits{S.}},
\oauthor{\bsnm{Sudo}, \binits{T.}},
\oauthor{\bsnm{Ahara}, \binits{K.}}:
Cubical {R}ipser: Software for computing persistent homology of image and
  volume data.
arXiv preprint arXiv:2005.12692
(2020)
\end{botherref}
\endbibitem

\bibitem{kalivsnik2019tropical}
\begin{barticle}
\bauthor{\bsnm{Kali{\v{s}}nik}, \binits{S.}}:
\batitle{Tropical coordinates on the space of persistence barcodes}.
\bjtitle{Foundations of Computational Mathematics}
\bvolume{19}(\bissue{1}),
\bfpage{101}--\blpage{129}
(\byear{2019})
\end{barticle}
\endbibitem

\bibitem{kannan2019persistent}
\begin{barticle}
\bauthor{\bsnm{Kannan}, \binits{H.}},
\bauthor{\bsnm{Saucan}, \binits{E.}},
\bauthor{\bsnm{Roy}, \binits{I.}},
\bauthor{\bsnm{Samal}, \binits{A.}}:
\batitle{Persistent homology of unweighted complex networks via discrete
  {M}orse theory}.
\bjtitle{Scientific reports}
\bvolume{9}(\bissue{1}),
\bfpage{13817}
(\byear{2019})
\end{barticle}
\endbibitem

\bibitem{karaguler2021survey}
\begin{botherref}
\oauthor{\bsnm{Karag{\"u}ler}, \binits{D.}}:
A survey on multidimensional persistence theory.
Master's thesis,
Middle East Technical University (Turkey)
(2021)
\end{botherref}
\endbibitem

\bibitem{karan2021time}
\begin{barticle}
\bauthor{\bsnm{Karan}, \binits{A.}},
\bauthor{\bsnm{Kaygun}, \binits{A.}}:
\batitle{Time series classification via topological data analysis}.
\bjtitle{Expert Systems with Applications}
\bvolume{183},
\bfpage{115326}
(\byear{2021}).
\doiurl{10.1016/j.eswa.2021.115326}
\end{barticle}
\endbibitem

\bibitem{keros2023spectral}
\begin{bchapter}
\bauthor{\bsnm{Keros}, \binits{A.}},
\bauthor{\bsnm{Subr}, \binits{K.}}:
\bctitle{Spectral coarsening with {H}odge {L}aplacians}.
In: \bbtitle{ACM SIGGRAPH 2023 Conference Proceedings},
pp. \bfpage{1}--\blpage{11}
(\byear{2023})
\end{bchapter}
\endbibitem

\bibitem{khovanov2000categorification}
\begin{botherref}
\oauthor{\bsnm{Khovanov}, \binits{M.}}:
A categorification of the {J}ones polynomial
(2000)
\end{botherref}
\endbibitem

\bibitem{kim2020pllay}
\begin{barticle}
\bauthor{\bsnm{Kim}, \binits{K.}},
\bauthor{\bsnm{Kim}, \binits{J.}},
\bauthor{\bsnm{Zaheer}, \binits{M.}},
\bauthor{\bsnm{Kim}, \binits{J.}},
\bauthor{\bsnm{Chazal}, \binits{F.}},
\bauthor{\bsnm{Wasserman}, \binits{L.}}:
\batitle{Pllay: Efficient topological layer based on persistent landscapes}.
\bjtitle{Advances in Neural Information Processing Systems}
\bvolume{33},
\bfpage{15965}--\blpage{15977}
(\byear{2020})
\end{barticle}
\endbibitem

\bibitem{kim2021generalized}
\begin{barticle}
\bauthor{\bsnm{Kim}, \binits{W.}},
\bauthor{\bsnm{M{\'e}moli}, \binits{F.}}:
\batitle{Generalized persistence diagrams for persistence modules over posets}.
\bjtitle{Journal of Applied and Computational Topology}
\bvolume{5}(\bissue{4}),
\bfpage{533}--\blpage{581}
(\byear{2021})
\end{barticle}
\endbibitem

\bibitem{kim2023persistence}
\begin{botherref}
\oauthor{\bsnm{Kim}, \binits{W.}},
\oauthor{\bsnm{M{\'e}moli}, \binits{F.}}:
Persistence over posets.
Notices of the American Mathematical Society
\textbf{70}(08)
(2023)
\end{botherref}
\endbibitem

\bibitem{kirchhoff1847ueber}
\begin{barticle}
\bauthor{\bsnm{Kirchhoff}, \binits{G.}}:
\batitle{Ueber die aufl{\"o}sung der gleichungen, auf welche man bei der
  untersuchung der linearen vertheilung galvanischer str{\"o}me gef{\"u}hrt
  wird}.
\bjtitle{Annalen der Physik}
\bvolume{148}(\bissue{12}),
\bfpage{497}--\blpage{508}
(\byear{1847})
\end{barticle}
\endbibitem

\bibitem{knill2018cohomology}
\begin{botherref}
\oauthor{\bsnm{Knill}, \binits{O.}}:
The cohomology for {Wu} characteristics.
arXiv preprint arXiv:1803.06788
(2018)
\end{botherref}
\endbibitem

\bibitem{kovacev2016using}
\begin{barticle}
\bauthor{\bsnm{Kovacev-Nikolic}, \binits{V.}},
\bauthor{\bsnm{Bubenik}, \binits{P.}},
\bauthor{\bsnm{Nikoli{\'c}}, \binits{D.}},
\bauthor{\bsnm{Heo}, \binits{G.}}:
\batitle{Using persistent homology and dynamical distances to analyze protein
  binding}.
\bjtitle{Statistical applications in genetics and molecular biology}
\bvolume{15}(\bissue{1}),
\bfpage{19}--\blpage{38}
(\byear{2016})
\end{barticle}
\endbibitem

\bibitem{krishnagopal2023topology}
\begin{barticle}
\bauthor{\bsnm{Krishnagopal}, \binits{S.}},
\bauthor{\bsnm{Bianconi}, \binits{G.}}:
\batitle{Topology and dynamics of higher-order multiplex networks}.
\bjtitle{Chaos, Solitons \& Fractals}
\bvolume{177},
\bfpage{114296}
(\byear{2023})
\end{barticle}
\endbibitem

\bibitem{kusano2016persistence}
\begin{bchapter}
\bauthor{\bsnm{Kusano}, \binits{G.}},
\bauthor{\bsnm{Hiraoka}, \binits{Y.}},
\bauthor{\bsnm{Fukumizu}, \binits{K.}}:
\bctitle{Persistence weighted {Gaussian} kernel for topological data analysis}.
In: \bbtitle{International Conference on Machine Learning},
pp. \bfpage{2004}--\blpage{2013}
(\byear{2016}).
\bcomment{PMLR}
\end{bchapter}
\endbibitem

\bibitem{ladyzhenskaya1969mathematical}
\begin{bbook}
\bauthor{\bsnm{Ladyzhenskaya}, \binits{O.A.}}:
\bbtitle{The Mathematical Theory of Viscous Incompressible Flow}.
\bpublisher{Gordon \& Breach},
\blocation{\hspace{0pt}}
(\byear{1969})
\end{bbook}
\endbibitem

\bibitem{le2025persistent}
\begin{barticle}
\bauthor{\bsnm{Le}, \binits{M.Q.}},
\bauthor{\bsnm{Taylor}, \binits{D.}}:
\batitle{Persistent homology with k-nearest-neighbor filtrations reveals
  topological convergence of pagerank}.
\bjtitle{Foundations of Data Science}
\bvolume{7}(\bissue{2}),
\bfpage{536}--\blpage{567}
(\byear{2025})
\end{barticle}
\endbibitem

\bibitem{le2018persistence}
\begin{botherref}
\oauthor{\bsnm{Le}, \binits{T.}},
\oauthor{\bsnm{Yamada}, \binits{M.}}:
Persistence {F}isher kernel: A {R}iemannian manifold kernel for persistence
  diagrams.
Advances in neural information processing systems
\textbf{31}
(2018)
\end{botherref}
\endbibitem

\bibitem{lee2017quantifying}
\begin{barticle}
\bauthor{\bsnm{Lee}, \binits{Y.}},
\bauthor{\bsnm{Barthel}, \binits{S.D.}},
\bauthor{\bsnm{D{\l}otko}, \binits{P.}},
\bauthor{\bsnm{Moosavi}, \binits{S.M.}},
\bauthor{\bsnm{Hess}, \binits{K.}},
\bauthor{\bsnm{Smit}, \binits{B.}}:
\batitle{Quantifying similarity of pore-geometry in nanoporous materials}.
\bjtitle{Nature communications}
\bvolume{8}(\bissue{1}),
\bfpage{1}--\blpage{8}
(\byear{2017})
\end{barticle}
\endbibitem

\bibitem{Leray1946sheaf}
\begin{botherref}
\oauthor{\bsnm{Leray}, \binits{J.}}:
L'anneau d'homologie d'une repr\'esentation.
C.R. Acad. Sci. Paris
(222),
1366--1368
(1946)
\end{botherref}
\endbibitem

\bibitem{li2024singular}
\begin{botherref}
\oauthor{\bsnm{Li}, \binits{J.}},
\oauthor{\bsnm{Muranov}, \binits{Y.}},
\oauthor{\bsnm{Wu}, \binits{J.}},
\oauthor{\bsnm{Yau}, \binits{S.-T.}}:
On singular homology theories of digraphs and quivers
(2024)
\end{botherref}
\endbibitem

\bibitem{li2024primitive}
\begin{botherref}
\oauthor{\bsnm{Li}, \binits{J.}},
\oauthor{\bsnm{Muranov}, \binits{Y.}},
\oauthor{\bsnm{Wu}, \binits{J.}},
\oauthor{\bsnm{Yau}, \binits{S.-T.}}:
Primitive path homology.
arXiv preprint arXiv:2411.18955
(2024)
\end{botherref}
\endbibitem

\bibitem{li2017twisted}
\begin{botherref}
\oauthor{\bsnm{LI}, \binits{J.}},
\oauthor{\bsnm{VERSHININ}, \binits{V.}},
\oauthor{\bsnm{WU}, \binits{J.}}:
Twisted simplicial groups and twisted homology of categories.
Homology, Homotopy \& Applications
\textbf{19}(2)
(2017)
\end{botherref}
\endbibitem

\bibitem{liang1994knots}
\begin{barticle}
\bauthor{\bsnm{Liang}, \binits{C.}},
\bauthor{\bsnm{Mislow}, \binits{K.}}:
\batitle{Knots in proteins}.
\bjtitle{Journal of the American Chemical Society}
\bvolume{116}(\bissue{24}),
\bfpage{11189}--\blpage{11190}
(\byear{1994})
\end{barticle}
\endbibitem

\bibitem{Lieutier2014HarmonicForms}
\begin{botherref}
\oauthor{\bsnm{Lieutier}, \binits{A.}}:
Talk: Persistent harmonic forms
(2014)
\end{botherref}
\endbibitem

\bibitem{lim2020hodge}
\begin{barticle}
\bauthor{\bsnm{Lim}, \binits{L.-H.}}:
\batitle{Hodge {L}aplacians on graphs}.
\bjtitle{Siam Review}
\bvolume{62}(\bissue{3}),
\bfpage{685}--\blpage{715}
(\byear{2020})
\end{barticle}
\endbibitem

\bibitem{lin2019weighted}
\begin{botherref}
\oauthor{\bsnm{Lin}, \binits{Y.}},
\oauthor{\bsnm{Ren}, \binits{S.}},
\oauthor{\bsnm{Wang}, \binits{C.}},
\oauthor{\bsnm{Wu}, \binits{J.}}:
Weighted path homology of weighted digraphs and persistence.
arXiv preprint arXiv:1910.09891
(2019)
\end{botherref}
\endbibitem

\bibitem{liu2023interaction}
\begin{botherref}
\oauthor{\bsnm{Liu}, \binits{J.}},
\oauthor{\bsnm{Chen}, \binits{D.}},
\oauthor{\bsnm{Wei}, \binits{G.-W.}}:
Interaction homotopy and interaction homology.
arXiv preprint arXiv:2311.16322
(2023)
\end{botherref}
\endbibitem

\bibitem{liu2024persistent}
\begin{botherref}
\oauthor{\bsnm{Liu}, \binits{J.}},
\oauthor{\bsnm{Chen}, \binits{D.}},
\oauthor{\bsnm{Wei}, \binits{G.-W.}}:
Persistent interaction topology in data analysis.
arXiv preprint arXiv:2404.11799
(2024)
\end{botherref}
\endbibitem

\bibitem{liu2024algebraic}
\begin{barticle}
\bauthor{\bsnm{Liu}, \binits{J.}},
\bauthor{\bsnm{Li}, \binits{J.}},
\bauthor{\bsnm{Wu}, \binits{J.}}:
\batitle{The algebraic stability for persistent {L}aplacians}.
\bjtitle{Homology, Homotopy and Applications}
\bvolume{26}(\bissue{2}),
\bfpage{297}--\blpage{323}
(\byear{2024})
\end{barticle}
\endbibitem

\bibitem{liu2024chatgpt}
\begin{barticle}
\bauthor{\bsnm{Liu}, \binits{J.}},
\bauthor{\bsnm{Shen}, \binits{L.}},
\bauthor{\bsnm{Wei}, \binits{G.-W.}}:
\batitle{{ChatGPT} for computational topology}.
\bjtitle{Foundations of data science}
\bvolume{6}(\bissue{2}),
\bfpage{221}--\blpage{250}
(\byear{2024})
\end{barticle}
\endbibitem

\bibitem{liu2024persistentTangles}
\begin{botherref}
\oauthor{\bsnm{Liu}, \binits{J.}},
\oauthor{\bsnm{Shen}, \binits{L.}},
\oauthor{\bsnm{Wei}, \binits{G.-W.}}:
Persistent {K}hovanov homology of tangles.
arXiv preprint arXiv:2409.18312
(2024)
\end{botherref}
\endbibitem

\bibitem{liu2025topological}
\begin{botherref}
\oauthor{\bsnm{Liu}, \binits{J.}},
\oauthor{\bsnm{Shen}, \binits{L.}},
\oauthor{\bsnm{Chen}, \binits{D.}},
\oauthor{\bsnm{Wei}, \binits{G.-W.}}:
Topological sequence analysis of genomes: Delta complex approaches.
arXiv preprint arXiv:2507.05452
(2025)
\end{botherref}
\endbibitem

\bibitem{liu2023persistent}
\begin{barticle}
\bauthor{\bsnm{Liu}, \binits{R.}},
\bauthor{\bsnm{Liu}, \binits{X.}},
\bauthor{\bsnm{Wu}, \binits{J.}}:
\batitle{Persistent path-spectral (pps) based machine learning for
  protein--ligand binding affinity prediction}.
\bjtitle{Journal of chemical information and modeling}
\bvolume{63}(\bissue{3}),
\bfpage{1066}--\blpage{1075}
(\byear{2023})
\end{barticle}
\endbibitem

\bibitem{liu2021neighborhood}
\begin{bchapter}
\bauthor{\bsnm{Liu}, \binits{X.}},
\bauthor{\bsnm{Xia}, \binits{K.}}:
\bctitle{Neighborhood complex based machine learning ({NCML}) models for drug
  design}.
In: \bbtitle{International Workshop on Interpretability of Machine Intelligence
  in Medical Image Computing, and Topological Data Analysis and Its
  Applications for Medical Data},
pp. \bfpage{87}--\blpage{97}
(\byear{2021})
\end{bchapter}
\endbibitem

\bibitem{liu2021hypergraph}
\begin{barticle}
\bauthor{\bsnm{Liu}, \binits{X.}},
\bauthor{\bsnm{Wang}, \binits{X.}},
\bauthor{\bsnm{Wu}, \binits{J.}},
\bauthor{\bsnm{Xia}, \binits{K.}}:
\batitle{Hypergraph-based persistent cohomology ({HPC}) for molecular
  representations in drug design}.
\bjtitle{Briefings in Bioinformatics}
\bvolume{22}(\bissue{5}),
\bfpage{411}
(\byear{2021})
\end{barticle}
\endbibitem

\bibitem{liu2021persistent}
\begin{barticle}
\bauthor{\bsnm{Liu}, \binits{X.}},
\bauthor{\bsnm{Feng}, \binits{H.}},
\bauthor{\bsnm{Wu}, \binits{J.}},
\bauthor{\bsnm{Xia}, \binits{K.}}:
\batitle{Persistent spectral hypergraph based machine learning ({PSH-ML}) for
  protein-ligand binding affinity prediction}.
\bjtitle{Briefings in Bioinformatics}
\bvolume{22}(\bissue{5}),
\bfpage{127}
(\byear{2021})
\end{barticle}
\endbibitem

\bibitem{liu2022dowker}
\begin{barticle}
\bauthor{\bsnm{Liu}, \binits{X.}},
\bauthor{\bsnm{Feng}, \binits{H.}},
\bauthor{\bsnm{Wu}, \binits{J.}},
\bauthor{\bsnm{Xia}, \binits{K.}}:
\batitle{Dowker complex based machine learning ({DCML}) models for
  protein-ligand binding affinity prediction}.
\bjtitle{PLoS computational biology}
\bvolume{18}(\bissue{4}),
\bfpage{1009943}
(\byear{2022})
\end{barticle}
\endbibitem

\bibitem{liu2022hom}
\begin{barticle}
\bauthor{\bsnm{Liu}, \binits{X.}},
\bauthor{\bsnm{Feng}, \binits{H.}},
\bauthor{\bsnm{Wu}, \binits{J.}},
\bauthor{\bsnm{Xia}, \binits{K.}}:
\batitle{Hom-complex-based machine learning ({HCML}) for the prediction of
  protein-protein binding affinity changes upon mutation}.
\bjtitle{Journal of chemical information and modeling}
\bvolume{62}(\bissue{17}),
\bfpage{3961}--\blpage{3969}
(\byear{2022})
\end{barticle}
\endbibitem

\bibitem{liu2024computing}
\begin{barticle}
\bauthor{\bsnm{Liu}, \binits{X.}},
\bauthor{\bsnm{Feng}, \binits{H.}},
\bauthor{\bsnm{Wu}, \binits{J.}},
\bauthor{\bsnm{Xia}, \binits{K.}}:
\batitle{Computing hypergraph homology}.
\bjtitle{Foundations of Data Science}
\bvolume{6}(\bissue{2}),
\bfpage{172}--\blpage{194}
(\byear{2024})
\end{barticle}
\endbibitem

\bibitem{liu2024intcomplex}
\begin{botherref}
\oauthor{\bsnm{Liu}, \binits{X.}},
\oauthor{\bsnm{Liu}, \binits{R.}},
\oauthor{\bsnm{Li}, \binits{J.}},
\oauthor{\bsnm{Wu}, \binits{R.}},
\oauthor{\bsnm{Wu}, \binits{J.}}:
Intcomplex for high-order interactions.
arXiv preprint arXiv:2412.02806
(2024)
\end{botherref}
\endbibitem

\bibitem{liu2025manifold}
\begin{botherref}
\oauthor{\bsnm{Liu}, \binits{X.}},
\oauthor{\bsnm{Su}, \binits{Z.}},
\oauthor{\bsnm{Shi}, \binits{Y.}},
\oauthor{\bsnm{Tong}, \binits{Y.}},
\oauthor{\bsnm{Wang}, \binits{G.}},
\oauthor{\bsnm{Wei}, \binits{G.-W.}}:
Manifold topological deep learning for biomedical data.
arXiv preprint arXiv:2503.00175
(2025)
\end{botherref}
\endbibitem

\bibitem{lloyd2016quantum}
\begin{barticle}
\bauthor{\bsnm{Lloyd}, \binits{S.}},
\bauthor{\bsnm{Garnerone}, \binits{S.}},
\bauthor{\bsnm{Zanardi}, \binits{P.}}:
\batitle{Quantum algorithms for topological and geometric analysis of data}.
\bjtitle{Nature communications}
\bvolume{7}(\bissue{1}),
\bfpage{10138}
(\byear{2016})
\end{barticle}
\endbibitem

\bibitem{love2023topological}
\begin{barticle}
\bauthor{\bsnm{Love}, \binits{E.R.}},
\bauthor{\bsnm{Filippenko}, \binits{B.}},
\bauthor{\bsnm{Maroulas}, \binits{V.}},
\bauthor{\bsnm{Carlsson}, \binits{G.}}:
\batitle{Topological convolutional layers for deep learning}.
\bjtitle{Journal of Machine Learning Research}
\bvolume{24}(\bissue{59}),
\bfpage{1}--\blpage{35}
(\byear{2023})
\end{barticle}
\endbibitem

\bibitem{lutgehetmann2020computing}
\begin{barticle}
\bauthor{\bsnm{L{\"u}tgehetmann}, \binits{D.}},
\bauthor{\bsnm{Govc}, \binits{D.}},
\bauthor{\bsnm{Smith}, \binits{J.P.}},
\bauthor{\bsnm{Levi}, \binits{R.}}:
\batitle{Computing persistent homology of directed flag complexes}.
\bjtitle{Algorithms}
\bvolume{13}(\bissue{1}),
\bfpage{19}
(\byear{2020})
\end{barticle}
\endbibitem

\bibitem{maria2014gudhi}
\begin{bchapter}
\bauthor{\bsnm{Maria}, \binits{C.}},
\bauthor{\bsnm{Boissonnat}, \binits{J.-D.}},
\bauthor{\bsnm{Glisse}, \binits{M.}},
\bauthor{\bsnm{Yvinec}, \binits{M.}}:
\bctitle{The {GUDHI} library: simplicial complexes and persistent homology}.
In: \bbtitle{Mathematical Software--ICMS 2014: 4th International Congress,
  Seoul, South Korea, August 5-9, 2014. Proceedings 4},
pp. \bfpage{167}--\blpage{174}
(\byear{2014}).
\bcomment{Springer}
\end{bchapter}
\endbibitem

\bibitem{maroulas2022bayesian}
\begin{barticle}
\bauthor{\bsnm{Maroulas}, \binits{V.}},
\bauthor{\bsnm{Micucci}, \binits{C.P.}},
\bauthor{\bsnm{Nasrin}, \binits{F.}}:
\batitle{Bayesian topological learning for classifying the structure of
  biological networks}.
\bjtitle{Bayesian Analysis}
\bvolume{17}(\bissue{3}),
\bfpage{711}--\blpage{736}
(\byear{2022}).
\doiurl{10.1214/21-BA1270}
\end{barticle}
\endbibitem

\bibitem{maroulas2019nonparametricPDF}
\begin{barticle}
\bauthor{\bsnm{Maroulas}, \binits{V.}},
\bauthor{\bsnm{Mike}, \binits{J.L.}},
\bauthor{\bsnm{Oballe}, \binits{C.}}:
\batitle{Nonparametric estimation of probability density functions of random
  persistence diagrams}.
\bjtitle{Journal of Machine Learning Research}
\bvolume{20}(\bissue{196}),
\bfpage{1}--\blpage{49}
(\byear{2019})
\end{barticle}
\endbibitem

\bibitem{maroulas2020bayesian}
\begin{barticle}
\bauthor{\bsnm{Maroulas}, \binits{V.}},
\bauthor{\bsnm{Nasrin}, \binits{F.}},
\bauthor{\bsnm{Oballe}, \binits{C.}}:
\batitle{A {Bayesian} framework for persistent homology}.
\bjtitle{SIAM Journal on Mathematics of Data Science}
\bvolume{2}(\bissue{1}),
\bfpage{48}--\blpage{74}
(\byear{2020}).
\doiurl{10.1137/19M1268719}
\end{barticle}
\endbibitem

\bibitem{mayer1942new}
\begin{barticle}
\bauthor{\bsnm{Mayer}, \binits{W.}}:
\batitle{A new homology theory}.
\bjtitle{Annals of Mathematics}
\bvolume{43}(\bissue{2}),
\bfpage{370}--\blpage{380}
(\byear{1942})
\end{barticle}
\endbibitem

\bibitem{mccleary2001user}
\begin{bbook}
\bauthor{\bsnm{McCleary}, \binits{J.}}:
\bbtitle{A User's Guide to Spectral Sequences}
vol. \bseriesno{58}.
\bpublisher{Cambridge University Press},
\blocation{\hspace{0pt}}
(\byear{2001})
\end{bbook}
\endbibitem

\bibitem{memoli2022persistent}
\begin{barticle}
\bauthor{\bsnm{M{\'e}moli}, \binits{F.}},
\bauthor{\bsnm{Wan}, \binits{Z.}},
\bauthor{\bsnm{Wang}, \binits{Y.}}:
\batitle{Persistent {L}aplacians: Properties, algorithms and implications}.
\bjtitle{SIAM Journal on Mathematics of Data Science}
\bvolume{4}(\bissue{2}),
\bfpage{858}--\blpage{884}
(\byear{2022})
\end{barticle}
\endbibitem

\bibitem{mendoza2017parallel}
\begin{botherref}
\oauthor{\bsnm{Mendoza-Smith}, \binits{R.}},
\oauthor{\bsnm{Tanner}, \binits{J.}}:
Parallel multi-scale reduction of persistent homology filtrations.
arXiv preprint arXiv:1708.04710
(2017)
\end{botherref}
\endbibitem

\bibitem{meng2021persistent}
\begin{barticle}
\bauthor{\bsnm{Meng}, \binits{Z.}},
\bauthor{\bsnm{Xia}, \binits{K.}}:
\batitle{Persistent spectral-based machine learning ({PerSpect ML}) for
  protein-ligand binding affinity prediction}.
\bjtitle{Science advances}
\bvolume{7}(\bissue{19}),
\bfpage{5329}
(\byear{2021})
\end{barticle}
\endbibitem

\bibitem{meng2020weighted}
\begin{barticle}
\bauthor{\bsnm{Meng}, \binits{Z.}},
\bauthor{\bsnm{Anand}, \binits{D.V.}},
\bauthor{\bsnm{Lu}, \binits{Y.}},
\bauthor{\bsnm{Wu}, \binits{J.}},
\bauthor{\bsnm{Xia}, \binits{K.}}:
\batitle{Weighted persistent homology for biomolecular data analysis}.
\bjtitle{Scientific reports}
\bvolume{10}(\bissue{1}),
\bfpage{2079}
(\byear{2020})
\end{barticle}
\endbibitem

\bibitem{mileyko2011probability}
\begin{barticle}
\bauthor{\bsnm{Mileyko}, \binits{Y.}},
\bauthor{\bsnm{Mukherjee}, \binits{S.}},
\bauthor{\bsnm{Harer}, \binits{J.}}:
\batitle{Probability measures on the space of persistence diagrams}.
\bjtitle{Inverse Problems}
\bvolume{27}(\bissue{12}),
\bfpage{124007}
(\byear{2011})
\end{barticle}
\endbibitem

\bibitem{miller1984interpretations}
\begin{barticle}
\bauthor{\bsnm{Miller}, \binits{B.P.}}:
\batitle{Interpretations from {H}elmholtz' theorem in classical
  electromagnetism}.
\bjtitle{American Journal of Physics}
\bvolume{52}(\bissue{10}),
\bfpage{948}--\blpage{950}
(\byear{1984})
\end{barticle}
\endbibitem

\bibitem{millett2013identifying}
\begin{barticle}
\bauthor{\bsnm{Millett}, \binits{K.C.}},
\bauthor{\bsnm{Rawdon}, \binits{E.J.}},
\bauthor{\bsnm{Stasiak}, \binits{A.}},
\bauthor{\bsnm{Su{\l}kowska}, \binits{J.I.}}:
\batitle{Identifying knots in proteins}.
\bjtitle{Biochemical Society Transactions}
\bvolume{41}(\bissue{2}),
\bfpage{533}--\blpage{537}
(\byear{2013})
\end{barticle}
\endbibitem

\bibitem{mischaikow2013morse}
\begin{barticle}
\bauthor{\bsnm{Mischaikow}, \binits{K.}},
\bauthor{\bsnm{Nanda}, \binits{V.}}:
\batitle{Morse theory for filtrations and efficient computation of persistent
  homology}.
\bjtitle{Discrete \& Computational Geometry}
\bvolume{50},
\bfpage{330}--\blpage{353}
(\byear{2013})
\end{barticle}
\endbibitem

\bibitem{mischaikow1999conley}
\begin{bbook}
\bauthor{\bsnm{Mischaikow}, \binits{K.}},
\bauthor{\bsnm{Mrozek}, \binits{M.}},
\bauthor{\bsnm{Zgliczy{\'n}ski}, \binits{P.}}:
\bbtitle{Conley Index Theory}
vol. \bseriesno{47}.
\bpublisher{Springer},
\blocation{\hspace{0pt}}
(\byear{1999})
\end{bbook}
\endbibitem

\bibitem{mitchell2024topological}
\begin{barticle}
\bauthor{\bsnm{Mitchell}, \binits{E.C.}},
\bauthor{\bsnm{Story}, \binits{B.}},
\bauthor{\bsnm{Boothe}, \binits{D.}},
\bauthor{\bsnm{Franaszczuk}, \binits{P.J.}},
\bauthor{\bsnm{Maroulas}, \binits{V.}}:
\batitle{A topological deep learning framework for neural spike decoding}.
\bjtitle{Biophysical Journal}
\bvolume{123}(\bissue{17}),
\bfpage{2781}--\blpage{2789}
(\byear{2024}).
\doiurl{10.1016/j.bpj.2024.01.025}
\end{barticle}
\endbibitem

\bibitem{monod2020tropical}
\begin{barticle}
\bauthor{\bsnm{Monod}, \binits{A.}},
\bauthor{\bsnm{Kali{\v{s}}nik}, \binits{S.}},
\bauthor{\bsnm{Pati{\~n}o-Galindo}, \binits{J.{\'A}.}},
\bauthor{\bsnm{Crawford}, \binits{L.}}:
\batitle{Tropical sufficient statistics for persistent homology}.
\bjtitle{SIAM Journal on Applied Algebra and Geometry}
\bvolume{3}(\bissue{2}),
\bfpage{337}--\blpage{371}
(\byear{2019}).
\doiurl{10.1137/17M1148037}
\end{barticle}
\endbibitem

\bibitem{montagna2024topological}
\begin{botherref}
\oauthor{\bsnm{Montagna}, \binits{M.}},
\oauthor{\bsnm{Scardapane}, \binits{S.}},
\oauthor{\bsnm{Telyatnikov}, \binits{L.}}:
Topological deep learning with state-space models: A mamba approach for
  simplicial complexes.
arXiv preprint arXiv:2409.12033
(2024)
\end{botherref}
\endbibitem

\bibitem{morozov2007dionysus}
\begin{botherref}
\oauthor{\bsnm{Morozov}, \binits{D.}}:
Dionysus, a C++ library for computing persistent homology
(2007)
\end{botherref}
\endbibitem

\bibitem{morrey1956variational}
\begin{barticle}
\bauthor{\bsnm{Morrey}, \binits{C.B.}}:
\batitle{A variational method in the theory of harmonic integrals, ii}.
\bjtitle{American Journal of Mathematics}
\bvolume{78}(\bissue{1}),
\bfpage{137}--\blpage{170}
(\byear{1956})
\end{barticle}
\endbibitem

\bibitem{morse1925relations}
\begin{barticle}
\bauthor{\bsnm{Morse}, \binits{M.}}:
\batitle{Relations between the critical points of a real function of n
  independent variables}.
\bjtitle{Transactions of the American Mathematical Society}
\bvolume{27}(\bissue{3}),
\bfpage{345}--\blpage{396}
(\byear{1925})
\end{barticle}
\endbibitem

\bibitem{munch2015probabilistic}
\begin{barticle}
\bauthor{\bsnm{Munch}, \binits{E.}},
\bauthor{\bsnm{Turner}, \binits{K.}},
\bauthor{\bsnm{Bendich}, \binits{P.}},
\bauthor{\bsnm{Mukherjee}, \binits{S.}},
\bauthor{\bsnm{Mattingly}, \binits{J.}},
\bauthor{\bsnm{Harer}, \binits{J.}}:
\batitle{Probabilistic {Fr\'e}chet means for time varying persistence
  diagrams}.
\bjtitle{Electronic Journal of Statistics}
\bvolume{9},
\bfpage{1173}--\blpage{1204}
(\byear{2015})
\end{barticle}
\endbibitem

\bibitem{nguyen2019mathematical}
\begin{barticle}
\bauthor{\bsnm{Nguyen}, \binits{D.D.}},
\bauthor{\bsnm{Cang}, \binits{Z.}},
\bauthor{\bsnm{Wu}, \binits{K.}},
\bauthor{\bsnm{Wang}, \binits{M.}},
\bauthor{\bsnm{Cao}, \binits{Y.}},
\bauthor{\bsnm{Wei}, \binits{G.-W.}}:
\batitle{Mathematical deep learning for pose and binding affinity prediction
  and ranking in {D3R Grand Challenges}}.
\bjtitle{Journal of computer-aided molecular design}
\bvolume{33},
\bfpage{71}--\blpage{82}
(\byear{2019})
\end{barticle}
\endbibitem

\bibitem{nguyen2020mathdl}
\begin{barticle}
\bauthor{\bsnm{Nguyen}, \binits{D.D.}},
\bauthor{\bsnm{Gao}, \binits{K.}},
\bauthor{\bsnm{Wang}, \binits{M.}},
\bauthor{\bsnm{Wei}, \binits{G.-W.}}:
\batitle{Mathdl: mathematical deep learning for {D3R} {Grand Challenge} 4}.
\bjtitle{Journal of computer-aided molecular design}
\bvolume{34},
\bfpage{131}--\blpage{147}
(\byear{2020})
\end{barticle}
\endbibitem

\bibitem{nielsen2010quantum}
\begin{bbook}
\bauthor{\bsnm{Nielsen}, \binits{M.A.}},
\bauthor{\bsnm{Chuang}, \binits{I.L.}}:
\bbtitle{Quantum Computation and Quantum Information}.
\bpublisher{Cambridge university press},
\blocation{\hspace{0pt}}
(\byear{2010})
\end{bbook}
\endbibitem

\bibitem{oballe2022bayesian}
\begin{barticle}
\bauthor{\bsnm{Oballe}, \binits{C.}},
\bauthor{\bsnm{Cherne}, \binits{A.}},
\bauthor{\bsnm{Boothe}, \binits{D.}},
\bauthor{\bsnm{Kerick}, \binits{S.}},
\bauthor{\bsnm{Franaszczuk}, \binits{P.J.}},
\bauthor{\bsnm{Maroulas}, \binits{V.}}:
\batitle{Bayesian topological signal processing}.
\bjtitle{Discrete and Continuous Dynamical Systems - Series S}
\bvolume{15}(\bissue{4}),
\bfpage{797}--\blpage{817}
(\byear{2022}).
\doiurl{10.3934/dcdss.2021084}
\end{barticle}
\endbibitem

\bibitem{ohtsuki2001quantum}
\begin{bbook}
\bauthor{\bsnm{Ohtsuki}, \binits{T.}}:
\bbtitle{Quantum Invariants: a Study of Knots, 3-manifolds, and Their Sets}
vol. \bseriesno{29}.
\bpublisher{World Scientific},
\blocation{\hspace{0pt}}
(\byear{2001})
\end{bbook}
\endbibitem

\bibitem{ozsvath2004holomorphic}
\begin{barticle}
\bauthor{\bsnm{Ozsv{\'a}th}, \binits{P.}},
\bauthor{\bsnm{Szab{\'o}}, \binits{Z.}}:
\batitle{Holomorphic disks and knot invariants}.
\bjtitle{Advances in Mathematics}
\bvolume{186}(\bissue{1}),
\bfpage{58}--\blpage{116}
(\byear{2004})
\end{barticle}
\endbibitem

\bibitem{panagiotou2020knot}
\begin{barticle}
\bauthor{\bsnm{Panagiotou}, \binits{E.}},
\bauthor{\bsnm{Kauffman}, \binits{L.H.}}:
\batitle{Knot polynomials of open and closed curves}.
\bjtitle{Proceedings of the Royal Society A}
\bvolume{476}(\bissue{2240}),
\bfpage{20200124}
(\byear{2020})
\end{barticle}
\endbibitem

\bibitem{panagiotou2020topological}
\begin{barticle}
\bauthor{\bsnm{Panagiotou}, \binits{E.}},
\bauthor{\bsnm{Plaxco}, \binits{K.W.}}:
\batitle{A topological study of protein folding kinetics}.
\bjtitle{Topol. Geom. Biopolym. AMS Contemp. Math. Ser}
\bvolume{746},
\bfpage{223}--\blpage{233}
(\byear{2020})
\end{barticle}
\endbibitem

\bibitem{panagiotou2019topological}
\begin{barticle}
\bauthor{\bsnm{Panagiotou}, \binits{E.}},
\bauthor{\bsnm{Millett}, \binits{K.C.}},
\bauthor{\bsnm{Atzberger}, \binits{P.J.}}:
\batitle{Topological methods for polymeric materials: characterizing the
  relationship between polymer entanglement and viscoelasticity}.
\bjtitle{Polymers}
\bvolume{11}(\bissue{3}),
\bfpage{437}
(\byear{2019})
\end{barticle}
\endbibitem

\bibitem{papamarkou2022random}
\begin{botherref}
\oauthor{\bsnm{Papamarkou}, \binits{T.}},
\oauthor{\bsnm{Nasrin}, \binits{F.}},
\oauthor{\bsnm{Lawson}, \binits{A.}},
\oauthor{\bsnm{Gong}, \binits{N.}},
\oauthor{\bsnm{Rios}, \binits{O.}},
\oauthor{\bsnm{Maroulas}, \binits{V.}}:
A random persistence diagram generator.
Statistics and Computing
\textbf{32}(88)
(2022).
\doiurl{10.1007/s11222-022-10141-y}
\end{botherref}
\endbibitem

\bibitem{papamarkou2024position}
\begin{barticle}
\bauthor{\bsnm{Papamarkou}, \binits{T.}},
\bauthor{\bsnm{Birdal}, \binits{T.}},
\bauthor{\bsnm{Bronstein}, \binits{M.M.}},
\bauthor{\bsnm{Carlsson}, \binits{G.}},
\bauthor{\bsnm{Curry}, \binits{J.}},
\bauthor{\bsnm{Gao}, \binits{Y.}},
\bauthor{\bsnm{Hajij}, \binits{M.}},
\bauthor{\bsnm{Kwitt}, \binits{R.}},
\bauthor{\bsnm{Li{\`o}}, \binits{P.}},
\bauthor{\bsnm{Di~Lorenzo}, \binits{P.}},
\bauthor{\bsnm{Maroulas}, \binits{V.}},
\bauthor{\bsnm{Miolane}, \binits{N.}},
\bauthor{\bsnm{Nasrin}, \binits{F.}},
\bauthor{\bsnm{Natesan~Ramamurthy}, \binits{K.}},
\bauthor{\bsnm{Rieck}, \binits{B.}},
\bauthor{\bsnm{Scardapane}, \binits{S.}},
\bauthor{\bsnm{Schaub}, \binits{M.T.}},
\bauthor{\bsnm{Veli{\v{c}}kovi{\'c}}, \binits{P.}},
\bauthor{\bsnm{Wang}, \binits{B.}},
\bauthor{\bsnm{Wang}, \binits{Y.}},
\bauthor{\bsnm{Wei}, \binits{G.-W.}},
\bauthor{\bsnm{Zamzmi}, \binits{G.}}:
\batitle{Position: Topological deep learning is the new frontier for relational
  learning}.
\bjtitle{Proceedings of Machine Learning Research}
\bvolume{235},
\bfpage{39529}--\blpage{39555}
(\byear{2024}).
\bcomment{41st International Conference on Machine Learning (ICML 2024)}
\end{barticle}
\endbibitem

\bibitem{patel2018generalized}
\begin{barticle}
\bauthor{\bsnm{Patel}, \binits{A.}}:
\batitle{Generalized persistence diagrams}.
\bjtitle{Journal of Applied and Computational Topology}
\bvolume{1}(\bissue{3}),
\bfpage{397}--\blpage{419}
(\byear{2018})
\end{barticle}
\endbibitem

\bibitem{perea2016persistent}
\begin{bchapter}
\bauthor{\bsnm{Perea}, \binits{J.A.}}:
\bctitle{Persistent homology of toroidal sliding window embeddings}.
In: \bbtitle{2016 IEEE International Conference on Acoustics, Speech and Signal
  Processing (icassp)},
pp. \bfpage{6435}--\blpage{6439}
(\byear{2016}).
\bcomment{IEEE}
\end{bchapter}
\endbibitem

\bibitem{perea2015sliding}
\begin{barticle}
\bauthor{\bsnm{Perea}, \binits{J.A.}},
\bauthor{\bsnm{Harer}, \binits{J.}}:
\batitle{Sliding windows and persistence: an application of topological methods
  to signal analysis}.
\bjtitle{Foundations of computational mathematics}
\bvolume{15},
\bfpage{799}--\blpage{838}
(\byear{2015})
\end{barticle}
\endbibitem

\bibitem{petri2013topological}
\begin{barticle}
\bauthor{\bsnm{Petri}, \binits{G.}},
\bauthor{\bsnm{Scolamiero}, \binits{M.}},
\bauthor{\bsnm{Donato}, \binits{I.}},
\bauthor{\bsnm{Vaccarino}, \binits{F.}}:
\batitle{Topological strata of weighted complex networks}.
\bjtitle{PloS one}
\bvolume{8}(\bissue{6}),
\bfpage{66506}
(\byear{2013})
\end{barticle}
\endbibitem

\bibitem{poelke2017hodge}
\begin{botherref}
\oauthor{\bsnm{Poelke}, \binits{K.}}:
Hodge-type decompositions for piecewise constant vector fields on simplicial
  surfaces and solids with boundary.
PhD thesis
(2017)
\end{botherref}
\endbibitem

\bibitem{poelke2016boundary}
\begin{barticle}
\bauthor{\bsnm{Poelke}, \binits{K.}},
\bauthor{\bsnm{Polthier}, \binits{K.}}:
\batitle{Boundary-aware {Hodge} decompositions for piecewise constant vector
  fields}.
\bjtitle{Computer-Aided Design}
\bvolume{78},
\bfpage{126}--\blpage{136}
(\byear{2016})
\end{barticle}
\endbibitem

\bibitem{pun2022persistent}
\begin{barticle}
\bauthor{\bsnm{Pun}, \binits{C.S.}},
\bauthor{\bsnm{Lee}, \binits{S.X.}},
\bauthor{\bsnm{Xia}, \binits{K.}}:
\batitle{Persistent-homology-based machine learning: a survey and a comparative
  study}.
\bjtitle{Artificial Intelligence Review}
\bvolume{55}(\bissue{7}),
\bfpage{5169}--\blpage{5213}
(\byear{2022})
\end{barticle}
\endbibitem

\bibitem{pun2020weighted}
\begin{barticle}
\bauthor{\bsnm{Pun}, \binits{C.S.}},
\bauthor{\bsnm{Yong}, \binits{B.Y.S.}},
\bauthor{\bsnm{Xia}, \binits{K.}}:
\batitle{Weighted-persistent-homology-based machine learning for {RNA}
  flexibility analysis}.
\bjtitle{PloS one}
\bvolume{15}(\bissue{8}),
\bfpage{0237747}
(\byear{2020})
\end{barticle}
\endbibitem

\bibitem{puzyn2010recent}
\begin{botherref}
\oauthor{\bsnm{Puzyn}, \binits{T.}},
\oauthor{\bsnm{Leszczynski}, \binits{J.}},
\oauthor{\bsnm{Cronin}, \binits{M.T.}}:
Recent advances in {QSAR} studies: methods and applications
(2010)
\end{botherref}
\endbibitem

\bibitem{qiu2023persistent}
\begin{barticle}
\bauthor{\bsnm{Qiu}, \binits{Y.}},
\bauthor{\bsnm{Wei}, \binits{G.-W.}}:
\batitle{Persistent spectral theory-guided protein engineering}.
\bjtitle{Nature computational science}
\bvolume{3}(\bissue{2}),
\bfpage{149}--\blpage{163}
(\byear{2023})
\end{barticle}
\endbibitem

\bibitem{reeb1946points}
\begin{barticle}
\bauthor{\bsnm{Reeb}, \binits{G.}}:
\batitle{Sur les points singuliers d'une forme de pfaff completement integrable
  ou d'une fonction numerique [on the singular points of a completely
  integrable pfaff form or of a numerical function]}.
\bjtitle{Comptes Rendus Acad. Sciences Paris}
\bvolume{222},
\bfpage{847}--\blpage{849}
(\byear{1946})
\end{barticle}
\endbibitem

\bibitem{reininghaus2015stable}
\begin{bchapter}
\bauthor{\bsnm{Reininghaus}, \binits{J.}},
\bauthor{\bsnm{Huber}, \binits{S.}},
\bauthor{\bsnm{Bauer}, \binits{U.}},
\bauthor{\bsnm{Kwitt}, \binits{R.}}:
\bctitle{A stable multi-scale kernel for topological machine learning}.
In: \bbtitle{Proceedings of the IEEE Conference on Computer Vision and Pattern
  Recognition},
pp. \bfpage{4741}--\blpage{4748}
(\byear{2015})
\end{bchapter}
\endbibitem

\bibitem{ren2020stability}
\begin{botherref}
\oauthor{\bsnm{Ren}, \binits{S.}},
\oauthor{\bsnm{Wu}, \binits{J.}}:
The stability of persistent homology of hypergraphs.
arXiv preprint arXiv:2002.02237
(2020)
\end{botherref}
\endbibitem

\bibitem{ren2018weighted}
\begin{barticle}
\bauthor{\bsnm{Ren}, \binits{S.}},
\bauthor{\bsnm{Wu}, \binits{C.}},
\bauthor{\bsnm{Wu}, \binits{J.}}:
\batitle{Weighted persistent homology}.
\bjtitle{The Rocky Mountain Journal of Mathematics}
\bvolume{48}(\bissue{8}),
\bfpage{2661}--\blpage{2687}
(\byear{2018})
\end{barticle}
\endbibitem

\bibitem{ren2021discrete}
\begin{botherref}
\oauthor{\bsnm{Ren}, \binits{S.}},
\oauthor{\bsnm{Wang}, \binits{C.}},
\oauthor{\bsnm{Wu}, \binits{C.}},
\oauthor{\bsnm{Wu}, \binits{J.}}:
On the discrete {Morse} functions for hypergraphs.
arXiv preprint arXiv:2108.02384
(2021)
\end{botherref}
\endbibitem

\bibitem{reuter2006laplace}
\begin{barticle}
\bauthor{\bsnm{Reuter}, \binits{M.}},
\bauthor{\bsnm{Wolter}, \binits{F.-E.}},
\bauthor{\bsnm{Peinecke}, \binits{N.}}:
\batitle{Laplace-{B}eltrami spectra as `{Shape-DNA}' of surfaces and solids}.
\bjtitle{Computer-Aided Design}
\bvolume{38}(\bissue{4}),
\bfpage{342}--\blpage{366}
(\byear{2006})
\end{barticle}
\endbibitem

\bibitem{ribando2024graph}
\begin{barticle}
\bauthor{\bsnm{Ribando-Gros}, \binits{E.}},
\bauthor{\bsnm{Wang}, \binits{R.}},
\bauthor{\bsnm{Chen}, \binits{J.}},
\bauthor{\bsnm{Tong}, \binits{Y.}},
\bauthor{\bsnm{Wei}, \binits{G.-W.}}:
\batitle{Combinatorial and {H}odge {L}aplacians: Similarity and difference}.
\bjtitle{SIAM Review}
\bvolume{66}(\bissue{3}),
\bfpage{575}--\blpage{601}
(\byear{2024})
\end{barticle}
\endbibitem

\bibitem{russold2022persistent}
\begin{botherref}
\oauthor{\bsnm{Russold}, \binits{F.}}:
Persistent sheaf cohomology.
arXiv preprint arXiv:2204.13446
(2022)
\end{botherref}
\endbibitem

\bibitem{rustamov2007laplace}
\begin{bchapter}
\bauthor{\bsnm{Rustamov}, \binits{R.M.}}, \betal:
\bctitle{Laplace-{B}eltrami eigenfunctions for deformation invariant shape
  representation}.
In: \bbtitle{Symposium on Geometry Processing},
vol. \bseriesno{257},
pp. \bfpage{225}--\blpage{233}
(\byear{2007})
\end{bchapter}
\endbibitem

\bibitem{schenck2022algebraic}
\begin{bbook}
\bauthor{\bsnm{Schenck}, \binits{H.}}:
\bbtitle{Algebraic Foundations for Applied Topology and Data Analysis}.
\bpublisher{Springer},
\blocation{\hspace{0pt}}
(\byear{2022})
\end{bbook}
\endbibitem

\bibitem{schlick2021knot}
\begin{barticle}
\bauthor{\bsnm{Schlick}, \binits{T.}},
\bauthor{\bsnm{Zhu}, \binits{Q.}},
\bauthor{\bsnm{Dey}, \binits{A.}},
\bauthor{\bsnm{Jain}, \binits{S.}},
\bauthor{\bsnm{Yan}, \binits{S.}},
\bauthor{\bsnm{Laederach}, \binits{A.}}:
\batitle{To knot or not to knot: multiple conformations of the {SARS-CoV-2}
  frameshifting {RNA} element}.
\bjtitle{Journal of the American Chemical Society}
\bvolume{143}(\bissue{30}),
\bfpage{11404}--\blpage{11422}
(\byear{2021})
\end{barticle}
\endbibitem

\bibitem{schubert1954numerische}
\begin{barticle}
\bauthor{\bsnm{Schubert}, \binits{H.}}:
\batitle{{\"U}ber eine numerische knoteninvariante}.
\bjtitle{Mathematische Zeitschrift}
\bvolume{61}(\bissue{1}),
\bfpage{245}--\blpage{288}
(\byear{1954})
\end{barticle}
\endbibitem

\bibitem{schwarz2006hodge}
\begin{bbook}
\bauthor{\bsnm{Schwarz}, \binits{G.}}:
\bbtitle{Hodge Decomposition - A Method for Solving Boundary Value Problems}.
\bpublisher{Springer},
\blocation{\hspace{0pt}}
(\byear{2006})
\end{bbook}
\endbibitem

\bibitem{segal1968classifying}
\begin{barticle}
\bauthor{\bsnm{Segal}, \binits{G.}}:
\batitle{Classifying spaces and spectral sequences}.
\bjtitle{Publications Math{\'e}matiques de l'IH{\'E}S}
\bvolume{34},
\bfpage{105}--\blpage{112}
(\byear{1968})
\end{barticle}
\endbibitem

\bibitem{shen2024evolutionary}
\begin{barticle}
\bauthor{\bsnm{Shen}, \binits{L.}},
\bauthor{\bsnm{Liu}, \binits{J.}},
\bauthor{\bsnm{Wei}, \binits{G.-W.}}:
\batitle{Evolutionary {K}hovanov homology}.
\bjtitle{AIMS Mathematics}
\bvolume{9}(\bissue{9}),
\bfpage{26139}--\blpage{26165}
(\byear{2024})
\end{barticle}
\endbibitem

\bibitem{shen2024persistent}
\begin{barticle}
\bauthor{\bsnm{Shen}, \binits{L.}},
\bauthor{\bsnm{Liu}, \binits{J.}},
\bauthor{\bsnm{Wei}, \binits{G.-W.}}:
\batitle{Persistent {M}ayer homology and persistent {M}ayer {L}aplacian}.
\bjtitle{Foundations of Data Science}
\bvolume{6}(\bissue{4}),
\bfpage{584}--\blpage{612}
(\byear{2024})
\end{barticle}
\endbibitem

\bibitem{shen2024knot}
\begin{barticle}
\bauthor{\bsnm{Shen}, \binits{L.}},
\bauthor{\bsnm{Feng}, \binits{H.}},
\bauthor{\bsnm{Li}, \binits{F.}},
\bauthor{\bsnm{Lei}, \binits{F.}},
\bauthor{\bsnm{Wu}, \binits{J.}},
\bauthor{\bsnm{Wei}, \binits{G.-W.}}:
\batitle{Knot data analysis using multiscale {Gauss} link integral}.
\bjtitle{Proceedings of the National Academy of Sciences}
\bvolume{121}(\bissue{42}),
\bfpage{2408431121}
(\byear{2024})
\end{barticle}
\endbibitem

\bibitem{shepard1985cellular}
\begin{bbook}
\bauthor{\bsnm{Shepard}, \binits{A.D.}}:
\bbtitle{A Cellular Description of the Derived Category of a Stratified Space}.
\bpublisher{Brown University},
\blocation{\hspace{0pt}}
(\byear{1985})
\end{bbook}
\endbibitem

\bibitem{shonkwiler2009poincare}
\begin{botherref}
\oauthor{\bsnm{Shonkwiler}, \binits{C.}}:
Poincar{\'e} duality angles on {R}iemannian manifolds with boundary.
PhD thesis,
University of Pennsylvania
(2009)
\end{botherref}
\endbibitem

\bibitem{silver2019knot}
\begin{barticle}
\bauthor{\bsnm{Silver}, \binits{D.S.}},
\bauthor{\bsnm{Williams}, \binits{S.G.}}:
\batitle{Knot invariants from {L}aplacian matrices}.
\bjtitle{Journal of Knot Theory and Its Ramifications}
\bvolume{28}(\bissue{09}),
\bfpage{1950058}
(\byear{2019})
\end{barticle}
\endbibitem

\bibitem{singh2007topological}
\begin{barticle}
\bauthor{\bsnm{Singh}, \binits{G.}},
\bauthor{\bsnm{M{\'e}moli}, \binits{F.}},
\bauthor{\bsnm{Carlsson}, \binits{G.E.}}, \betal:
\batitle{Topological methods for the analysis of high dimensional data sets and
  {3D} object recognition}.
\bjtitle{PBG@ Eurographics}
\bvolume{2},
\bfpage{091}--\blpage{100}
(\byear{2007})
\end{barticle}
\endbibitem

\bibitem{singh2023topological}
\begin{barticle}
\bauthor{\bsnm{Singh}, \binits{Y.}},
\bauthor{\bsnm{Farrelly}, \binits{C.M.}},
\bauthor{\bsnm{Hathaway}, \binits{Q.A.}},
\bauthor{\bsnm{Leiner}, \binits{T.}},
\bauthor{\bsnm{Jagtap}, \binits{J.}},
\bauthor{\bsnm{Carlsson}, \binits{G.E.}},
\bauthor{\bsnm{Erickson}, \binits{B.J.}}:
\batitle{Topological data analysis in medical imaging: current state of the
  art}.
\bjtitle{Insights into Imaging}
\bvolume{14}(\bissue{1}),
\bfpage{58}
(\byear{2023})
\end{barticle}
\endbibitem

\bibitem{som2020pi}
\begin{bchapter}
\bauthor{\bsnm{Som}, \binits{A.}},
\bauthor{\bsnm{Choi}, \binits{H.}},
\bauthor{\bsnm{Ramamurthy}, \binits{K.N.}},
\bauthor{\bsnm{Buman}, \binits{M.P.}},
\bauthor{\bsnm{Turaga}, \binits{P.}}:
\bctitle{Pi-net: A deep learning approach to extract topological persistence
  images}.
In: \bbtitle{Proceedings of the IEEE/CVF Conference on Computer Vision and
  Pattern Recognition Workshops},
pp. \bfpage{834}--\blpage{835}
(\byear{2020})
\end{bchapter}
\endbibitem

\bibitem{song2025multiscale}
\begin{barticle}
\bauthor{\bsnm{Song}, \binits{R.}},
\bauthor{\bsnm{Li}, \binits{F.}},
\bauthor{\bsnm{Wu}, \binits{J.}},
\bauthor{\bsnm{Lei}, \binits{F.}},
\bauthor{\bsnm{Wei}, \binits{G.-W.}}:
\batitle{Multi-scale jones polynomial and persistent jones polynomial for knot
  data analysis}.
\bjtitle{AIMS Mathematics}
\bvolume{10}(\bissue{1}),
\bfpage{1463}--\blpage{1487}
(\byear{2025})
\end{barticle}
\endbibitem

\bibitem{spanier1949mayer}
\begin{barticle}
\bauthor{\bsnm{Spanier}, \binits{E.H.}}:
\batitle{The {Mayer} homology theory}.
\bjtitle{Bull. Amer. Math. Soc.}
\bvolume{55}(\bissue{12}),
\bfpage{102}--\blpage{112}
(\byear{1949})
\end{barticle}
\endbibitem

\bibitem{stolz2017persistent}
\begin{botherref}
\oauthor{\bsnm{Stolz}, \binits{B.J.}},
\oauthor{\bsnm{Harrington}, \binits{H.A.}},
\oauthor{\bsnm{Porter}, \binits{M.A.}}:
Persistent homology of time-dependent functional networks constructed from
  coupled time series.
Chaos: An Interdisciplinary Journal of Nonlinear Science
\textbf{27}(4)
(2017)
\end{botherref}
\endbibitem

\bibitem{su2024hodgerna}
\begin{barticle}
\bauthor{\bsnm{Su}, \binits{Z.}},
\bauthor{\bsnm{Tong}, \binits{Y.}},
\bauthor{\bsnm{Wei}, \binits{G.-W.}}:
\batitle{Hodge decomposition of single-cell {RNA} velocity}.
\bjtitle{Journal of chemical information and modeling}
\bvolume{64}(\bissue{8}),
\bfpage{3558}--\blpage{3568}
(\byear{2024})
\end{barticle}
\endbibitem

\bibitem{su2024hodge}
\begin{bchapter}
\bauthor{\bsnm{Su}, \binits{Z.}},
\bauthor{\bsnm{Tong}, \binits{Y.}},
\bauthor{\bsnm{Wei}, \binits{G.-W.}}:
\bctitle{Hodge decomposition of vector fields in {Cartesian} grids}.
In: \bbtitle{SIGGRAPH Asia 2024 Conference Papers},
pp. \bfpage{1}--\blpage{10}
(\byear{2024})
\end{bchapter}
\endbibitem

\bibitem{su2024persistent}
\begin{barticle}
\bauthor{\bsnm{Su}, \binits{Z.}},
\bauthor{\bsnm{Tong}, \binits{Y.}},
\bauthor{\bsnm{Wei}, \binits{G.-W.}}:
\batitle{Persistent de {R}ham-{H}odge {L}aplacians in {E}ulerian representation
  for manifold topological learning}.
\bjtitle{AIMS Mathematics}
\bvolume{9}(\bissue{10}),
\bfpage{27438}--\blpage{27470}
(\byear{2024})
\end{barticle}
\endbibitem

\bibitem{su2024topology}
\begin{botherref}
\oauthor{\bsnm{Su}, \binits{Z.}},
\oauthor{\bsnm{Tong}, \binits{Y.}},
\oauthor{\bsnm{Wei}, \binits{G.-W.}}:
Topology-preserving {Hodge} decomposition in the {Eulerian} representation.
arXiv preprint arXiv:2408.14356
(2024)
\end{botherref}
\endbibitem

\bibitem{sulkowska2012conservation}
\begin{barticle}
\bauthor{\bsnm{Sulkowska}, \binits{J.I.}},
\bauthor{\bsnm{Rawdon}, \binits{E.J.}},
\bauthor{\bsnm{Millet}, \binits{K.C.}},
\bauthor{\bsnm{Onuchic}, \binits{J.N.}},
\bauthor{\bsnm{Stasiak}, \binits{A.}}:
\batitle{Conservation of complex knotting and slipknotting patterns in
  proteins}.
\bjtitle{Biophysical Journal}
\bvolume{102}(\bissue{3}),
\bfpage{253}
(\byear{2012})
\end{barticle}
\endbibitem

\bibitem{sumners2020role}
\begin{bchapter}
\bauthor{\bsnm{Sumners}, \binits{D.}}:
\bctitle{The role of knot theory in {DNA} research}.
In: \bbtitle{Geometry and Topology},
pp. \bfpage{297}--\blpage{318}.
\bpublisher{CRC Press},
\blocation{\hspace{0pt}}
(\byear{2020})
\end{bchapter}
\endbibitem

\bibitem{suwayyid2024persistentB}
\begin{barticle}
\bauthor{\bsnm{Suwayyid}, \binits{F.}},
\bauthor{\bsnm{Wei}, \binits{G.-W.}}:
\batitle{Persistent {Dirac} of paths on digraphs and hypergraphs}.
\bjtitle{Foundations of data science (Springfield, Mo.)}
\bvolume{6}(\bissue{2}),
\bfpage{124}
(\byear{2024})
\end{barticle}
\endbibitem

\bibitem{suwayyid2024persistent}
\begin{barticle}
\bauthor{\bsnm{Suwayyid}, \binits{F.}},
\bauthor{\bsnm{Wei}, \binits{G.-W.}}:
\batitle{Persistent {M}ayer {Dirac}}.
\bjtitle{Journal of Physics: Complexity}
\bvolume{5}(\bissue{4}),
\bfpage{045005}
(\byear{2024})
\end{barticle}
\endbibitem

\bibitem{takens2006detecting}
\begin{bchapter}
\bauthor{\bsnm{Takens}, \binits{F.}}:
\bctitle{Detecting strange attractors in turbulence}.
In: \bbtitle{Dynamical Systems and Turbulence, Warwick 1980: Proceedings of a
  Symposium Held at the University of Warwick 1979/80},
pp. \bfpage{366}--\blpage{381}
(\byear{2006}).
\bcomment{Springer}
\end{bchapter}
\endbibitem

\bibitem{tauzin2021giotto}
\begin{barticle}
\bauthor{\bsnm{Tauzin}, \binits{G.}},
\bauthor{\bsnm{Lupo}, \binits{U.}},
\bauthor{\bsnm{Tunstall}, \binits{L.}},
\bauthor{\bsnm{P{\'e}rez}, \binits{J.B.}},
\bauthor{\bsnm{Caorsi}, \binits{M.}},
\bauthor{\bsnm{Medina-Mardones}, \binits{A.M.}},
\bauthor{\bsnm{Dassatti}, \binits{A.}},
\bauthor{\bsnm{Hess}, \binits{K.}}:
\batitle{giotto-tda:: A topological data analysis toolkit for machine learning
  and data exploration}.
\bjtitle{Journal of Machine Learning Research}
\bvolume{22}(\bissue{39}),
\bfpage{1}--\blpage{6}
(\byear{2021})
\end{barticle}
\endbibitem

\bibitem{tierny2017topology}
\begin{barticle}
\bauthor{\bsnm{Tierny}, \binits{J.}},
\bauthor{\bsnm{Favelier}, \binits{G.}},
\bauthor{\bsnm{Levine}, \binits{J.A.}},
\bauthor{\bsnm{Gueunet}, \binits{C.}},
\bauthor{\bsnm{Michaux}, \binits{M.}}:
\batitle{The topology toolkit}.
\bjtitle{IEEE transactions on visualization and computer graphics}
\bvolume{24}(\bissue{1}),
\bfpage{832}--\blpage{842}
(\byear{2017})
\end{barticle}
\endbibitem

\bibitem{townsend2020representation}
\begin{barticle}
\bauthor{\bsnm{Townsend}, \binits{J.}},
\bauthor{\bsnm{Micucci}, \binits{C.P.}},
\bauthor{\bsnm{Hymel}, \binits{J.H.}},
\bauthor{\bsnm{Maroulas}, \binits{V.}},
\bauthor{\bsnm{Vogiatzis}, \binits{K.D.}}:
\batitle{Representation of molecular structures with persistent homology for
  machine learning applications in chemistry}.
\bjtitle{Nature communications}
\bvolume{11}(\bissue{1}),
\bfpage{3230}
(\byear{2020})
\end{barticle}
\endbibitem

\bibitem{turner2024extended}
\begin{barticle}
\bauthor{\bsnm{Turner}, \binits{K.}},
\bauthor{\bsnm{Robins}, \binits{V.}},
\bauthor{\bsnm{Morgan}, \binits{J.}}:
\batitle{The extended persistent homology transform of manifolds with
  boundary}.
\bjtitle{Journal of Applied and Computational Topology}
\bvolume{8}(\bissue{7}),
\bfpage{2111}--\blpage{2154}
(\byear{2024})
\end{barticle}
\endbibitem

\bibitem{turner2014frechet}
\begin{barticle}
\bauthor{\bsnm{Turner}, \binits{K.}},
\bauthor{\bsnm{Mileyko}, \binits{Y.}},
\bauthor{\bsnm{Mukherjee}, \binits{S.}},
\bauthor{\bsnm{Harer}, \binits{J.}}:
\batitle{Fréchet means for distributions of persistence diagrams}.
\bjtitle{Discrete \& Computational Geometry}
\bvolume{52}(\bissue{1}),
\bfpage{44}--\blpage{70}
(\byear{2014})
\end{barticle}
\endbibitem

\bibitem{van2019kepler}
\begin{barticle}
\bauthor{\bsnm{Van~Veen}, \binits{H.J.}},
\bauthor{\bsnm{Saul}, \binits{N.}},
\bauthor{\bsnm{Eargle}, \binits{D.}},
\bauthor{\bsnm{Mangham}, \binits{S.W.}}:
\batitle{Kepler {Mapper}: a flexible {Python} implementation of the {Mapper}
  algorithm.}
\bjtitle{Journal of Open Source Software}
\bvolume{4}(\bissue{42}),
\bfpage{1315}
(\byear{2019})
\end{barticle}
\endbibitem

\bibitem{vietoris1927hoheren}
\begin{barticle}
\bauthor{\bsnm{Vietoris}, \binits{L.}}:
\batitle{{\"U}ber den h{\"o}heren zusammenhang kompakter r{\"a}ume und eine
  klasse von zusammenhangstreuen abbildungen}.
\bjtitle{Mathematische Annalen}
\bvolume{97}(\bissue{1}),
\bfpage{454}--\blpage{472}
(\byear{1927})
\end{barticle}
\endbibitem

\bibitem{vipond2020multiparameter}
\begin{barticle}
\bauthor{\bsnm{Vipond}, \binits{O.}}:
\batitle{Multiparameter persistence landscapes}.
\bjtitle{Journal of Machine Learning Research}
\bvolume{21}(\bissue{61}),
\bfpage{1}--\blpage{38}
(\byear{2020})
\end{barticle}
\endbibitem

\bibitem{wagner2023slice}
\begin{bchapter}
\bauthor{\bsnm{Wagner}, \binits{H.}}:
\bctitle{Slice, simplify and stitch: Topology-preserving simplification scheme
  for massive voxel data}.
In: \bbtitle{39th International Symposium on Computational Geometry (SoCG
  2023)}
(\byear{2023}).
\bcomment{Schloss Dagstuhl-Leibniz-Zentrum f{\"u}r Informatik}
\end{bchapter}
\endbibitem

\bibitem{wagner2011efficient}
\begin{bchapter}
\bauthor{\bsnm{Wagner}, \binits{H.}},
\bauthor{\bsnm{Chen}, \binits{C.}},
\bauthor{\bsnm{Vu{\c{c}}ini}, \binits{E.}}:
\bctitle{Efficient computation of persistent homology for cubical data}.
In: \bbtitle{Topological Methods in Data Analysis and Visualization II: Theory,
  Algorithms, and Applications},
pp. \bfpage{91}--\blpage{106}.
\bpublisher{Springer},
\blocation{\hspace{0pt}}
(\byear{2011})
\end{bchapter}
\endbibitem

\bibitem{wang2016object}
\begin{barticle}
\bauthor{\bsnm{Wang}, \binits{B.}},
\bauthor{\bsnm{Wei}, \binits{G.-W.}}:
\batitle{Object-oriented persistent homology}.
\bjtitle{Journal of computational physics}
\bvolume{305},
\bfpage{276}--\blpage{299}
(\byear{2016})
\end{barticle}
\endbibitem

\bibitem{wang2023persistent}
\begin{barticle}
\bauthor{\bsnm{Wang}, \binits{R.}},
\bauthor{\bsnm{Wei}, \binits{G.-W.}}:
\batitle{Persistent path {L}aplacian}.
\bjtitle{Foundations of data science (Springfield, Mo.)}
\bvolume{5}(\bissue{1}),
\bfpage{26}
(\byear{2023})
\end{barticle}
\endbibitem

\bibitem{wang2021mechanisms}
\begin{barticle}
\bauthor{\bsnm{Wang}, \binits{R.}},
\bauthor{\bsnm{Chen}, \binits{J.}},
\bauthor{\bsnm{Wei}, \binits{G.-W.}}:
\batitle{Mechanisms of {SARS-CoV}-2 evolution revealing vaccine-resistant
  mutations in {E}urope and {A}merica}.
\bjtitle{The journal of physical chemistry letters}
\bvolume{12}(\bissue{49}),
\bfpage{11850}--\blpage{11857}
(\byear{2021})
\end{barticle}
\endbibitem

\bibitem{wang2020persistent}
\begin{barticle}
\bauthor{\bsnm{Wang}, \binits{R.}},
\bauthor{\bsnm{Nguyen}, \binits{D.D.}},
\bauthor{\bsnm{Wei}, \binits{G.-W.}}:
\batitle{Persistent spectral graph}.
\bjtitle{International Journal for Numerical Methods in Biomedical Engineering}
\bvolume{36}(\bissue{9}),
\bfpage{3376}
(\byear{2020})
\end{barticle}
\endbibitem

\bibitem{wang2021hermes}
\begin{barticle}
\bauthor{\bsnm{Wang}, \binits{R.}},
\bauthor{\bsnm{Zhao}, \binits{R.}},
\bauthor{\bsnm{Ribando-Gros}, \binits{E.}},
\bauthor{\bsnm{Chen}, \binits{J.}},
\bauthor{\bsnm{Tong}, \binits{Y.}},
\bauthor{\bsnm{Wei}, \binits{G.-W.}}:
\batitle{{HERMES}: Persistent spectral graph software}.
\bjtitle{Foundations of Data Science}
\bvolume{3}(\bissue{1}),
\bfpage{67}--\blpage{97}
(\byear{2021})
\end{barticle}
\endbibitem

\bibitem{wang2021computing}
\begin{barticle}
\bauthor{\bsnm{Wang}, \binits{S.}},
\bauthor{\bsnm{Chern}, \binits{A.}}:
\batitle{Computing minimal surfaces with differential forms}.
\bjtitle{ACM Transactions on Graphics (TOG)}
\bvolume{40}(\bissue{4}),
\bfpage{1}--\blpage{14}
(\byear{2021})
\end{barticle}
\endbibitem

\bibitem{wang2025join}
\begin{botherref}
\oauthor{\bsnm{Wang}, \binits{Y.}},
\oauthor{\bsnm{Liu}, \binits{X.}},
\oauthor{\bsnm{Zhang}, \binits{Y.}},
\oauthor{\bsnm{Wang}, \binits{X.}},
\oauthor{\bsnm{Xia}, \binits{K.}}:
Join persistent homology ({JPH})-based machine learning for
  metalloprotein--ligand binding affinity prediction.
Journal of Chemical Information and Modeling
(2025)
\end{botherref}
\endbibitem

\bibitem{warner1983foundations}
\begin{bbook}
\bauthor{\bsnm{Warner}, \binits{F.W.}}:
\bbtitle{Foundations of Differentiable Manifolds and Lie Groups}
vol. \bseriesno{94}.
\bpublisher{Springer},
\blocation{\hspace{0pt}}
(\byear{1983})
\end{bbook}
\endbibitem

\bibitem{wee2022persistent}
\begin{barticle}
\bauthor{\bsnm{Wee}, \binits{J.}},
\bauthor{\bsnm{Xia}, \binits{K.}}:
\batitle{Persistent spectral based ensemble learning ({PerSpect-EL}) for
  protein-protein binding affinity prediction}.
\bjtitle{Briefings in Bioinformatics}
\bvolume{23}(\bissue{2}),
\bfpage{024}
(\byear{2022})
\end{barticle}
\endbibitem

\bibitem{wee2023persistent}
\begin{barticle}
\bauthor{\bsnm{Wee}, \binits{J.}},
\bauthor{\bsnm{Bianconi}, \binits{G.}},
\bauthor{\bsnm{Xia}, \binits{K.}}:
\batitle{Persistent {Dirac} for molecular representation}.
\bjtitle{Scientific Reports}
\bvolume{13}(\bissue{1}),
\bfpage{11183}
(\byear{2023})
\end{barticle}
\endbibitem

\bibitem{wei2025persistentsheaf}
\begin{barticle}
\bauthor{\bsnm{Wei}, \binits{X.}},
\bauthor{\bsnm{Wei}, \binits{G.-W.}}:
\batitle{Persistent sheaf {L}aplacians}.
\bjtitle{Foundations of Data Science}
\bvolume{7}(\bissue{2}),
\bfpage{446}--\blpage{463}
(\byear{2025})
\end{barticle}
\endbibitem

\bibitem{wei2025persistent}
\begin{barticle}
\bauthor{\bsnm{Wei}, \binits{X.}},
\bauthor{\bsnm{Wei}, \binits{G.-W.}}:
\batitle{Persistent topological {Laplacians} -- a survey}.
\bjtitle{Mathematics}
\bvolume{13}(\bissue{2}),
\bfpage{208}
(\byear{2025})
\end{barticle}
\endbibitem

\bibitem{weil1949numbers}
\begin{barticle}
\bauthor{\bsnm{Weil}, \binits{A.}}:
\batitle{Numbers of solutions of equations in finite fields}.
\bjtitle{Bulletin of the American Mathematical Society}
\bvolume{55}(\bissue{5}),
\bfpage{497}--\blpage{508}
(\byear{1949})
\end{barticle}
\endbibitem

\bibitem{wu2023metabolomic}
\begin{barticle}
\bauthor{\bsnm{Wu}, \binits{S.}},
\bauthor{\bsnm{Liu}, \binits{X.}},
\bauthor{\bsnm{Dong}, \binits{A.}},
\bauthor{\bsnm{Gragnoli}, \binits{C.}},
\bauthor{\bsnm{Griffin}, \binits{C.}},
\bauthor{\bsnm{Wu}, \binits{J.}},
\bauthor{\bsnm{Yau}, \binits{S.-T.}},
\bauthor{\bsnm{Wu}, \binits{R.}}:
\batitle{The metabolomic physics of complex diseases}.
\bjtitle{Proceedings of the National Academy of Sciences}
\bvolume{120}(\bissue{42}),
\bfpage{2308496120}
(\byear{2023})
\end{barticle}
\endbibitem

\bibitem{xia2014persistent}
\begin{barticle}
\bauthor{\bsnm{Xia}, \binits{K.}},
\bauthor{\bsnm{Wei}, \binits{G.-W.}}:
\batitle{Persistent homology analysis of protein structure, flexibility, and
  folding}.
\bjtitle{International journal for numerical methods in biomedical engineering}
\bvolume{30}(\bissue{8}),
\bfpage{814}--\blpage{844}
(\byear{2014})
\end{barticle}
\endbibitem

\bibitem{xia2015multidimensional}
\begin{barticle}
\bauthor{\bsnm{Xia}, \binits{K.}},
\bauthor{\bsnm{Wei}, \binits{G.-W.}}:
\batitle{Multidimensional persistence in biomolecular data}.
\bjtitle{Journal of computational chemistry}
\bvolume{36}(\bissue{20}),
\bfpage{1502}--\blpage{1520}
(\byear{2015})
\end{barticle}
\endbibitem

\bibitem{xia2015multiresolution}
\begin{botherref}
\oauthor{\bsnm{Xia}, \binits{K.}},
\oauthor{\bsnm{Zhao}, \binits{Z.}},
\oauthor{\bsnm{Wei}, \binits{G.-W.}}:
Multiresolution persistent homology for excessively large biomolecular
  datasets.
The Journal of chemical physics
\textbf{143}(13)
(2015)
\end{botherref}
\endbibitem

\bibitem{xia2015multiresolution2}
\begin{barticle}
\bauthor{\bsnm{Xia}, \binits{K.}},
\bauthor{\bsnm{Zhao}, \binits{Z.}},
\bauthor{\bsnm{Wei}, \binits{G.-W.}}:
\batitle{Multiresolution topological simplification}.
\bjtitle{Journal of Computational Biology}
\bvolume{22}(\bissue{9}),
\bfpage{887}--\blpage{891}
(\byear{2015})
\end{barticle}
\endbibitem

\bibitem{xia2015persistent}
\begin{barticle}
\bauthor{\bsnm{Xia}, \binits{K.}},
\bauthor{\bsnm{Feng}, \binits{X.}},
\bauthor{\bsnm{Tong}, \binits{Y.}},
\bauthor{\bsnm{Wei}, \binits{G.W.}}:
\batitle{Persistent homology for the quantitative prediction of fullerene
  stability}.
\bjtitle{Journal of computational chemistry}
\bvolume{36}(\bissue{6}),
\bfpage{408}--\blpage{422}
(\byear{2015})
\end{barticle}
\endbibitem

\bibitem{yang2021clebsch}
\begin{barticle}
\bauthor{\bsnm{Yang}, \binits{S.}},
\bauthor{\bsnm{Xiong}, \binits{S.}},
\bauthor{\bsnm{Zhang}, \binits{Y.}},
\bauthor{\bsnm{Feng}, \binits{F.}},
\bauthor{\bsnm{Liu}, \binits{J.}},
\bauthor{\bsnm{Zhu}, \binits{B.}}:
\batitle{Clebsch gauge fluid}.
\bjtitle{ACM Transactions on Graphics (TOG)}
\bvolume{40}(\bissue{4}),
\bfpage{1}--\blpage{11}
(\byear{2021})
\end{barticle}
\endbibitem

\bibitem{yegnesh2016persistence}
\begin{botherref}
\oauthor{\bsnm{Yegnesh}, \binits{K.}}:
Persistence and sheaves.
arXiv preprint arXiv:1612.03522
(2016)
\end{botherref}
\endbibitem

\bibitem{Yin2023FluidCohomology}
\begin{botherref}
\oauthor{\bsnm{Yin}, \binits{H.}},
\oauthor{\bsnm{Nabizadeh}, \binits{M.S.}},
\oauthor{\bsnm{Wu}, \binits{B.}},
\oauthor{\bsnm{Wang}, \binits{S.}},
\oauthor{\bsnm{Chern}, \binits{A.}}:
Fluid cohomology.
ACM Trans. Graph.
\textbf{42}(4)
(2023)
\end{botherref}
\endbibitem

\bibitem{yoon2018cellular}
\begin{botherref}
\oauthor{\bsnm{Yoon}, \binits{H.R.}}:
Cellular sheaves and cosheaves for distributed topological data analysis.
PhD thesis,
University of Pennsylvania
(2018)
\end{botherref}
\endbibitem

\bibitem{zhang2022twisted}
\begin{barticle}
\bauthor{\bsnm{Zhang}, \binits{M.M.}},
\bauthor{\bsnm{Li}, \binits{J.Y.}},
\bauthor{\bsnm{Wu}, \binits{J.}}:
\batitle{The twisted homology of simplicial set}.
\bjtitle{Acta Mathematica Sinica, English Series}
\bvolume{38}(\bissue{10}),
\bfpage{1781}--\blpage{1802}
(\byear{2022})
\end{barticle}
\endbibitem

\bibitem{zhao2018protein}
\begin{barticle}
\bauthor{\bsnm{Zhao}, \binits{R.}},
\bauthor{\bsnm{Cang}, \binits{Z.}},
\bauthor{\bsnm{Tong}, \binits{Y.}},
\bauthor{\bsnm{Wei}, \binits{G.-W.}}:
\batitle{Protein pocket detection via convex hull surface evolution and
  associated reeb graph}.
\bjtitle{Bioinformatics}
\bvolume{34}(\bissue{17}),
\bfpage{830}--\blpage{837}
(\byear{2018})
\end{barticle}
\endbibitem

\bibitem{zhao20193d}
\begin{barticle}
\bauthor{\bsnm{Zhao}, \binits{R.}},
\bauthor{\bsnm{Desbrun}, \binits{M.}},
\bauthor{\bsnm{Wei}, \binits{G.-W.}},
\bauthor{\bsnm{Tong}, \binits{Y.}}:
\batitle{{3D} {Hodge} decompositions of edge- and face-based vector fields}.
\bjtitle{ACM Transactions on Graphics (TOG)}
\bvolume{38}(\bissue{6}),
\bfpage{1}--\blpage{13}
(\byear{2019})
\end{barticle}
\endbibitem

\bibitem{zhao2020rham}
\begin{barticle}
\bauthor{\bsnm{Zhao}, \binits{R.}},
\bauthor{\bsnm{Wang}, \binits{M.}},
\bauthor{\bsnm{Chen}, \binits{J.}},
\bauthor{\bsnm{Tong}, \binits{Y.}},
\bauthor{\bsnm{Wei}, \binits{G.-W.}}:
\batitle{The de {R}ham--{H}odge analysis and modeling of biomolecules}.
\bjtitle{Bulletin of mathematical biology}
\bvolume{82},
\bfpage{1}--\blpage{38}
(\byear{2020})
\end{barticle}
\endbibitem

\bibitem{zheng2024towards}
\begin{barticle}
\bauthor{\bsnm{Zheng}, \binits{J.}},
\bauthor{\bsnm{Feng}, \binits{Z.}},
\bauthor{\bsnm{Ekstrom}, \binits{A.D.}}:
\batitle{Towards analysis of multivariate time series using topological data
  analysis}.
\bjtitle{Mathematics}
\bvolume{12}(\bissue{11}),
\bfpage{1727}
(\byear{2024})
\end{barticle}
\endbibitem

\bibitem{zhou2019pymesh}
\begin{botherref}
\oauthor{\bsnm{Zhou}, \binits{Q.}}:
Pymesh-geometry processing library for {Python}.
Software available for download at https://github.com/PyMesh/PyMesh
\textbf{7}
(2019)
\end{botherref}
\endbibitem

\bibitem{zhou2021mapper}
\begin{bchapter}
\bauthor{\bsnm{Zhou}, \binits{Y.}},
\bauthor{\bsnm{Chalapathi}, \binits{N.}},
\bauthor{\bsnm{Rathore}, \binits{A.}},
\bauthor{\bsnm{Zhao}, \binits{Y.}},
\bauthor{\bsnm{Wang}, \binits{B.}}:
\bctitle{Mapper interactive: a scalable, extendable, and interactive toolbox
  for the visual exploration of high-dimensional data}.
In: \bbtitle{2021 IEEE 14th Pacific Visualization Symposium (PacificVis)},
pp. \bfpage{101}--\blpage{110}
(\byear{2021}).
\bcomment{IEEE}
\end{bchapter}
\endbibitem

\bibitem{zia2025persistent}
\begin{barticle}
\bauthor{\bsnm{Zia}, \binits{M.}},
\bauthor{\bsnm{Jones}, \binits{B.}},
\bauthor{\bsnm{Feng}, \binits{H.}},
\bauthor{\bsnm{Wei}, \binits{G.-W.}}:
\batitle{Persistent directed flag {Laplacian} ({PDFL})-based machine learning
  for protein--ligand binding affinity prediction}.
\bjtitle{Journal of Chemical Theory and Computation}
\bvolume{21}(\bissue{8}),
\bfpage{4276}--\blpage{4285}
(\byear{2025})
\end{barticle}
\endbibitem

\bibitem{zielinski2019persistence}
\begin{bchapter}
\bauthor{\bsnm{Zieli{\'n}ski}, \binits{B.}},
\bauthor{\bsnm{Lipi{\'n}ski}, \binits{M.}},
\bauthor{\bsnm{Juda}, \binits{M.}},
\bauthor{\bsnm{Zeppelzauer}, \binits{M.}},
\bauthor{\bsnm{D{\l}otko}, \binits{P.}}:
\bctitle{Persistence bag-of-words for topological data analysis}.
In: \bbtitle{Proceedings of the 28th International Joint Conference on
  Artificial Intelligence},
pp. \bfpage{4489}--\blpage{4495}
(\byear{2019})
\end{bchapter}
\endbibitem

\bibitem{zielinski2021persistence}
\begin{barticle}
\bauthor{\bsnm{Zieli{\'n}ski}, \binits{B.}},
\bauthor{\bsnm{Lipi{\'n}ski}, \binits{M.}},
\bauthor{\bsnm{Juda}, \binits{M.}},
\bauthor{\bsnm{Zeppelzauer}, \binits{M.}},
\bauthor{\bsnm{D{\l}otko}, \binits{P.}}:
\batitle{Persistence codebooks for topological data analysis}.
\bjtitle{Artificial Intelligence Review}
\bvolume{54},
\bfpage{1969}--\blpage{2009}
(\byear{2021})
\end{barticle}
\endbibitem

\bibitem{zomorodian2010fast}
\begin{barticle}
\bauthor{\bsnm{Zomorodian}, \binits{A.}}:
\batitle{Fast construction of the {Vietoris-Rips} complex}.
\bjtitle{Computers \& Graphics}
\bvolume{34}(\bissue{3}),
\bfpage{263}--\blpage{271}
(\byear{2010})
\end{barticle}
\endbibitem

\bibitem{zomorodian2005computing}
\begin{barticle}
\bauthor{\bsnm{Zomorodian}, \binits{A.}},
\bauthor{\bsnm{Carlsson}, \binits{G.}}:
\batitle{Computing persistent homology}.
\bjtitle{Discrete \& Computational Geometry}
\bvolume{33}(\bissue{2}),
\bfpage{249}--\blpage{274}
(\byear{2005})
\end{barticle}
\endbibitem

\bibitem{zomorodian2008localized}
\begin{barticle}
\bauthor{\bsnm{Zomorodian}, \binits{A.}},
\bauthor{\bsnm{Carlsson}, \binits{G.}}:
\batitle{Localized homology}.
\bjtitle{Computational Geometry}
\bvolume{41}(\bissue{3}),
\bfpage{126}--\blpage{148}
(\byear{2008})
\end{barticle}
\endbibitem

\end{thebibliography}

\end{document}